\documentclass[dvipsnames]{amsart}

\newcounter{puzzles}
\usepackage{pdfsync}\synctex=1
\usepackage{fullpage,subfiles,bibentry,dynkin-diagrams}

\newcommand\plus{{\hspace{-0.1em}{\scalebox{0.7}{$+$}}\hspace{-0em}}}
\newcommand\moins{{\hspace{-0.1em}{\scalebox{0.7}{$-$}}\hspace{-0em}}}
\AtEndDocument{
\bibliographystyle{amsalpha}
\bibliography{multC2.bib}}

\usepackage{hyperref,enumerate}
\usepackage{amsfonts,amssymb,stmaryrd,amsmath}
\usepackage{multirow}

\usepackage{comment}
\usepackage{graphicx,tikz,tabularx,tabulary}
\usetikzlibrary{cd,arrows}\usetikzlibrary{backgrounds}
\usetikzlibrary{matrix,arrows,decorations.markings}
\usepackage{amsthm,mathtools,frcursive,caption}
\usepackage{xcolor}

\usepackage[utf8]{inputenc}
\usepackage[T1]{fontenc}

\definecolor{darkgray}{rgb}{0.33, 0.33, 0.33}

\usepackage[bbgreekl]{mathbbol}
\DeclareSymbolFontAlphabet{\amsmathbb}{AMSb}

\newcommand{\puz}[1]{\ifnum \arabic{puzzles} = 1 \input{#1} \fi}

\newcommand\longto{{\longrightarrow}}

\newcommand\PP{\mathbb P}
\newcommand\CC{\mathbb{C}}\newcommand\QQ{\mathbb{Q}}
\newcommand\NN{\mathbb{N}}\newcommand\ZZ{\mathbb{Z}}

\newcommand{\parenthese}[1]{(#1)}    
\newcommand\scal[1]{\langle #1\rangle}
\newcommand{\drapeau}[2]{{\prescript{}{#1}{X_{#2}}}}

\newcommand{\euler}[1]{\text{\Large $\chi$}_{\hspace{-1pt}\vspace{1pt}\scalebox{0.6}{$#1$}} \hspace{-1pt}}

\makeatletter
\def\revddots{\mathinner{\mkern1mu\raise\p@\vbox{\kern7\p@\hbox{.}}\mkern2mu\raise4\p@\hbox{.}\mkern2mu\raise7\p@\hbox{.}\mkern1mu}}
\makeatother

\newcommand\lr{Littlewood-Richardson }

\newcommand\Pl{P} 
\newcommand{\GB}{{G/B}}
\newcommand{\Fl}{\mathrm{Fl}}
\newcommand{\GP}{{G/P}}
\newcommand\bole{{\leq}}
\newcommand\wbole{{\leq_{\rm w}}}

\newcommand\Sub{{\mathcal {S}}}
\newcommand\ccup{{\cup}}

\newcommand\C{{\mathbb C}}

\newcommand\Z{{\mathbb Z}}
\newcommand\p{{\mathbb P}}

\newcommand\hG{{\hat G}}
\newcommand\hT{{\hat T}}\newcommand\hL{{\hat L}}
\newcommand{\hl}{{\hat l}}

\newcommand\hB{{\hat B}}\newcommand\hP{{\hat P}}

\newcommand\hepsilon{{\hat \epsilon}}
\newcommand\hzeta{{\hat \zeta}}
\newcommand\hW{{\widehat W}}
\newcommand\hPl{{\hat \Pl}}
\newcommand\hg{{\hat g}}

\newcommand\hx{{\hat x}}
\newcommand\hz{{\hat z}}
\newcommand\hWP{{W^\hPl}}
\newcommand\hGB{{\hG/\hB}}\newcommand\hGP{{\hG/\hPl}}
\newcommand\hgamma{{\hat \gamma}}
\newcommand\hpi{{\hat \pi}}
\newcommand\homega{{\hat \omega}}
\newcommand\htheta{{\hat \theta}}

\newcommand\halpha{{\hat \alpha}}   

\newcommand\cI{{\mathcal{I}}}
\newcommand\cT{{\mathcal{T}}}
\newcommand\cL{{\mathcal{L}}}\newcommand\cF{{\mathcal{F}}}
\newcommand\cD{{\mathcal{D}}}

\newcommand\cP{{\mathcal{P}}}
\newcommand\cO{{\mathcal O}}\newcommand\cY{{\mathcal{Y}}}

\newcommand\cX{{\mathcal X}}

\newcommand\cXPv{{\cX^\cP_v}}
\newcommand\cYv{{\cY_v}}
\newcommand\hcYv{{\cY_\hv}}
\newcommand\hcYvG{{\cY_\hv^G}}

\newcommand\cK{{\mathcal{K}}}

\newcommand{\bv}{{\pmb{v}}} \newcommand{\bY}{{\pmb{Y}}}  \newcommand{\bYv}{{\bY(\bbbv)}}
\newcommand{\bcY}{{\pmb{\cY}}} \newcommand{\bcD}{{\pmb{\cD}}}
\newcommand{\bbbv}{{\pmb{\bbv}}}  \newcommand{\bbby}{{\pmb{\bby}}}  

\newcommand{\bcL}{{\pmb{\cL}}}
\newcommand{\bpi}{{\pmb{\pi}}} 
\newcommand\bcYv{{\pmb{\cY_v}}}
\newcommand\bhv{{\pmb{\hv}}}

\newcommand\hv{{\hat v}}  \newcommand\hbv{{\hat \bv}}
  \newcommand\hby{{\hat \by}}
\newcommand\bgamma{{\pmb{\gamma}}}
\newcommand{\ggamma}{\color{ForestGreen} \gamma}
\newcommand{\ghgamma}{\color{ForestGreen} \hgamma}
\newcommand{\gbbgamma}{\color{ForestGreen} \bbgamma}

\newcommand\bepsilon{{\pmb{\epsilon}}}
\newcommand\hbepsilon{{\hat{\bepsilon}}}
\newcommand\by{{\pmb{y}}}
\newcommand\hbcYvG{{\bcY_\hbv^G}}

\newcommand\hbpi{{\hat \bpi}}\newcommand\hbgamma{{\hat \bgamma}}

\newcommand\Ho{{\operatorname H}^0}

\newcommand\oH{{\operatorname H}}

\newcommand\Pic{{\operatorname{Pic}}}\newcommand\Cl{{\operatorname{Cl}}}

\newcommand\Span{{\operatorname{Span}}}
\newcommand\reg{{\operatorname{reg}}}
\newcommand\Tr{{\operatorname{Tr}}}\newcommand\Frac{{\operatorname{Frac}}}
\newcommand\Spec{{\operatorname{Spec}}}
\newcommand\CaCl{{\operatorname{CaCl}}}
\newcommand\Supp{{\operatorname{Supp}}}
\newcommand\Hom{{\operatorname{Hom}}}
\newcommand\codim{{\operatorname{codim}}}
\newcommand\BKdeg{{\operatorname{BK-deg}}}
\newcommand\GL{{\operatorname{GL}}}  \newcommand\PSL{{\operatorname{PSL}}}
\newcommand\SL{{\operatorname{SL}}}  
\newcommand\Spin{{\operatorname{Spin}}}
\newcommand\OG{{\operatorname{OG}}}

\newcommand\Gr{{\operatorname{Gr}}}

\newcommand\bbG{{\amsmathbb G}} \newcommand\bbL{{\amsmathbb L}}
\newcommand\bbT{{\amsmathbb T}} \newcommand\bbB{{\amsmathbb B}}
\newcommand\bbP{{\amsmathbb P}} \newcommand\bbGB{{\amsmathbb{G/B}}} 
\newcommand\bbu{{\mathbb u}}  \newcommand\bbv{{\mathbb v}}
\newcommand\bby{{\mathbb y}}
\newcommand\bbW{{\amsmathbb W}}\newcommand\bbPhi{{\mathbb
    \Phi}}\newcommand\bbDelta{{\mathbb{\Delta}}}

\newcommand\hc{\mbox{\cursive h}}

\newtheorem{notation}{Notation}
\newtheorem{exple}{Example}
\newtheorem{defi}{Definition}
\newtheorem{rema}{Remark}

\newtheorem{prop}{Proposition}

\newtheorem{theo}[prop]{Theorem}
\newtheorem{lemma}[prop]{Lemma}
\newtheorem{coro}[prop]{Corollary}

\newenvironment{proof*}[1]
  {\renewcommand{\proofname}{#1}%
   \begin{proof}}
  {\end{proof}}

\newcommand\inv{{^{-1}}}
\newcommand\invv{^{-1}}

\newcommand{\tikzline}[3]{
\begin{tikzcd}
    #1 \arrow[d,no head,"#3"] \\
    #2
\end{tikzcd}
}
\newcommand{\tikzlineleft}[3]{
\begin{tikzcd}
    #1 \arrow[d,no head,swap,"#3"] \\
    #2
\end{tikzcd}
}
\newcommand{\tikzdashedline}[3]{
\begin{tikzcd}
    #1 \arrow[d,no head,dashed,"#3"] \\
    #2
\end{tikzcd}
}
\newcommand{\tikzdashedlineleft}[3]{
\begin{tikzcd}
    #1 \arrow[d,no head,dashed,swap,"#3"] \\
    #2
\end{tikzcd}
}
\newcommand{\tikzdashedlineright}[3]{
\begin{tikzcd}
    #1 \arrow[d,no head,dashed,"#3"] \\
    #2
\end{tikzcd}
}
\newcommand{\tikzv}[5]{  
\begin{tikzcd}[row sep=small, column sep=tiny, ampersand replacement=\&]
#1 \arrow[rd,dashed,no head,"#4" left] \&\& #2 \arrow[ld,no head,"#5"] \\
\& #3
\end{tikzcd}
}
\newcommand{\tikzhat}[5]{  
\begin{tikzcd}[row sep=small, column sep=tiny, ampersand replacement=\&]
\& #1 \arrow[ld,no head,"#4" left] \arrow[rd,dashed,no head,"#5"] \\
#2 \&\& #3
\end{tikzcd}
}

\newcommand\DynkinNodeSize{2mm}
\newcommand\DynkinArrowLength{3mm}
\tikzset{
  dnode/.style={
    circle,
    inner sep=0pt,
    minimum size=\DynkinNodeSize,
    fill=white,
    draw},
  middlearrow/.style={
    decoration={markings,
      mark=at position 0.6 with
      {\draw (0:0mm) -- +(+135:\DynkinArrowLength); \draw (0:0mm) -- +(-135:\DynkinArrowLength);},
    },
    postaction={decorate}
  },
  leftrightarrow/.style={
    decoration={markings,
      mark=at position 0.999 with
      {
      \draw (0:0mm) -- +(+135:\DynkinArrowLength); \draw (0:0mm) -- +(-135:\DynkinArrowLength);
      },
      mark=at position 0.001 with
      {
      \draw (0:0mm) -- +(+45:\DynkinArrowLength); \draw (0:0mm) -- +(-45:\DynkinArrowLength);
      },
    },
    postaction={decorate}
  },
  sedge/.style={
  },
  dedge/.style={
    middlearrow,
    double distance=0.5mm,
  },
  tedge/.style={
    middlearrow,
    double distance=1.0mm+\pgflinewidth,
    postaction={draw}, 
  },
  infedge/.style={
    leftrightarrow,
    double distance=0.5mm,
  }
}

\tikzset{
    master/.style={
        execute at end picture={
            \coordinate (lower right) at (current bounding box.south east);
            \coordinate (upper left) at (current bounding box.north west);
        }
    },
    slave/.style={
        execute at end picture={
            \pgfresetboundingbox
            \path (upper left) rectangle (lower right);
        }
    }
}

\newcommand\TikzLambdaLG{
\begin{tikzpicture}[scale=0.7,master]
      \foreach \y in {1,...,3} {\foreach \x in {1,...,\y}
{\pgfmathtruncatemacro{\XnameA}{\x-1} 
\pgfmathtruncatemacro{\YnameA}{\y} 
        \node[dnode] (o\XnameA\YnameA) at
        (\x-\y,5-\y-\x+1) {};
      }}; 

\foreach \y in {0,...,3}
{\pgfmathtruncatemacro{\YnameA}{\y} 
\node[dnode,fill=red] (o\YnameA \YnameA) at
         (1,5-2*\y) {};
      };

\foreach \y in {0,1,2} {\foreach \x in {0,...,\y}
{\pgfmathtruncatemacro{\XnameA}{\x} 
\pgfmathtruncatemacro{\YnameA}{\y} 
\pgfmathtruncatemacro{\YnameB}{1+\y} 
        \draw (o\XnameA\YnameA) -- (o\XnameA\YnameB);
      }}; 
\foreach \y in {1,...,3} {\foreach \x in {1,...,\y}
{\pgfmathtruncatemacro{\XnameA}{-1+\x} 
\pgfmathtruncatemacro{\YnameA}{\y} 
\pgfmathtruncatemacro{\XnameB}{\x} 
        \draw(o\XnameA\YnameA) -- (o\XnameB\YnameA);
      }};

    \node[dnode] (1) at (-2,-2) {};
     \node[dnode] (2) at (-1,-2) {};
   \node[dnode] (3) at (0,-2) {};
  \node[dnode,fill=black] (4) at (1,-2) {};

    \path (1) edge[sedge] (2);
 \path (2) edge[sedge] (3); 
 \path (4) edge[dedge] (3);

\def\shift{5}
 \foreach \y in {1,2,...,3} {\foreach \x in {1,...,\y}
{\pgfmathtruncatemacro{\XnameA}{\x} 
\pgfmathtruncatemacro{\YnameA}{\y} 
        \node[dnode] (\XnameA\YnameA) at
        (\shift+1+\x-\y,6-\y-\x) {};
}
      }; 
\foreach \y in {1,2,...,4}{\node[dnode,fill=red] (\y) at
        (\shift+2,\y-3) {};
};
\foreach \y in {1,2,...,3}{
\pgfmathtruncatemacro{\YnameB}{1+\y};
\draw (\y) -- (\YnameB);
};

\foreach \y in {1,2} {\foreach \x in {1,...,\y}
{\pgfmathtruncatemacro{\XnameA}{\x} 
\pgfmathtruncatemacro{\YnameA}{\y} 
\pgfmathtruncatemacro{\YnameB}{1+\y} 
\draw (\XnameA \YnameA) -- (\XnameA \YnameB);
}};
\foreach \y in {2,3} {\foreach \x in {2,...,\y}
{\pgfmathtruncatemacro{\XnameA}{\x} 
\pgfmathtruncatemacro{\YnameA}{\y} 
\pgfmathtruncatemacro{\XnameB}{-1+\x} 
\draw(\XnameA \YnameA) -- (\XnameB\YnameA);
}};

\draw[green]  (2) -- (33);
\draw[dotted] (o01) -- (11);
    \end{tikzpicture}
}

\newlength\mylength
\begin{document}

\title{Reduction for branching multiplicities}

\author[P.-E. Chaput]{Pierre-Emmanuel Chaput}
\address{Universit\'e de Lorraine, CNRS, Institut \'Elie Cartan de Lorraine, UMR 7502, Vandoeu\-vre-l\`es-Nancy, F-54506, France}
\email{pierre-emmanuel.chaput@univ-lorraine.fr}

\author[N. Ressayre]{Nicolas Ressayre}
\address{Université Claude Bernard Lyon I, Institut Camille Jordan (ICJ), UMR CNRS 5208, 43 boulevard du 11 novembre 1918,
69622 Villeurbanne CEDEX}
\email{ressayre@math.univ-lyon1.fr}

\thanks{This work was supported by the ANR GeoLie project, of the French Agence Nationale de la Recherche.}

\begin{abstract}
A reduction formula for the branching coefficients of 
the restrictions
of representations of a semisimple group to a semisimple subgroup is proved in \cite{kt,dw, ressayre-birational}.
This formula holds when the highest weights of the representations belong to a codimension $1$ face of the Horn cone, which by
\cite{GITEigen} corresponds to some Schubert coefficient equal to $1$. We prove a similar reduction formula when this
Schubert coefficient is equal to $2$, and show some properties of the class of the branch divisor corresponding to a generically finite
morphism naturally defined in this context.
\end{abstract}

\maketitle

\section{Introduction}

Fix an inclusion $G \subset \hG$ of complex connected reductive groups. 
We are interested in the branching problem for the decomposition of irreducible
$\hG$-modules as representations of $G$.
The appearing multiplicities are nonnegative integers parametrized by
the pairs of dominant weights for $G$ and $\hG$ respectively. 
The support of this multiplicity function is known to be a finitely
generated semigroup.
It generates a convex polyhedral cone $\Gamma_\QQ(G\subset\hG)$, that
we call the {\it Horn cone}.
The set of codimension one faces of this cone is in bijection with the
set of Levi-movable pairs (see Section~\ref{sec:BKprod}) of Schubert
classes in some projective homogeneous spaces $\GP \subset \hGP$ such
that the Schubert coefficient $c$ of the first
Schubert class with the pullback of the second Schubert class in $\GP$ is equal to one.
See \cite{GITEigen}.

Moreover, the
multiplicities on those faces satisfy reduction rules: if $L,\hL$ denote
the Levi subgroups of $P,\hP$, the branching multiplicities for the inclusion
$G \subset \hG$ on such a face are equal to branching multiplicities
for the inclusion $L \subset \hL$. See \cite{ressayre-birational}.

Actually, to any pair of Schubert classes such that the Schubert coefficient $c$ is positive corresponds an
inequality satisfied by $\Gamma_\QQ(G\subset\hG)$. Our aim is to study the associated reduction
rules for the multiplicities on the associated face when $c=2$.
Although most of our results are general we consider in this introduction the
case of the tensor product decomposition for the linear group.

\medskip

Fix an $n$-dimensional vector space $V$.
Let $\Lambda_n^+=\{(\lambda_1\geq\cdots\geq\lambda_n\geq
0\,:\,\lambda_i\in\NN\}$ denote the set of partitions. 
For $\lambda\in\Lambda_n^+$, let $S^\lambda V$ be the corresponding Schur module,
that is the irreducible $\GL(V)$-module of highest weight $\sum \lambda_i \epsilon_i$ (notation as in \cite{bourb}). 
The Littlewood-Richardson coefficients (or LR coefficients for short) $c_{\lambda,\mu}^\nu$ are
defined by
\begin{equation}
  \label{eq:17}
  S^\lambda V\otimes S^\mu V \simeq \bigoplus_{\nu\in\Lambda_n^+}
  \C^{c_{\lambda,\mu}^\nu} \otimes S^\nu V,
\end{equation}
where $\C^{c_{\lambda,\mu}^\nu}$ is a multiplicity space.

Let $1\leq r\leq n-1$ and $\Gr(r,n)$ be the Grassmannian of $r$-dimensional
linear subspaces of $V$. Recall that the Schubert basis $(\sigma^I)$ of the Chow
ring $A^*(\Gr(r,n),\ZZ)$ is parametrized by the subsets $I$ of
$\{1,\dots,n\}$ with $r$ elements. 
The Schubert coefficients $c_{I,J}^K$ are defined by
\begin{equation}
  \label{eq:15}
  \sigma^I \sigma^J=\sum_K c_{I,J}^K \sigma^K.
\end{equation}
Actually, $c_{I,J}^K$ is also a LR coefficient by \cite{Lesieur}, but
this coincidence is specific to the type $A$.

Given a partition $\lambda$ and a subset $I$, let $\lambda_I$ be the
partition whose 
parts are $\lambda_i$ with $i \in I$. Let also $\overline I$
denote the complementary subset $\{1,\dots,n\} \setminus I$. 

\begin{theo} \cite{kt,dw} 
\label{introtheo:Horn}
\begin{enumerate}
\item Let $1\leq r\leq n-1$ and $I,J,K \subset \{1,\ldots,n\}$ be subsets with $r$ elements such that $c_{I,J}^K\neq 0$.
 For any $(\lambda,\mu,\nu)\in(\Lambda_n^+)^3$, if 
 $c_{\lambda,\mu}^\nu\neq 0$ then
 \begin{equation}
   \label{eq:166}
   |\lambda_I|+|\mu_J|\geq |\nu_K|\,.
 \end{equation}
\item Conversely, fix $(\lambda,\mu,\nu)\in(\Lambda_n^+)^3$ such that
  $|\lambda|+|\mu|=|\nu|$. 
If for any $1\leq r\leq n-1$ and $I,J,K$ of cardinal $r$ such that $c_{I,J}^K=1$,
inequality~\eqref{eq:166} holds, then $c_{\lambda,\mu}^\nu\neq 0$.
\end{enumerate}
\end{theo}

The reduction result that we mentioned above, corresponding
to the semisimple part $\SL_r \times \SL_{n-r} \subset \SL_n$ of the Levi subgroup is the following:

\begin{theo} \cite{KTT:factorLR,Roth,dw,ressayre-birational} 
\label{introtheo:birational}
Assume that $c_{I,J}^K=1$. Let $\lambda,\mu,\nu$ be partitions such that
\begin{equation}
\label{equ:intro}
|\lambda_I|+|\mu_J|=|\nu_K|\,.
\end{equation}
Then
\begin{equation}
c_{\lambda,\mu}^\nu = \  
c_{\lambda_I,\mu_J}^{\nu_K} \cdot c_{\lambda_{\overline
    I},\mu_{\overline J}}^{\nu_{\overline K}}\, .
\label{eq:DWred}
\end{equation}
\end{theo}

Formula \eqref{eq:DWred} is a multiplicativity property. 
Let us first report on a similar property for  Belkale-Kumar
\cite{bk} coefficients (BK coefficients for short).
See Section~\ref{sec:BKprod} for a short description of these numbers
which are all Schubert coefficients or zero.
Consider an inclusion $P \subset Q$ of parabolic subgroups of a reductive algebraic group $G$,
and the corresponding fibration $G/P \to G/Q$.
Richmond \cite{richmond} proved that any  BK coefficient $d$ of $G/P$
is the product of two  BK coefficients for $G/Q$ and $Q/P$.
In type A, this implies that a non-zero BK coefficient for any two steps flag
manifold is a product of two LR coefficients: $d=c_1c_2$.

If moreover $c_1=1$, Theorem~\ref{introtheo:birational} implies that
$c_2$ itself is the product of two LR coefficients:
$c_2=c_2'c_2''$. Thus $d=c_1c_2'c_2''$ is the product of \emph{three}
LR coefficients.
This is the content of \cite[Theorem 3]{kp}, which even more generally states that on a
$k$-step flag variety, any BK coefficient can be factorized as
a product of $\frac{k(k-1)}{2}$ LR coefficients. Unfortunately, this assertion
needs $c_1=1$ and is not correct in general, as we show in Remark~\ref{rema-erreur}.
Our original motivation was to correct this result. We get such a
correction if $c_1=2$.

\medskip

In \cite{ressayre-birational},
the variety $Y \subset \Gr(r,n) \times (G/B)^3$ is defined, where
a quadruple $(V,X_1,X_2,X_3)$ belongs to $Y$ if and only if $V \in \Gr(r,n)$
belongs to the intersection
of the three Schubert varieties defined by $I,J,K$ and the three flags $X_1,X_2,X_3$.

The hypothesis $c_{I,J}^K=1$ implies that the projection $\Pi:Y \longto (G/B)^3$ is birational.
By the projection formula, since $Y$ has rational singularities, it follows that
the sections of line bundles on $(G/B)^3$ identify to those on $Y$. Let $C \subset (G/B)^3$ be a
product of flag varieties under $L$. In \cite{ressayre-birational}, it is shown that
taking $G$-invariants
and restricting to $C$ is the same as restricting to $C$ and
then taking $L$-invariants, which leads to \eqref{eq:DWred}.

\bigskip

Fix now $r$ and $I,J,K \subset \{1,\ldots,n\}$ of cardinal $r$ such that
$c_{I,J}^K=2$. We can still define $\Pi:Y \longto (G/B)^3$ in a similar way,
but it is now generically finite of degree $2$.
We show in Subsection~\ref{sec:genfinite} that we can associate to
this generically finite morphism a branch divisor
in $(G/B)^3$ and that the pushforward of the structure sheaf $\cO_Y$
is expressed in terms of this branch divisor.
In Theorem~\ref{theo:branch}~, we determine the triple
$(\alpha,\beta,\gamma)$ of dominant weights defining the branch divisor
class.
This is done using the morphism $Y \to \Gr(r,n)$ which describes $Y$ as a relative product
of Schubert varieties over
the Grassmannian. This shows that $Y$ is normal and allows computing the
ramification divisor as the relative canonical
sheaf, adapting in this relative setting previous computations of the canonical sheaf
of Schubert varieties
(\cite[Theorem 4.2]{Raman} and \cite[Proposition 4.4]{perrin-small}),
see Section~\ref{sub:canonical}.

\medskip

Under  an assumption of Levi-movability  (see Sections~\ref{sec:BKprod}
and~\ref{sec:eigenconeBPi}), which is automatically satisfied in
type~A, we get an alternative characterization of
this triple in Proposition~\ref{pro:posC}. 
Namely, $(\alpha,\beta,\gamma)$ is the unique minimal element of the
set of triples
$(\lambda,\mu,\nu)\in(\Lambda_n^+)^3$ such that 
$$
  |\lambda_I|+|\mu_J|=|\nu_K|\ \mbox{ and }\ 
0\neq c_{\lambda,\mu}^\nu < \ c_{\lambda_I,\mu_J}^{\nu_K} \cdot
c_{\lambda_{\overline I},\mu_{\overline J}}^{\nu_{\overline K}}.
$$

The interested reader can find at \cite{programmation} a program allowing to compute
$(\alpha,\beta,\gamma)$ and its generalisations to any classical group
(and under the more general assumption $c_{I,J}^K\neq 0$).

Arguing as in \cite{ressayre-birational}, we deduce from our description of $\Pi_* \cO_Y$
our general result, Theorem \ref{main-theo}. It holds for any inclusion $G \subset \hG$.
For $\hG= G \times G$ and $G/P$ a Grassmannian, it states the following modification of
\eqref{eq:DWred}:

\begin{theo}
\label{main-theo-A}
Assume  $c_{I,J}^K=2$. Let $\lambda,\mu,\nu$ be partitions such that
$|\lambda_I|+|\mu_J|=|\nu_K|$.
Then
\begin{equation}
\label{main-formula}
c_{\lambda,\mu}^\nu + c_{\lambda-\alpha,\mu-\beta}^{\nu-\gamma} =   
c_{\lambda_I,\mu_J}^{\nu_K} \cdot
c_{\lambda_{\overline I},\mu_{\overline J}}^{\nu_{\overline K}}\, .
\end{equation}
\end{theo} 

Observe that although the Horn cone is defined by inequations of the form 
\eqref{equ:intro} for $I,J,K$ such that $c_{I,J}^K=1$ inside $(\Lambda_n^+)^3$, the inequalities coming from
the fact $\lambda,\mu,\nu \in \Lambda^+_n$ make it possible that an element $(\lambda,\mu,\nu)$
in the Horn cone satisfies an equality
$|\lambda_I|+|\mu_J|=|\nu_K|$ for some triple $(I,J,K)$ such that $c_{I,J}^K=2$ but for no triple such
that $c_{I,J}^K=1$: see Remark \ref{rem:facevsfacereg}. For such an element, Theorem \ref{introtheo:Horn}
cannot be applied but Theorem \ref{main-theo-A} can.

\medskip

Under the assumption of Levi-movability, the mentioned
characterization of $(\alpha,\beta,\gamma)$ implies that it
satisfies~\eqref{equ:intro}. 
Then, an immediate induction expresses the multiplicity  as an alternating sum of similar
coefficients for the reduction $L
\subset \hL$. This is Corollary~\ref{cor:alt} below.
In the setting of Theorem~\ref{main-theo-A}, we get: 

\begin{theo} 
\label{introtheo:LRalt}
Assume  $c_{I,J}^K=2$. Let $\lambda,\mu,\nu$ be partitions such that
$|\lambda_I|+|\mu_J|=|\nu_K|$.
Then,
\begin{equation}
c_{\lambda,\mu}^\nu = \sum_{k\geq 0}(-1)^k   
c_{\lambda_I-k\alpha_I,\mu_J-k\beta_J}^{\nu_K-k\gamma_K} \cdot c_{\lambda_{\overline
    I}-k\alpha_{\overline I},\mu_{\overline J}-k\beta_{\overline J}}^{\nu_{\overline K}-k\gamma_{\overline K}}\, .
\label{eq:LRalt}
\end{equation}
\end{theo}

\bigskip

Given three subsets $I,J,K$ as above, it is a difficult task to describe the face it defines in the
Horn cone: to the best of our knowledge, even the dimension of such a face is not known in general,
and our experience is that a computer will only give very limited information related to this problem. 
From a theoretical point of view as well as a computational point of view, both finding linear equations defining this face
and points on it is challenging. In Section~\ref{sub:big}, we describe completely the faces corresponding to some relevant examples.

We also believe that the study of the ramification or branch divisor of $Y \longto (G/B)^3$ might be
interesting in itself. In general, we have a quite poor understanding of the geometry of these divisors,
although we know their classes as Weil divisors, by Theorem \ref{theo:branch}. In the example treated in
Section \ref{sub:g36}, we have a very explicit and nice description of the ramification and
branch divisors in terms of
projective geometry in $\PP^2$. Moreover, we observe that they are irreducible.

Finally, we give in Section \ref{expleB} examples in type $B$ that show that the Levi-movability
assumption in Theorem \ref{introtheo:LRalt} is necessary.

\bigskip

Our methods do not extend easily to the case when the multiplicity
$c_{I,J}^K$ is more than $2$. Indeed, a key point in our arguments is
that any finite degree 2 morphism is cyclic, allowing to compute the
pushforward of the structural sheaf.
Example \ref{exam:c=3} shows that the situation is deeply different when $c_{I,J}^K >2$.
Maybe in the case of multiplicity $3$, \cite{miranda} could be of some help. It might also be
helpful considering a deformation of $\Pi$ which gives a cyclic covering. Another natural possibility would be to
study directly the combinatorics as done in \cite{kp} with the new understanding we have of the geometry involved.

\bigskip \noindent {\bf Acknowledgement.}
The authors are partially supported by the French National Agency
(Project GeoLie ANR-15-CE40-0012). 

\tableofcontents

\section{Main results}
\label{sec:main}

\subsection{Notation}
\label{sec:notations}

In this section, we introduce most of our needed notation.

\subsubsection{Algebraic groups}
Come back to a general inclusion $G \subset \hG$,
or even more generally a finite morphism $G \to \hG$ of connected reductive algebraic
groups. We also assume 
that $G$ and $\hG$ are simply connected.

Fix a maximal torus $T$  and a Borel subgroup $B$ such that $T\subset B$. 
Let $\hB$ be a Borel subgroup of $\hG$ containing $B$ and let $\hT$ be
a maximal torus of $\hB$ containing $T$. 
We denote by $X(T)$ and $X(\hT)$ the groups of characters of $T$ and
$\hT$ respectively, and
by $X(T)^+$ and $X(\hT)^+$ the subsets of dominant weights.

Given two dominant weights $\zeta$ and $\hzeta$ for $G$ and $\hG$
respectively,  we denote by $V_\zeta$ and $V_\hzeta$ the 
highest weight modules of $G$ and $\hG$ respectively.
In this paper, the following integers are what we are mostly interested in:
$$
m_{G \subset \hG}(\zeta,\hzeta) = \dim(V_\zeta \otimes V_\hzeta)^G\,,
$$
where the exponent $G$ means subset of $G$-invariant vectors. 
Indeed, they encode the branching law since
$$
V_\hzeta=\oplus_{\zeta\in X(T)^+} \C^{m_{G \subset \hG}(\zeta,\hzeta)} \otimes V_\zeta^*\,,
$$
as $G$-modules, where $V_\zeta^*$ denotes the dual $G$-module of
$V_\zeta$.\\

We choose a one-parameter subgroup $\tau$ of $T$, which is also a one-parameter subgroup of $\hT$, and thus defines
parabolic subgroups $\Pl=P(\tau) \subset G$ and $\hPl=\hP(\tau) \subset \hG$.
We thus have the following inclusions:
\begin{equation}
\label{equ:inclusions}
\begin{tikzcd}
     T \arrow[hookrightarrow]{r} \arrow[hookrightarrow]{d} & 
     B \arrow[hookrightarrow]{r} \arrow[hookrightarrow]{d} &
     \Pl \arrow[hookrightarrow]{r} \arrow[hookrightarrow]{d} &
     G \arrow[hookrightarrow]{d} \\
     \hT \arrow[hookrightarrow]{r} & \hB \arrow[hookrightarrow]{r} & \hPl \arrow[hookrightarrow]{r} & \hG
\end{tikzcd}
\end{equation}
Denote by $\iota:G/\Pl \longto \hG/\hPl$ the inclusion morphism.
We denote by $L$ the Levi factor of $\Pl$ containing $T$, and similarly for $\hL$.
We denote by $S$ the neutral component of the center of $L$.
Given a character $\gamma$ of $T$ or $\hT$, we denote by $\gamma_{|S}$
its restriction to the torus $S$.

Denote by $W$ and $W_P$ the Weyl groups of $G$ and $\Pl$ respectively,
and by $W^\Pl$ the set of minimal length representatives of the cosets
in $W/W_\Pl$, and similarly $\hW$, $\hWP$ and $W_\hP$.

Consider the group $\bbG=G\times \hG$ and use bbolded letters to refer to
this group. Namely
$\bbT=T\times \hT$, $\bbB=B\times \hB$, $\bbP=P\times\hP$,
$\bbL=L\times\hL$ and $\bbW=W\times\hW$.  
The set $\bbPhi$ of roots of $\bbG$ is 
$$
\bbPhi=(\Phi\times\{0\})\cup (\{0\}\times\hat\Phi),
$$ 
where $\Phi$ and $\hat\Phi$ denote the sets of roots for $G$ and
$\hG$ respectively.

Let $\Delta$, $\hat\Delta$ and $\bbDelta$ be the sets of simple roots for
$G$, $\hat G$ and $\bbG$ respectively.
For $\bbalpha\in\bbDelta$, we denote by $\varpi_\bbalpha$ the
associated fundamental weight. Let $\rho_\bbG$ (resp. $\rho_\bbL$) denote the half sum of positive roots of
$\bbG$ (resp. $\bbL$). Set also $\rho^\bbL=\rho_\bbG-\rho_\bbL$.
Similarly, we use $\omega_\alpha$, $\rho_\hG\dots$

\subsubsection{Schubert classes}
\label{sec:schubert_classes}
Let $p:\GP \to \Spec(\C)$ be the structure morphism, and denote, for $\xi$ in the Chow group $A^*(\GP)$,
\begin{equation}
\label{equ:chi}
\euler{\GP}(\xi)=p_* \xi \in A^0(\Spec(\C)) \simeq \Z\,.
\end{equation}
Since $\GP$ is rationally connected, $A^{\dim(\GP)}=\Z[pt]$ and $\euler{\GP}(\xi)$
is the coefficient of $[pt]$ in the homogeneous part
of degree $\dim(\GP)$ of $\xi$.

\medskip

For $v \in W^P$, $\tau_v$ denotes the
class of the Schubert variety $\overline{BvP/P}$ in $A^{\dim(G/P)-\ell(v)}(G/P)$.
Let $w_0$ be the longest element of $W$
and $w_{0,P}$ be the longest element of $W_P$. Poincaré duality takes the nice form
\begin{equation}
\label{equ:poincare}
\euler{\GP}(\tau_v \cdot \tau_w) = \delta_{v^\vee,w} \mbox{ with } v^\vee:=w_0vw_{0,P}\,.
\end{equation}
Similarly, $(\tau_\hv)_{\hv \in W^\hP}$ denotes the Schubert basis of $A^*(\hG/\hP)$.
When $P=B$ and $u \in W$, we denote by $\sigma_u \in A^*(\GB)$ the Chow class of $\overline{BuB/B} \subset \GB$.
Moreover, we define $\sigma^u := \sigma_{w_0u}$ and $\tau^v := \tau_{w_0vw_{0,P}}$.

Let 
$$
\begin{array}{cccl}
  \delta\,:&\GP&\longto&\bbG/\bbP\\
&gP/P&\longmapsto&(gP/P,g\hP/\hP)
\end{array}
$$
denote the small diagonal map and consider the pullback
$\delta^*\,:\,A^*(\bbG/\bbP)\longto A^*(\GP)$.
For $\bbv=(v,\hv)\in W^\bbP$, observe that
\begin{equation}
\label{equ:delta*}
\delta^*(\tau_{\bbv})=\tau_{v} \cdot \iota^*(\tau_{\hv})\,.
\end{equation}

For $\bbv=(v,\hv)\in W^\bbP=W^P\times W^\hP$ we set
\begin{equation}
\label{eq:def_c}
c(\bbv)=c(v,\hv):=\euler{\GP}(\delta^*(\tau_{\bbv}))\,, \mbox{ so that } \iota^*(\tau_\hv)=\sum_{v\in W^P}c(v,\hv)\tau^{v}\,.
\end{equation}
Since $\tau_{\bbv}$ has degree $\dim \GP + \dim \hGP - \ell(\bbv)$, we have
$c(\bbv)=0$ unless $\ell(\bbv)=\dim(\hG/\hP)(=\ell(v)+\ell(\hv))$.

\subsubsection{Bruhat orders}
A covering relation in the left weak Bruhat order is a pair $(v,w)$ such that there exists a simple root $\alpha$ with
$v = s_\alpha w$ and $\ell(v)=\ell(w)+1$. We use the graph
\tikzlineleft{v}{w}{\alpha}
to depict such a relation.

\begin{figure}
  \begin{center}
    \begin{tikzpicture}[scale=1.5,el/.style = {inner sep=2pt, align=left, sloped}]
      \node[dnode,fill=black] (a) [label=right:{$e$}] at (1,0) {}; 
      \node[dnode,fill=black] (b) [label=left:{$s_1$}] at (0,1) {}; 
      \node[dnode,fill=black] (c) [label=right:{$s_2$}] at (2,1) {}; 
      \node[dnode,fill=black] (d) [label=left:{$s_2s_1$}] at (0,2) {};
      \node[dnode,fill=black] (f) [label=right:{$s_1s_2s_1=s_2s_1s_2$}] at (1,3) {}; 
      \node[dnode,fill=black] (e) [label=right:{$s_1s_2$}] at (2,2)
      {};
       \draw (a) -- node[left] {$\alpha_1$} (b);
\draw (b) -- node[left] {$\alpha_2$} (d);
\draw (d) -- node[left] {$\alpha_1$} (f);
\draw (a) -- node[right] {$\alpha_2$} (c);
\draw (c) -- node[right] {$\alpha_1$} (e);
\draw (e) -- node[right] {$\alpha_2$} (f);
\draw[dashed] (b) -- node[el, below left] {{\tiny $\alpha_1+\alpha_2$}} (e);
\draw[dashed] (c) -- node[el,above left] {{\tiny $\alpha_1+\alpha_2$}}
(d);

\def\shift{6} 

\node[dnode,fill=black] (aa) [label=right:{$e$}] at (\shift+1,0) {}; 
      \node[dnode,fill=black] (bb) [label=left:{$s_1$}] at (\shift,1) {}; 
      \node[dnode,fill=black] (cc) [label=right:{$s_2$}] at (\shift+2,1) {}; 
      \node[dnode,fill=black] (dd) [label=left:{$s_2s_1$}] at (\shift+0,2) {};
      \node[dnode,fill=black] (ff) [label=right:{$s_1s_2s_1=s_2s_1s_2$}] at (\shift+1,3) {}; 
      \node[dnode,fill=black] (ee) [label=right:{$s_1s_2$}] at (\shift+2,2)
      {};
       \draw (aa) -- node[left] {\color{ForestGreen} $\alpha_1$} (bb);
\draw (bb) -- node[left] {\color{ForestGreen} $\alpha_1+\alpha_2$} (dd);
\draw (dd) -- node[left] {\color{ForestGreen} $\alpha_2$} (ff);
\draw (aa) -- node[right] {\color{ForestGreen} $\alpha_2$} (cc);
\draw (cc) -- node[right] {\color{ForestGreen} $\alpha_1+\alpha_2$} (ee);
\draw (ee) -- node[right] {\color{ForestGreen} $\alpha_1$} (ff);
\draw[dashed] (bb) -- node[el, below left] {{\tiny \color{ForestGreen} $\alpha_2$}} (ee);
\draw[dashed] (cc) -- node[el,above left] {{\tiny \color{ForestGreen} $\alpha_1$}}
(dd);

    \end{tikzpicture}
  \end{center}

\caption{Bruhat graphs for $\SL_3$}
\label{fig:Bruhatgraph}
\end{figure}
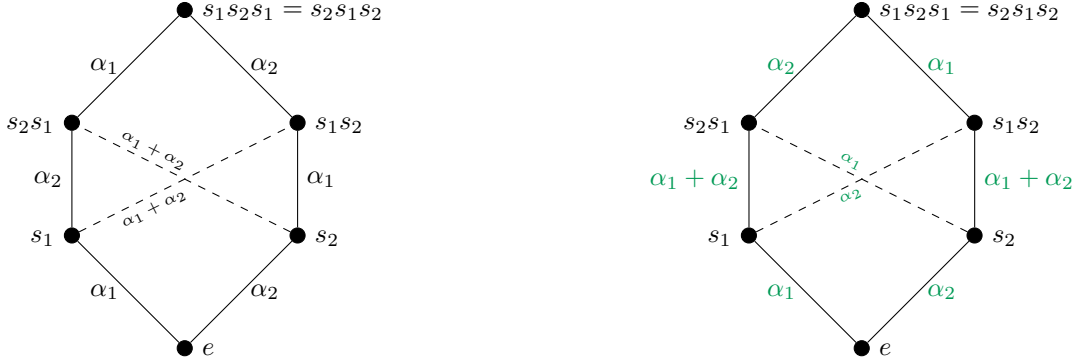

A covering relation in the strong Bruhat order is a pair $(v,w)$ such that there exists a positive root $\beta$ with
$v=s_\beta w$ and $\ell(v)=\ell(w)+1$. We  use the graph
\tikzdashedlineleft{v}{w}{\beta} to depict such a relation.
For example the graph on the left on
Figure~\ref{fig:Bruhatgraph}  is the Bruhat graph of $\SL_3$.

We also define the (weak) twisted Bruhat graph by labelling the
previous edges
by $w\inv\alpha$ and $w\inv \beta$ in place of $\alpha$ and $\beta$. 
We use  green color for these twisted labels to avoid any confusion. 
Hence    
\tikzdashedlineleft{v}{w}{\beta}
means $v=s_\beta w$ and    
\tikzdashedlineright{v}{w}{\ggamma}
means $v=ws_\gamma$, with in both cases $\ell(v)=\ell(w)+1$.
For $\SL_3(\CC)$ we get the right graph of Figure~\ref{fig:Bruhatgraph}.

We denote by $v\bole w$ the strong Bruhat order and by $v\wbole w$ the left weak Bruhat order.




\subsubsection{Chevalley formula}

Given $\zeta\in X(T)$, we denote by $\cL_\GB(\zeta)$ the $G$-linearized line
bundle on $\GB$ such that $B$ acts with weight $-\zeta$ on the fiber
over $B/B$. Similarly, define $\cL_\hGB(\hzeta)$. 
For any $v$ in $W^P$ and any character $\zeta$ of $P$, we have the Chevalley formula, proved for example in
\cite{Prag:mult}:
\begin{equation}
\label{equa:chevalley}
c_1(\cL_\GP(\zeta)) \cup \tau_v = \sum
\scal{\zeta,\gamma^\vee} \tau_{v'}\,,
\end{equation}
where the sum runs over the
covering relations
\tikzdashedline{v}{v'}{\ggamma}
in the Bruhat graph of $G/P$.
Here, $c_1$ denotes the first Chern class.

\subsection{Our reduction result for the branching problem}
\label{sub:reduction}

In \cite{ressayre-birational}, the following result is proved:

\begin{theo}
\label{theo1}
Let $\bbv=(v,\hv) \in W^\bbP$ be such that $c(\bbv)=1$. 
Let $\zeta \in X(T)^+$ and $\hzeta \in X(\hT)^+$ be such that
$v\inv(\zeta)_{|S} + \hv\inv(\hzeta)_{|S}$ is trivial. Then $v\inv(\zeta)$
and $\hv\inv(\hzeta)$ are dominant weights for $L$ and $\hL$
respectively. Moreover
$$
m_{G \subset \hG}(\zeta,\hzeta)=m_{L \subset \hL}(v\inv(\zeta),\hv\inv(\hzeta))\,.
$$
\end{theo}

The conclusion of Theorem~\ref{theo1} does not hold when
$c(\bbv)=2$. Our reduction result is a modification of this
conclusion which holds when $c(\bbv)=2$.

Consider the incidence variety
\begin{equation}
Y(\bbv)=\{(x,gB/B,\hg\hB/\hB)\, : \, x\in g \overline{BvP/P} \mbox{ and } \iota(x) \in \hg \overline{\hB\hv\hP/\hP}\}
\subset \GP \times \bbG/\bbB \,,
\label{eq:defYinc}
\end{equation}
endowed with its projection $\Pi\,:\,Y(\bbv) \longto \bbG/\bbB$
(in the sequel $\iota$ will be omitted and $\GP$ will be considered as a subvariety of $\hGP$).
Consider the ramification divisor class $K_\Pi$ of $\Pi$ (see
Section~\ref{sec:canonical} for details) and the branch divisor class
$[B_\Pi]=\Pi_*(K_\Pi)$ in $\Cl(\bbG/\bbB)$ (note that $\Pi$ is proper since $Y$ is
projective). Since $\bbG/\bbB$ is smooth, $[B_\Pi]$ is an element
of the Picard group and has an expression as
$$
[B_\Pi]=\cL_\GB(\sum_{\alpha\in\Delta}n_\alpha\varpi_\alpha)\otimes \cL_\hGB(\sum_{\halpha\in\hat\Delta}n_\halpha\varpi_\halpha)
$$
for some well defined integers $n_\alpha$ and $n_\halpha$.

Assume now that $c(\bbv)=2$. Then, we prove in Proposition~\ref{prop:mor2} that the integers
$n_\alpha$ and $n_\halpha$ are even.  Define $\theta$ and $\hat\theta$ in $X(T)$ and  $X(\hT)$ by setting
\begin{equation}
  \label{eq:2}
  \theta=\sum_\alpha\frac{n_\alpha} 2\varpi_\alpha\qquad\hat\theta=\sum_\halpha\frac{n_\halpha}2{\varpi}_\halpha\,.
\end{equation}

\bigskip

\begin{theo}
\label{main-theo}
Let $\bbv=(v,\hv)\in W^\bbP$ be such that $c(\bbv)=2$. 
Let $\zeta \in X(T)^+$ and $\hzeta \in X(\hT)^+$ be such that
$v\inv(\zeta)_{|S} + \hv\inv(\hzeta)_{|S}$ is trivial. 
Then
$$
m_{G \subset \hG}(\zeta,\hzeta) + m_{G \subset \hG}(\zeta-\theta,\hzeta-\hat\theta) =
m_{L \subset \hL}(v\inv\zeta,\hv\inv\hzeta)\,.
$$
\end{theo} 

Theorem \ref{main-theo} is proved in Section \ref{sub:proof-reduction}, using the same idea
as in \cite{ressayre-birational} and a formula for the sheaf $\Pi_* \cO_{Y(\bbv)}$ which
is proved in Section \ref{sec:finite} for any generically finite morphism of degree $2$.

Consider now the case when $\hG=G\times G$ corresponding to the tensor
product decomposition.
Then $\hP=P\times P$. Assume moreover that $P$ is cominuscule. 
Then Theorem~\ref{main-theo} can be improved by expressing $m_{G
  \subset G\times G}(\zeta,\hzeta) $ as an alternating sum of
multiplicities for the tensor product for $L$ as in Theorem~\ref{introtheo:LRalt} (see Corollary~\ref{cor:alt}).

\subsection{A formula for the branch divisor class}
\label{sec:branchdiv}

In this section, $P$ and $\hat P$ are any parabolic subgroups of $G$
and $\hat G$ respectively as in Diagram~\eqref{equ:inclusions}, but not
necessarily associated to a given one-parameter subgroup of $T$.
We let $K_\Pi:=K_{Y(\bbv)} - \Pi^* K_{\bbG/\bbB}$. If $c(\bbv) \neq 0$, it is the ramification divisor class of $\Pi$,
and we still call it in this way when $c(\bbv)=0$ although $\Pi$ is contractant in this case.
Under the assumption
$c(\bbv)=2$, $\Pi_* K_\Pi$ is twice the class of the branch divisor. More details are provided in Sections
\ref{sec:canonical} and \ref{sec:finite}.
For $\bbbeta\in\bbPhi$, we set:
\begin{equation}
\label{equa:def_h}
\hc(\bbbeta)=\scal{\rho_\bbG,\bbbeta^\vee}
\end{equation}

\begin{theo}
\label{theo:branch}
Let $\bbv\in W^\bbP$ such that $\ell(\bbv)=\dim(\hGP)$.  Consider the map
$$
\Pi\,:Y(\bbv)\longto
\bbG/\bbB\,,
$$
and its ramification divisor class $K_\Pi\in \Cl(Y)$.  Then, in
$\Cl(\bbG/\bbB)$, we have
$$
\Pi_*K_\Pi=\sum_{\bbv'} \left (
c_1(\cL(2\rho^\hL_{|T})) \cdot \delta^*(\tau_{\bbv'})+2c(\bbv)
-\sum_{\bbv''} (\hc(\bbgamma)+1) c(\bbv'') 
\right )
\cL_\bbGB(\varpi_\bbalpha) \,,              
$$
where the sums over $\bbv'$ and $\bbv''$ respectively run over the covering relations
\tikzlineleft{\bbv'}{\bbv}{\bbalpha}
and
\tikzdashedline{\bbv'}{\bbv''}{\gbbgamma}
in the Bruhat graph of $\bbG/\bbP$,
with $\bbalpha \in \bbDelta$ and
$\bbgamma \in \Phi^+(\bbG)$ is the twisted label. 
Here $c_1$ denotes the first Chern class. 
\end{theo}
\noindent
The first term $c_1(\cL(2\rho^\hL_{|T})) \cdot \delta^*(\tau_{\bbv'})$ is an integer
as an element of $A^{\dim(\GP)}(\GP) \simeq \Z$. We emphasise that it is an intersection number: see
\eqref{equ:delta*}.

\medskip

Let's outline the proof of this result, which is spread over several Sections of the paper.
The term $\Pi_* \Pi^* K_{\bbG/\bbB}$ is computed via the projection formula. To compute
$\Pi_* K_{Y(\bbv)}$, we observe that projection on $G/P$ expresses $Y(\bbv)$ as a
fibered product of two relative Schubert varieties. Our first task is to express
the canonical class of a relative Schubert variety, extending previous works on the absolute
case (see e.g. \cite[Theorem 4.2]{Raman} and \cite[Proposition 4.4]{perrin-small}).
This is done using the fact that (relative) Schubert varieties are normal and Cohen-Macaulay
and using the Bott-Samelson resolution (Proposition \ref{pro:kpi}).

Then, using the fibered product structure of $Y(\bbv)$, Proposition \ref{prop:KY}
expresses the canonical class of $Y(\bbv)$. We get the expression \eqref{eq:315} of $K_\Pi$
as a sum of line bundles coming from $\bbG/\bbB$ and $G/P$, and Weil divisors of the form
$[Y(\bbv')]$ for elements $\bbv'$ covered by $\bbv$ in the Bruhat order. To push forward
to $\bbG/\bbB$ those
Weil divisors as well as the line bundle coming from $G/P$, we express
the class of $Y(\bbv)$ in the Chow ring of $\bbG/\bbB \times G/P$,
thanks to a result by Brion valid for any
orbit closure of a spherical subgroup, see Propositions \ref{pro:incidence}
and \ref{pro:classY}.

\medskip

\begin{notation}
\label{nota:tensor_product_case}
In the case when $G\subset \hG=G^{n-1}$ and $\bbG/\bbP=(G/P)^n$ for some integer $n \geq 2$,
write $c_1^\alpha$ for the coefficient of $\Pi_*(K_\Pi)$ at ${\mathcal L}_{\bbG/\bbB}(\varpi_\alpha,0,\dots,0)$.
We let
$(\zeta,\hzeta)=(\zeta_1,\zeta_2,\dots,\zeta_n)$, $\bbv=(v_1,v_2,\dots,v_n)$,
$(\theta,\htheta)=(\theta^1,\theta^2,\dots,\theta^n)$ and
$m_{G\subset \hG}(\zeta,\hzeta)=m_G(\zeta_1,\zeta_2,\dots,\zeta_n)$.

Even more specifically, we call ``tensor product case'' the case when $n=3$, since then $m_G(\zeta_1,\zeta_2,\zeta_3)$ is the multiplicity
of $V_{\zeta_1}^*$ in $V_{\zeta_2} \otimes V_{\zeta_3}$. In this case, we set $\Gamma(G)=\Gamma(G\subset\hG)$.
\end{notation}

\begin{coro}
\label{coro:branch}
Let $n \geq 2$ be an integer.
Assume that $\hG=G^{n-1}$ and $\hP=P^{n-1}$. Let $\bbv=(v_1,\dots,v_n)\in {(W^P)}^n$ such that $\ell(\bbv)=(n-1)\dim(\GP)$ and
$c(\bbv)\neq 0$. Consider the map
$$
\Pi\,:Y(\bbv)\longto
\bbG/\bbB\,,
$$
and its ramification divisor class $K_\Pi\in \Cl(Y)$.  
Fix $\alpha\in\Delta$ and denote by $c_1^\alpha$ the coefficient of
$\Pi_*K_\Pi$ along $\cL((\varpi_\alpha,0,\dots,0))$.
If \tikzlineleft{v_1'}{v_1}{\alpha} is a covering relation in the
Bruhat graph of $G/P$ then 
$$
 c_1^\alpha=(n-1) c_1(\cL(2\rho^L))\tau_{v_1'}\tau_{v_2}\dots\tau_{v_n} +2 c(\bbv) - \sum_{\bbv''}
 (\hc(\bbgamma)+1) c(\bbv'')\,,
$$
where $\bbv'=(v_1',v_2,\dots,v_n)$ and the sum runs over the
covering relations
\tikzdashedline{\bbv'}{\bbv''}{\gbbgamma}.
Otherwise  $c_1^\alpha=0$.
\end{coro}

\begin{rema}
\label{rema:cominuscule}
In the case when moreover $G/P$ is cominuscule, strong and weak Bruhat order coincide, so
the covering relations
\tikzdashedline{\bbv'}{\bbv''}{\gbbgamma}
may be replaced by the covering relations
\tikzline{\bbv'}{\bbv''}{\gbbgamma}.
\end{rema}

\begin{proof}[Proof of Corollary~\ref{coro:branch}:]
 We only need to remark that since $\hG/\hP=(G/P)^{n-1}$, we have
 $\cL(2\rho^\hL_{|T})= (n-1) \cL_{\GP}(2\rho^L)$ and 
$\delta^*(\tau_{\bbv'})=\tau_{v_1'}\tau_{v_2}\dots\tau_{v_n}$.
\end{proof}

Theorem \ref{theo:branch} is proved in Section \ref{sub:branch}. Here we make one comment.
Let $\bbv'=(v',\hv') \in W^\bbP$ such that
$\ell(\bbv') = \dim(\hGP)+1$ like in the first sum of Theorem~\ref{theo:branch}.
Set
$$
C(\bbv') =g \cdot \overline{Bv' P/P} \cap \iota\inv(\hg \cdot \overline{B\hv' P/P})\,,
$$
for general elements $g \in G,\hg \in \hG$. Then $C(\bbv') \subset \GP$ is a curve.
The degree of the pullback of the
anticanonical bundle $\cL(2\rho^\hL_{|T})$ on $C(\bbv')$ is the first term
appearing in the expression of $\Pi_*K_\Pi$. By Chevalley's formula
applied to $c_1(\cL(2\rho^\hL_{|T}))\cdot \tau_{v_1'} $, we have an explicit formula
\begin{equation}
\label{equ:curve-degree}
c_1(\cL(2\rho^\hL_{|T})) \cdot \delta^*(\tau_{\bbv'})=
\sum_{\bbv''} \scal{2\rho^\hL_{|T},\gamma^\vee} c(\bbv'')\,,
\end{equation}
\vspace{-4mm}

\noindent
where the sum runs over the covering relations
\tikzdashedline{\bbv'}{\bbv''}{\ggamma}
with $\gamma \in \Phi^+(G)$.

\medskip

If $c(\bbv)=1$ then $\Pi$ is birational so
$\Pi_*(K_\Pi)=0$. Theorem~\ref{theo:branch} gives a funny but combinatorially
involved necessary condition for the equality $c(\bbv)=1$. In
Example~\ref{ex:poincare}, we check a family of examples. 

At the opposite, if $c(\bbv)\geq 2$ then $\Pi_*(K_\Pi)\neq 0$. 
Indeed, if $\Pi_*(K_\Pi)=0$ then by removing the singular locus of
$Y(\bbv)$ and the ramification divisor $R_\Pi$, we get an unramified
cover onto the complement in $\bbG/\bbB$ of a codimension two
subvariety. 
But, this complement is simply connected as  $\bbG/\bbB$ is by the Bruhat
decomposition.
This is a contradiction.

Assume now that  $c(\bbv)=0$. In this case, $\Pi$ is contractant, thus $\Pi_* \Pi^* K_{\bbG/\bbB}=0$.
However, we have no information on $\Pi_* K_{Y(\bbv)}$, and this can be non zero: see Example~\ref{exple:sl3}.
 
\begin{exple}
\label{ex:sl3}
Let $G=\SL_3,\hG=G \times G,P=B,\hP=B \times B$, and $\bbv=(s_2s_1, s_2s_1, s_1s_2)$. Then $\Pi_*K_\Pi=0$.
\end{exple}
\begin{proof} 
We use Corollary~\ref{coro:branch} to check this. Recall that the Bruhat graph of $SL_3$ is depicted in
Figure \ref{fig:Bruhatgraph}, and that the Poincaré dual of $\tau_{s_i}$ is $\tau_{w_0s_i}=\tau_{s_is_{3-i}}$.
Fix first $\bbalpha=(\alpha_1,0,0)$, which implies $\bbv'=(s_1s_2s_1,s_2s_1,s_1s_2)$.
Then $2 c_1(\cL(2\rho))\tau_{v_1'}\dots\tau_{v_3'} +2 c(\bbv)=2
c_1(\cL(2\rho))(\tau_{s_1}+\tau_{s_2}) +2\times 1=2\scal{2\rho,\alpha_1^\vee+\alpha_2^\vee}+2=10$.
The possibilities for $\bbgamma$ with $c(\bbv'')\neq 0$ are
$(0,\alpha_1+\alpha_2,0)$, $(0,0,\alpha_1+\alpha_2)$, $(\alpha_1,0,0)$
and $(\alpha_2,0,0)$.
In these cases $c(\bbv'')=1$ and 
$\hc(\bbgamma)=2$ twice and $\hc(\bbgamma)=1$ twice. We get
$c_1^{\alpha_1}=0$ as expected. 

Fix now $\bbalpha=(0,0,\alpha_2)$ with $\bbv'=(s_2s_1,s_2s_1,s_2s_1s_2)$.
Then $2 c_1(\cL(2\rho))\tau_{v_1'}\dots\tau_{v_3'} +2 c(\bbv)=2
c_1(\cL(2\rho))(\tau_{s_1}) +2\times
1=6$.
The possibilities for $\bbgamma$ with $c(\bbv'')\neq 0$ are
$(\alpha_1,0,0)$, $(0,\alpha_1,0)$
and $(0,0,\alpha_1)$
and give $3\times 2=6$ for the sum.

The remaining cases are either trivial or obtained by exchanging $v_1$ and
$v_2$. 
\end{proof}

\begin{exple}
\label{exple:sl3}
Let $G=\SL_3,\hG=G \times G,P=B,\hP=B \times B$, and $\bbv=(s_2s_1,s_1s_2s_1, s_1)$. 
We have $c_1^{\alpha_1}=-2$.
\end{exple}
\begin{proof} 
For $c_1^{\alpha_1}$, with the notation of Corollary~\ref{coro:branch}, we have
$\bbv'=(s_1s_2s_1,s_1s_2s_1,s_1)$ and the elements $\bbv''$ such that $c(\bbv'') \neq 0$
are the three triples
$(s_1s_2s_1,s_1s_2s_1,e),(s_1s_2s_1,s_1s_2,s_1),(s_1s_2,s_1s_2s_1,s_1)$.
The right hand-side of the formula given in Corollary~\ref{coro:branch} is thus 
$2c_1(\cL(2\rho))\tau_{s_1} +2 c(\bbv)-3 \times (1+1)=4+0-6=-2$.
\end{proof}

We now consider a family of examples again in the birational case:

\begin{exple}
\label{ex:poincare}
Let $\hG=G \times G$ and $\hP=P \times P$. Let $v \in W^P$ be arbitrary and let
$\bbv=(v,v^\vee,w^P)$. Then we have $\Pi_* R_\Pi=0$.
\end{exple}
\begin{proof}
We check that  Corollary~\ref{coro:branch} gives $\Pi_* R_\Pi=0$. 
Up to exchanging $v$ and $v^\vee$, the only case to consider is when
$\bbalpha=(\alpha,0,0)$ with $\alpha\in\Delta$,  $s_\alpha v\in W^P$
and $\ell(s_\alpha v)=\ell(v)+1$.
Set $\gamma_0=-v\inv \alpha$ and $\bbv'=(s_\alpha v,v^\vee,w^P)$.

Consider the Chevalley formula \eqref{equa:chevalley} for
$c_1(\cL(2\rho^L))\cdot \tau_{v^\vee}$. The only summand surviving after
multiplication by $\tau_{s_\alpha v}$ is
$\scal{2\rho^L,\gamma}\tau_{(s_\alpha v)^\vee}$. 
Moreover, $c(\bbv'')=1$ for this term and $\gamma=w_{0,P}(\gamma_0)$.
Finally $2 c_1(\cL(2\rho))\tau_{v_1'}\dots\tau_{v_3'} +2 c(\bbv)
=\scal{4\rho^L,w_{0,P}(\gamma_0^\vee)}+2=\scal{4\rho^L,\gamma_0^\vee}+2$.
 
The only case with $c(\bbv'')\neq 0$ and $\bbgamma$ in the first
factor is $\bbgamma=(\gamma_0,0,0)$. 
Its contribution in the sum is $\scal{\rho,\gamma_0^\vee}+1$.
Similarly, with $\bbgamma$ in the second factor, we get
$\scal{w_{0,P}(\rho),\gamma_0^\vee}+1$.

Assume now that $\gamma$ is in the third copy. It is a descent of
$w^P$, hence an element of $\Delta-\Delta(L)$. 
We have $\tau_{w^Ps_\gamma} =c_1(\cL(\varpi_\gamma))$, so
the associated coefficient $c(\bbv'')$ is
$\scal{\varpi_\gamma,\gamma_0^\vee}$ by the Chevalley formula. The term $\hc(\gamma)+1$ equals
$\scal{\rho,\gamma^\vee}+1=2$. The contribution of these terms is
$\sum_{\gamma\in \Delta-\Delta(L)}2
\scal{\varpi_\gamma,\gamma_0^\vee}=2\scal{\rho-\rho_L,\gamma_0^\vee}$.

Finally, 
$$c_1^\alpha = 4\scal{\rho^L,\gamma_0^\vee}+2-(\scal{\rho,\gamma_0^\vee}+1)-
(\scal{w_{0,P}(\rho),\gamma_0^\vee}+1)-2\scal{\rho-\rho_L,\gamma_0^\vee}
$$
which vanishes, since $\rho^L=\rho-\rho_L$ and $w_{0,P}(\rho)=\rho^L-\rho_L$.
\end{proof}

In Section \ref{sec:examples}, we give several examples of Theorem \ref{theo:branch} in the case $c(\bbv)=2$,
which is the case we are mainly interested in. Let us now give an
example with $c(\bbv)=3$:

\begin{exple}
\label{exam:c=3}
In $\Gr(4,8)$, let $u=[35681247]=(s_2s_3s_5s_4)^\vee,v=[24681357]=(s_2s_4s_6s_3s_5s_4)^\vee$ and $w=[24681357]=(s_2s_4s_3s_6s_5s_4)^\vee$.
Then $c(u,v,w)=c_{211,321}^{4321}=3$. We computed with Theorem~\ref{theo:branch} that
$$[B_\Pi]=\cL_{(\GB)^3}(4(\varpi_2+\varpi_4+\varpi_6),4(\varpi_2+\varpi_4+\varpi_6),4\varpi_3+6\varpi_6)\,.$$
In particular, $[B_\Pi]$ is not divisible by $3$, which shows that the following Lemma \ref{degree2} and Proposition
\ref{prop:mor2} don't have obvious analogues for morphisms of degree $3$.
\end{exple}

\section{Preliminaries in Algebraic Geometry}

\subsection{Recollections on Intersection theory}
\label{sec:intersection-theory}

We now recall useful notions in Intersection Theory and
give some details on degree two
morphisms and the associated ramification and branch divisors.

\medskip

Let $f:Y\longto X$ be a dominant morphism between irreducible
varieties of the same dimension. We say that $f$ is {\it
  generically finite}. 
The {\it degree of $f$} is defined to be $\deg(f)=[\CC(Y):\CC(X)]$.

Let $X$ be a quasi-projective irreducible variety. We denote by $\Cl(X)$ the Weil
divisor class group, by $\CaCl(X)$ the Cartier divisor class group and
by $\Pic(X)$ the Picard group. 
By \cite[Prop. II.6.15]{hartshorne}, we have a canonical isomorphism
$\CaCl(X)\simeq\Pic(X)$, mapping $D$ on $\cO(D)$.
There is also a morphism $c_1\,:\,\Pic(X)\longto \Cl(X)$ (see
\cite[p. 30]{fulton}) that is neither injective nor surjective. 
Nevertheless, if $X$ is normal, $c_1$ is injective and
$\CaCl(X)=\Pic(X)$ can be seen as a subgroup of $\Cl(X)$.
Moreover if $X$ is locally factorial (e.g. is smooth) then this
morphism is an isomorphism. 

We also consider the Chow ring $A^*(X)$ and identify $\Cl(X)$ with
$A^1(X)$.
Given two irreducible subvarieties $Z_1$ and $Z_2$ that intersect
transversally in $X$, we have
\begin{equation}
  \label{eq:Chowtransv}
  [Z_1\cap Z_2]=[Z_1] \cdot [Z_2]\mbox{ in }A^*(X).
\end{equation}
Here $ [Z_1\cap Z_2]$ denotes the sum of the classes of the
irreducible components of $Z_1\cap Z_2$.

\subsection{Canonical, ramification and branch divisors}
\label{sec:canonical}

Set $n=\dim(X)$.
If $X$ is smooth, the canonical divisor class $K_X\in\Pic(X)$ on $X$ is defined to be the
line bundle $\wedge^n T^*X$. 
Assume that $X$ is normal. 
Let $X_\reg$ denote the regular locus of $X$ and $TX_{\reg}$ denote
its tangent bundle. Then $\wedge^n T^*X_\reg$ is a line bundle $\cL$ on $X_\reg$.
The first Cern class $c_1(\cL)\in\Cl(X_\reg)=\Cl(X)$ is called the
canonical divisor class of $X$ and denoted by $K_X$. 
The canonical sheaf  ${\mathcal K}_X$ is defined to be the pushforward
of the sheaf $U\mapsto \cL(U)$ by the inclusion $X_\reg\subset X$.

Let $f:Y\longto X$ be a generically finite morphism between irreducible
varieties. 
Assume moreover that $Y$ is normal and $X$ is smooth. 
The determinant of the restriction of
$f$ to the smooth locus $Y_\reg$ of $Y$ define a Cartier divisor
$R_f$ in $Y_\reg$. By construction, the associated line bundle
$\cO(R_f)\in\Pic(Y_\reg)$ is $K_f:=K_Y-f^*K_X$.
Its closure, still denoted by $R_f$ has a class
$[R_f]$ in $\Cl(Y)$ called ramification divisor class of $f$.
We have
\begin{equation}
\label{eq:3}
[R_f]=K_Y-f^*K_X\in\Cl(Y),
\end{equation}
where $f^*\,:\Pic(X)\longto\Pic(Y)$ is the pullback of line bundles.  
Observe that $R_f$ is an effective Weil divisor. We write $R_f$ as
$R_f = \sum_{i \in I} a_i E_i + \sum_{j \in J} c_j R_j$, with prime Weil divisors $E_i,R_j$ such that the divisors $E_i$
are contracted and the divisors $R_j$ are not, and positive integers $a_i,c_j$. We define the Weil divisors
\begin{equation}
\label{equa:divisors}
E_f=\sum_{i \in I} a_i E_i \ , \ R^1_f=\sum_{j \in J} c_j R_j \,.
\end{equation}
We define the branch divisor $B_f:=f_*R^1_f=\sum_{j \in J} c_j\deg(f_{|R_j})f(R_j)$, and its
class $[B_f]=f_*[R_f]=\sum_{j \in J} c_j\deg(f_{|R_j})[f(R_j)]\in \Cl(X)$
(since $f$ is proper, $f_*:\Cl(Y) \longto \Cl(X)$ is well defined).

\subsection{Pushforward of the structure sheaf under a finite morphism of degree 2}
\label{sec:finite}

We consider a finite 
morphism $f\,:\,S\longto X$ 
and study $f_* \cO_S$. 
The next lemmas are well-known to specialists, see for example \cite[Lemma 1.1]{lazarsfeld} or \cite{pardini}. We provide
details for the singular case we are interested in.

As it is explained in \cite[p.9]{BK:book}, we have:
\begin{lemma}
\label{finite}
Let $f\,:\,S\longto X$ be a finite surjective morphism of degree $d$
and assume that $X$ is normal. Then there exists a splitting $f_* \cO_S = \cO_X \oplus L$ with $L$ a coherent sheaf on $X$.
\end{lemma}

We now consider the degree 2 case and exploit the fact that $f$ is a cyclic cover.

\begin{lemma}
\label{degree2}
Let $f\,:\,S\longto X$ be a finite surjective morphism of degree $2$
and assume that $S$ is normal and $X$ is normal and locally factorial.
Then, in the splitting $f_* \cO_S = \cO_X \oplus L$ as above, $L$ is a line bundle,
and $[L^{\otimes 2}] = -[B_f] \in \Pic(X) \simeq \Cl(X)$.
Moreover, $B_f$ and $R_f$ are reduced Weil divisors.
As Cartier divisors, we have $f^* \cO_X(B_f)=\cO_Y(2R_f)$.
\end{lemma}

\begin{proof}
Let us restrict ourselves over an affine open subset $U\subset X$ and set 
$A=\CC[U]$ and $B=\CC[f\inv(U)]$. 

Let $b\in B \setminus \Frac(A)$.
The ring $B$ being integral over $A$, there exists 
a monic polynomial $P$ in $A[T]$ of minimal degree such that
$P(b)=0$.
Then $P$ is irreducible in $A[T]$.
By assumption, $A$ is unique factorization domain. Then, $P$ is also irreducible as an element of
$\CC(X)[T]$, by \cite[Exercise 3.4.c]{eisenbud}.
It follows that $P$ is the minimal polynomial of $b\in \CC(S)$ over
$\CC(X)$.
But $[\C(S):\C(X)]=2$, hence $P$ has degree 2.

Write $P=T^2+a_1T+a_2$ with $a_1$ and $a_2$ in $A$. 
Up to changing $b$ by $b-\frac{a_1}2$, one may assume that $a_1=0$.
Assume that $-a_2=t^2c$ in $A$.
Set $b=tb'$. Then $b'^2-c=0$, and $b'\in\CC(X)$ is integral over
$A$. By normality of $S$, this implies that $b'\in B$. 
In particular, up to changing $b$, one may assume that the minimal
polynomial of $b$ is $T^2-a$ with some square free element $a\in A$.

\medskip
We claim  that $B=A \oplus Ab$.
Let $y=\alpha+\beta b\in B$ with $\alpha,\beta\in\Frac(A)$. 
The matrix of $x\mapsto xy$ is
$
\begin{pmatrix}
  \alpha&a\beta\\
\beta&\alpha
\end{pmatrix}.
$
Then $\Tr(y)=2\alpha$ and $N(y)=\alpha^2-a\beta^2$ are
integral over $A$. Hence $\alpha\in A$ and $a\beta^2\in A$.
Since $A$ is factorial and $a$ is square free, this implies that
$\beta\in A$.
This proves the claim.

\medskip
It follows that $f\inv(U)=\Spec(A[t]/(t^2-a))$ and $f$ corresponds to
the inclusion $A\subset A[t]/(t^2-a)$. 
Then $L$ is the kernel of $\Tr$, namely $tA \subset A[t]/(t^2-a)$. This is a free $A$-module.
This proves that $L$ is a line bundle.

A local computation as in \cite[Lemma 16.1]{barth} shows that $R_f$
is locally defined by $t$, so it is reduced. We get $f^* \cO_X(B_f)=\cO_S(2R_f)$, since
these locally free sheaves are generated by $a=t^2$.
Moreover, $a$ being square-free, is an equation of $f(R_f)\cap U$ in
$\CC[U]$.
It follows that  $L^{\otimes 2}(U)$ identifies 
with $aA=\cO(-B_f))(U) \subset \C(X)$. This proves that $L^{\otimes 2} \simeq \cO(-B_f)$.
Finally, $f$ restricts to an isomorphism from the hypersurface defined by the equation $t=0$ to $B_f \cap U$, proving that
the Weil divisor $B_f$ is reduced.
\end{proof}

\begin{rema}
\label{rema:c=3}
The degree $2$ is a special case since any finite morphism of degree $2$ is a cyclic covering.
Already in degree $3$, this is no longer true, and the above Lemma fails: see
 Example \ref{exam:c=3} below.
\end{rema}

\subsection{\texorpdfstring{Generically finite morphisms of degree $2$}{sub:gen-finite}}
\label{sec:genfinite}

We now consider the case of a degree 2 map, but not necessarily
finite: we let $f:Y \longto X$ be a generically finite morphism of degree $2$, with $X$ and $Y$ projective,
$Y$ irreducible normal, and $X$ normal and locally factorial. We let
\begin{tikzpicture}[baseline=-1.2mm]
    \node (Y) at (0,0) {$Y$};
    \node (S) at (1.5,0) {$S$};
    \node (X) at (3,0) {$X$};
    \draw[->] (Y) to node[midway,above]{$p$} (S);
    \draw[->] (Y) to[bend right] node[midway,below]{$f$} (X);
    \draw[->] (S) to node[midway,above]{$h$} (X);
\end{tikzpicture}
be the Stein factorization of $f$.

\begin{lemma}
\label{lem:zariski}
With the above notation, assume we have an irreducible divisor $C_S \subset S$.
Then there is a unique
irreducible component $C_Y$ of $p^{-1}(C_S)$ which maps onto $C_S$, the restriction $C_Y \stackrel{p}{\longto} C_S$
is birational, and we have the equivalence
$$
C_Y \subset R_f \Longleftrightarrow C_S \subset R_h\,.
$$
\end{lemma}
\begin{proof}
Since $p$ is surjective, there is an irreducible component $C_Y$ of $p\inv(C_S)$ mapping onto $C_S$.
This component is unique and $C_Y \stackrel{p}{\longto} C_S$ is
birational since $p$ has connected fibers.
Since $Y,S$ and $X$ are normal, for a generic point $y \in C_Y$, $Y,S,X$ are normal and locally factorial at
$y,p(y),f(y)$, respectively.
Since $C_Y \stackrel{p}{\longto} C_S$ is dominant, $T_yp$ is invertible,
by Zariski's Main Theorem again. Therefore $T_yf$ will be degenerate if and only if
$T_{p(y)}h$ is degenerate. In other words, $C_Y \subset R_f$ if and only if
$C_S \subset R_h$.
\end{proof}

In the next proposition, we use the prime divisors $E_i$ and the divisor $R^1_f$ introduced in Section \ref{sec:canonical}.

\begin{prop}
\label{prop:mor2}
Let $f:Y \longto X$ be a generically finite morphism of degree $2$, with $X$ and $Y$ projective,
$Y$ irreducible normal, and $X$ normal and locally factorial.
Then there is a line bundle $L$ on $X$ such that $f_* \cO_Y = \cO_X
\oplus L$. 
Moreover, in $\Pic(X)=\Cl(X)$,
we have $2L+f_*(K_Y)-f_*f^*(K_X)=0$ and $B_f=B_h$ is reduced.
As Weil divisors on $Y$, we have
\begin{equation}
\label{equa:f*}
f^* \cO_X(B_f)=2R^1_f+\sum_{i \in I} b_i E_i
\end{equation}
for some non-negative integers $b_i$.
\end{prop}
\begin{proof}
We keep the notation of the beginning of the section, and we want to prove that
$p_*R_f=R_h$.
By Lemma~\ref{degree2}, $R_h$ is reduced:
$R_h=\sum_{C_S\subset R_h}C_S$.
By Lemma \ref{lem:zariski}, for a component $C_S$ of $R_h$, there is a unique component $C_Y$ of $R_f$
such that $p(C_Y)=C_S$. Since $C_Y \stackrel{p}{\longto} C_S$ is birational, we get $p_*C_Y=C_S$.
Moreover, $p$ is an isomorphism at a generic point of $C_Y$, so $C_Y$ is a reduced component of $R_f$.
For a contracted component $C_Y'$ of $R_f$, we have $p_*[C_Y']=0$. Thus
$p_*R_f = \sum_{C_Y \subset R_f} p_*C_Y = \sum_{C_S\subset R_h}C_S=R_h$.
Applying $h_*$ and using Lemma \ref{degree2}, we get that $B_f=f_*R_f=h_* R_h = B_h$ is reduced.

By Lemma~\ref{finite}, let $L$ be such that $h_* \cO_S=\cO_X \oplus L$.
By Lemma~\ref{degree2}, $h_*[R_h]=-2L$, so $f_*[R_f]=[B_f]=-2[L]$.
It remains to prove that $f_*\cO_Y=h_*\cO_S$. This follows from
Zariski's Main Theorem that asserts that $p_*\cO_Y=\cO_S$.

Let us write the Weil divisor $f^*\cO_X(B_f)$ as $\sum_{i \in I} b_i E_i  + \sum_{j \in J} d_j R_j$.
Then \eqref{equa:f*} is equivalent to $b_i \geq 0$ for $i \in I$ and $d_j=2$ for $j \in J$.
We write the Cartier divisor $\cO_X(B_f)$ as
$(U_i,\xi_i)$, with $\xi_i$ a local equation of $B_f$ on $U_i$.
Then the Cartier divisor $f^* \cO_X(B_f)$ is defined by $(f\inv(U_i),\xi_i \circ f)$.
Since $\xi_i \circ f$ is regular, we have $b_i \geq 0$ for all $i$. Moreover, if $j \in J$, then
restricting to a neighborhood of a point of $R_j$ where $p$ is a local isomorphism
and using Lemma \ref{degree2}, we get $d_j=2$.
\end{proof}

\medskip
We now describe the branch locus as a set.

\begin{prop}
\label{pro:supportBPi}
Let $f:Y \longto X$ be a generically finite morphism of degree $2$, with $X$ and $Y$ projective,
$Y$ irreducible normal, and $X$ normal and locally factorial.
Then 
$$
\Supp(B_f)=\{x\in X\;:\; f\inv(x)\mbox{ is connected}\}
$$
\end{prop}
\begin{proof}
We keep our notation. First, we claim that 
$
B_h=B_f.
$
This is equivalent to saying that $h(R_h)=f(R^1_f)$.
These two varieties are closed in $X$ and of pure codimension one.
Lemma \ref{lem:zariski} yields a bijection between the irreducible components of $R_h$ and
those of $R^1_f$, implying that $h(R_h)$ and $f(R^1_f)$ have the same irreducible components,
proving the claim.

Observe now that $B_h$ is the set of $x\in X$ such that $h\inv(x)$ is
a point (see the proof of Lemma~\ref{degree2}). The fibers of $p$ being connected by
Zariski's Main Theorem, the proposition follows.
\end{proof}

\section{Preliminaries on Schubert varieties}
\label{sec:preSchub}

In this section, we fix a semisimple group $G$ with maximal torus $T$ and Borel
subgroup $B$ containing $T$. 
Let $P\supset B$ be a standard parabolic subgroup. We
are interested in the Schubert varieties in $G/P$, their relative
versions in $G/B\times G/P$ and their Bott-Samelson resolutions. 
We collect both combinatorial and geometric material.

\subsection{Reminder on Schubert varieties}
\label{sub:schubvar}

For $w\in W$ (resp. $w\in W^P$), denote by $\cX_w^B=\overline{BwB/B}$
(resp. $\cX_w^P=\overline{BwP/P}$) the associated
Schubert subvariety of $G/B$ (resp. $G/P$).
We use notation for the Bruhat orders introduced in Section \ref{sec:notations}.

\begin{lemma}
Let $v\in W^P$. Then
\begin{enumerate}
 \item \label{ass:normal}
The variety $\cXPv$ is normal. In particular,
   $\Pic(\cXPv)=\CaCl(\cXPv)$ embeds in $\Cl(\cXPv)$.

\item\label{ass:Cl} 
We have 
$$\Cl(\cXPv)=\bigoplus
\Z [\cX^P_{v'}],
$$
where the direct sum runs over the
covering relations
\tikzdashedline{v}{v'}{}
in the Bruhat graph of $G/P$.
\end{enumerate}
\end{lemma}

\begin{proof}
Good references for assertions~\ref{ass:normal} and~\ref{ass:Cl} are
\cite[XII, Lemme 75]{mathieu} and \cite[§3.2]{BK:book}.
\end{proof}

\subsection{Inversion sets}

\subsubsection{Inversion sets of Weyl group elements.}
\label{sec:Bruhatgraph}

As before, we use notation $W$, $W^P$, $W_P$, $\Delta$, $\Phi$ and $\Phi^+$.
The inversion set $\Phi(w)$ of an element $w$ of $W$ is defined to be
$\Phi(w)=\Phi^+\cap w\inv(\Phi^-)$ and it satisfies
$\ell(w)=\sharp\Phi(w)$. 
One gets a map
$$
\begin{array}{ccl}
  W&\longto&\Sub(\Phi^+)\\
w&\longmapsto&\Phi(w)
\end{array}
$$
($\Sub(\Phi^+)$ denotes the power set of $\Phi^+$) that is injective.

\medskip

The strong and left weak Bruhat orders are characterized by 
\begin{equation}
\label{equa:bruhat-orders}
\begin{array}{ccl}
v\bole w&\iff&\cX_v^B\subset \cX_w^B;\\
v\wbole w&\iff&\Phi(v)\subset \Phi(w).
\end{array}
\end{equation}
The second point is \cite[Proposition 1.3.1]{bb}. The
direct implication is a consequence of the easy fact that 
$w=s_\alpha v$ and $\ell(w)=\ell(v)+1$ imply that 
\begin{equation}
\Phi(w)=\Phi(v)\ccup\{v\inv(\alpha)\}.\label{eq:Invsalphav}
\end{equation}

A reduced expression $w=s_{\alpha_1}\dots s_{\alpha_\ell}$ corresponds
to a sequence $e \wbole s_{\alpha_\ell} \wbole s_{\alpha_{\ell-1}}s_{\alpha_\ell} \wbole \cdots \wbole w$ of
weak covering relations and hence to a sequence
$\emptyset\subset\cI_1\subset\cI_2\subset\cdots\subset\Phi(w)$ of
inversion sets of $\Phi^+$ with cardinality increasing by one at
each step. 
In particular, given such a reduced expression we get a numbering
\begin{equation}
\label{equ:gamma_i}
\Phi(w)=\{\gamma_1,\dots,\gamma_\ell\} \mbox{ with }
\gamma_i=s_{\alpha_\ell} \cdots s_{\alpha_{\ell-i+2}}(\alpha_{\ell-i+1})\,.
\end{equation}
Note that $w\inv=s_{\alpha_\ell}\dots s_{\alpha_1}$ is also a reduced
expression, thus
\begin{equation}
\label{equ:beta_i}
\Phi(w\inv)=\{\beta_1,\dots,\beta_\ell\} \mbox{ with }
\beta_i=s_{\alpha_1} \cdots s_{\alpha_{i-1}}(\alpha_i)\,.
\end{equation}

\begin{lemma}\label{lem:Bruhatorder}
  Let $v'\bole v$ in $W$ with $\ell(v)=\ell(v')+1$ be a covering
  relation in the Bruhat order. 
Fix a reduced expression $\bv=s_{\alpha_1}\dots s_{\alpha_\ell}$ of $v$. 
 Then there exists a unique $k$ such that  a reduced
expression of $v'$ is obtained from $\bv$ by deleting $s_{\alpha_k}$.
\end{lemma}
\begin{proof}
This is the content of \cite[Lemma 1.3.1]{bb}.
A geometric interpretation of this result can be obtained using the
Bott-Samelson resolution. Indeed, if the statement would be false, one
would get nontrivial finite fibers for the Bott-Samelson
resolution. By Zariski's Main Theorem, this contradicts the
normality of the Schubert varieties.
\end{proof}

\medskip

A useful result on the inversion sets is:

\begin{lemma}
See e.g. \cite[Corollary~1.3.22.3]{kumar}.
\label{lem:wrho}
Let $w \in W$. We have
$$
\sum_{\beta \in \Phi(w^{-1})} \hspace{-3mm} \beta \ = \rho - w(\rho)\,.
$$
\end{lemma}

\noindent
We also use the following variation, see \cite[Lemma 4.18]{perrin-small}. The roots $\beta_i$ are those defined by
\eqref{equ:beta_i} and the function $\hc$ by \eqref{equa:def_h}.

\begin{lemma}
  \label{lem:sumgamma}
Let $v'\bole v$ in $W$ with $\ell(v)=\ell(v')+1$.
Let $\bv=s_{\alpha_1}\dots s_{\alpha_\ell}$ be a reduced expression of
$v$ and $1\leq i\leq\ell$ such that $v'$ is obtained by deleting
$s_{\alpha_i}$.
Write $\Phi(v\inv)=\{\beta_1,\dots,\beta_\ell\}$ as in \eqref{equ:beta_i}.

Then, $v'=s_{\beta_i}v$ and
$\Phi(v'\inv)=\{\beta_1,\dots,\beta_{i-1},s_{\beta_i}(\beta_{i+1}),\dots,
s_{\beta_i}(\beta_\ell)\}$. Moreover, 
$$
\scal{\sum_{k=i}^\ell\beta_k,\;\beta_i^\vee}=\hc(-v\inv \beta_i)+1.
$$
\end{lemma}

\begin{proof}
One can easily check the first assertion. We provide an argument for the second assertion simpler than that in \cite{perrin-small}.
By Lemma~\ref{lem:wrho}, we have 
$$
\begin{array}{ll}
  v'(\rho)-v(\rho)&=\sum_{\theta\in \Phi(v\inv)}\theta-\sum_{\theta\in
                \Phi(v'\inv)}\theta\\
&=\beta_i+\sum_{k=i+1}^\ell(\beta_k-s_{\beta_i}(\beta_k))\\
&=\beta_i+\sum_{k=i+1}^\ell\scal{\beta_k,\beta_i^\vee}\beta_i.
\end{array}
$$ 
Set $S=\scal{\sum_{k=i}^\ell\beta_k,\;\beta_i^\vee}=2+\scal{\sum_{k=i+1}^\ell\beta_k,\;\beta_i^\vee}$. We get 
$$
\scal{ v'(\rho)-v(\rho),\beta_i^\vee}=2(S-1).
$$
On the other hand, 
$$
\scal{ v'(\rho)-v(\rho),\beta_i^\vee}=\scal{ s_\beta v(\rho)-v(\rho),\beta_i^\vee}=-2\scal{v(\rho),\beta_i^\vee},
$$
since $s_{\beta_i}(\zeta)=\zeta-\scal{\zeta,\beta_i^\vee}\beta_i$ for any weight $\zeta$.
The lemma follows.
\end{proof}

\subsubsection{Diamond lemmas.}

We first state the following diamond lemma in the weak Bruhat graph.

\begin{lemma}
  \label{lem:diam}
Let $\alpha,\beta,\gamma,\delta\in\Delta$ such that 
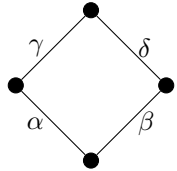
\begin{center}
  \begin{tikzpicture}
    \node[dnode,fill=black,inner sep=0pt] (a) at (0,1) {}; 
    \node[dnode,fill=black,inner sep=0pt] (b) at (-1,0) {}; 
     \node[dnode,fill=black,inner sep=0pt] (c) at (1,0) {};
     \node[dnode,fill=black,inner sep=0pt] (d) at (0,-1) {}; 
    \draw (a) --  node[right] {$\delta$} (c); 
     \draw (a) -- node[left] {$\gamma$} (b);
     \draw (b) -- node[left] {$\alpha$} (d); 
     \draw (c) -- node[right] {$\beta$} (d);
  \end{tikzpicture}
\end{center}
is a subgraph of the weak Bruhat graph.
Then $\alpha=\delta$, $\beta=\gamma$ and $s_\alpha s_\beta=s_\beta s_\alpha$.
\end{lemma}

\begin{proof}
Let $v$ denote the Weyl group element of the  bottom vertex.
We have $s_\gamma s_\alpha v=s_\delta s_\beta v$ whose length is
$\ell(v)+2$.
In particular  $s_\gamma s_\alpha =s_\delta s_\beta$ has length
2:
one may assume that $v$ is trivial. 
 Then, the four inversion sets of the vertices are
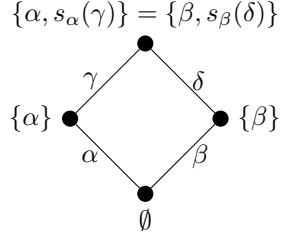
\begin{center}
  \begin{tikzpicture}
    \node[dnode,fill=black,inner sep=0pt] (a) [label=above:{$\{\alpha,s_\alpha(\gamma)\}=\{\beta,s_\beta(\delta)\}$}] at (0,1) {}; 
    \node[dnode,fill=black,inner sep=0pt] (b) [label=left:{$\{\alpha\}$}] at (-1,0) {}; 
     \node[dnode,fill=black,inner sep=0pt] (c) [label=right:{$\{\beta\}$}] at (1,0) {};
     \node[dnode,fill=black,inner sep=0pt] (d) [label=below:{$\emptyset$}] at (0,-1) {}; 
    \path (a) edge node[right] {$\delta$} (c); 
     \path (a) edge node[left] {$\gamma$} (b);
     \path (b) edge node[left] {$\alpha$} (d); 
     \path (c) edge node[right] {$\beta$} (d);
  \end{tikzpicture}
\end{center}
Since $\alpha\neq\beta$, we deduce 
that $\alpha=s_\beta(\delta)=\delta-\scal{\delta,\beta^\vee}\beta$.
Hence $\alpha,\beta$ and $\delta$ are linearly dependant simple
roots. Since $\alpha\neq \beta$ 
we deduce that $\alpha=\delta$. Now, the linear relation implies
$\scal{\delta,\beta^\vee}=0$, so $s_\delta$ and $s_\beta$ commute.
By symmetry this ends the proof.
\end{proof}

For $\alpha \in \Delta$, let $P_\alpha \subset G$ be the corresponding
minimal standard parabolic subgroup.
We prove a diamond lemma mixing strong and weak Bruhat orders:

\tikzcdset{column sep/tiny=0.1cm}
\begin{lemma}
\label{lemm:bruhat2}
Let $w \in W, \alpha \in \Delta, \gamma \in \Phi^+$ with $s_\alpha
w\neq ws_\gamma$. Let $P$ be a standard parabolic subgroup of $G$.
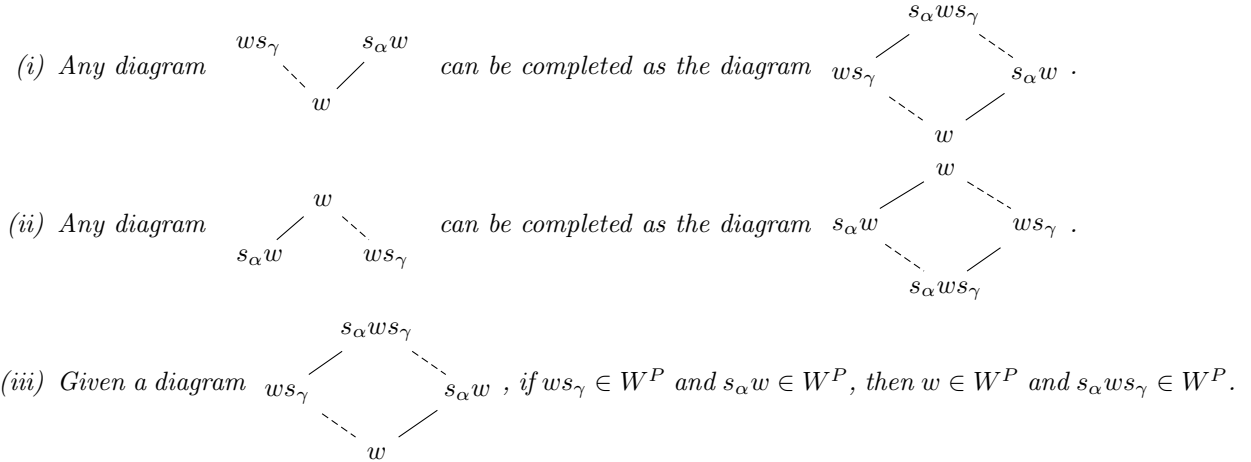
\begin{enumerate}[(i)]
 \item Any diagram
\tikzv{ws_\gamma}{s_\alpha w}{w}{}{}
can be completed as the diagram
\begin{tikzcd}[row sep=small, column sep=tiny]
& s_\alpha w s_\gamma \arrow[ld,no head] \arrow[rd,dashed,no head] \\
ws_\gamma \arrow[rd,dashed,no head] && s_\alpha w \arrow[ld,no head] \\
& w
\end{tikzcd}.

\item Any diagram
\tikzhat{w}{s_\alpha w}{w s_\gamma}{}{}
can be completed as the diagram
\begin{tikzcd}[row sep=small, column sep=tiny]
& w \arrow[ld,no head] \arrow[rd,dashed,no head] \\
s_\alpha w \arrow[rd,dashed,no head] && w s_\gamma \arrow[ld,no head] \\
& s_\alpha w s_\gamma
\end{tikzcd}.
\item
Given a diagram 
\begin{tikzcd}[row sep=small, column sep=tiny]
& s_\alpha w s_\gamma \arrow[ld,no head] \arrow[rd,dashed,no head] \\
ws_\gamma \arrow[rd,dashed,no head] && s_\alpha w \arrow[ld,no head] \\
& w
\end{tikzcd},
if $ws_\gamma \in W^P$ and $s_\alpha w \in W^P$, then $w \in W^P$ and $s_\alpha w s_\gamma \in W^P$.
\end{enumerate}
\end{lemma}
\begin{proof}
For $(i)$, since $\cX_{ws_\gamma}^B$ contains $\cX_w^B$
and not $\cX_{s_\alpha w}^B$, it cannot be
$P_\alpha$-stable. Hence  $\ell(s_\alpha w s_\beta) =\ell(w)+2$ and
the first assertion follows.

In the second assertion, assume for a contradiction that $\ell(s_\alpha ws_\gamma)=\ell(ws_\gamma)+1$. 
Since $\cX_{w}^B$ is $P_\alpha$-stable and contains
$\cX_{ws_\gamma}^B$, one must have $s_\alpha ws_\gamma=w$, a contradiction. Thus $\ell(s_\alpha
ws_\gamma)=\ell(ws_\gamma)-1$. 

To prove $(iii)$, let $s_{\alpha'} \in W_P$. Since
$s_\alpha w \in W^P$, $\ell(s_\alpha w s_{\alpha'})=\ell(w)+2$, so
$\ell(w s_{\alpha'}) = \ell(w)+1$. Thus, $w \in W^P$.
If $s_\alpha ws_\gamma\not\in W^P$ then $\dim(\cX_{s_\alpha ws_\gamma}^P)\leq \dim(\cX_{s_\alpha w}^P)=\dim(\cX_{ws_\gamma}^P)$. 
Then $\cX_{s_\alpha  ws_\gamma}^P=\cX_{s_\alpha w}^P=\cX_{ws_\gamma}^P$ and $s_\alpha w=ws_\gamma$, a contradiction.
\end{proof}

\subsection{Relative Schubert varieties and their Chow classes}

Fix $v\in W^P$. 
The usual Schubert varieties $\cX_v^B=\overline{B vB/B} \subset \GB$ and $\cX_v^P=\overline{B vP/P} \subset \GP$
are called \emph{absolute} Schubert varieties.
Consider the \emph{relative} Schubert variety
$\cY_v=G \times^B \overline{B vP/P}$. It is a closed
subvariety of $G/B\times G/P$ and, more precisely, it is the
$G$-orbit closure of $(B/B,vP/P)$ for the diagonal
$G$-action.  

The variety $\cY_v$ can be seen from different points of view:
\begin{enumerate}
\item $\cY_v=G\times^B \overline{B vP/P}$ and
  $\cY_v=G\times^P \overline{P v\inv B/B}$;
\item $\cY_v=\overline{G.(B/B,vP/P)}=\overline{G.(v\inv B/B,P/P)} \subset \GB\times\GP$;
\item $\cY_v=\{(gB/B,hP/P)\,:\, g\inv hP/P\in \overline{B vP/P}\}\subset \GB\times\GP$;
\item 
  $\cY_v=\{(gB/B,hP/P)\,:\, h\inv gB/B\in \overline{P v\inv B/B}\}\subset \GB\times\GP$.
\end{enumerate}

Recall that the Schubert classes $\sigma^w,\tau_u$ were defined in Section \ref{sec:schubert_classes}.
Think about $G/B\times G/P$ as a $G^2$-flag variety endowed
with the action of $G$ diagonally embedded in $G^2$. Then  $\cY_v$
is an orbit closure of a spherical subgroup in a flag variety. 
By \cite[Theorem~1.5]{Br:biBd}, we have in particular

\begin{prop}
\label{pro:incidence}
The Chow class of $\cY_v$ is given by:
$$
[\cY_v] = \sum
\sigma^{v'v\inv} \otimes \tau_{v'}\in A^*(G/B\times G/P),
$$
where the sum runs over the set of $v'\in W^P$ such that $v\wbole v'$.
\end{prop}

\subsection{Relative Bott-Samelson varieties}

 Resolutions of Schubert varieties can be produced as Bott-Samelson varieties. We recall this
construction in a relative context.

\begin{defi}
\label{def:bott}
Let $\bv=(\alpha_1,\dots,\alpha_\ell)$ be a sequence of elements in $\Delta$. We call (relative) Bott-Samelson 
variety the variety $\bcY_\bv$ above
$G/B$ defined by
$$\bcY_\bv := \big ( G \times P_{\alpha_1} \times \cdots \times P_{\alpha_\ell} \big ) \Big / B^{\ell+1}\,.$$
Here we consider the action of $B^{\ell+1}$ defined by
$$(b_0,b_1,\dots,b_\ell)\cdot (g,p_1,\dots,p_\ell) 
=(gb_0\inv,b_0p_1b_1\inv,b_1p_2b_2\inv,\dots,b_{\ell-1} p_\ell b_\ell^\inv)\,.$$
Let $v=s_{\alpha_1}\cdots s_{\alpha_\ell}$ and assume that $v\in W^P$ and $\ell(v)=\ell$.
As in \cite[\S2.2\parenthese{6}]{BK:book}, we may define a morphism $\omega$
from $\bcY_\bv$ to the relative Schubert variety $\cY_v$, by the
following diagram:
\begin{equation}
\label{dia:bott}
\begin{tikzpicture}[baseline=2.3cm]
\node (bcY) at (5,4) {$\bcY_\bv$};
\node[right][below=-5pt] at (5.67,4) {$:=$};
\node[right] at (6,4) {$\big ( G \times P_{\alpha_1} \times \cdots \times P_{\alpha_\ell} \big ) \Big / B^{\ell+1}$};
\node (cY) at (5,2) {$\cY_v$};
\node[right][below=-5pt] at (5.67,2) {$:=$};
\node[right] at (6,2) {$G \times^B \overline{B vP/P}$};
\node (GP) at (7,0) {$G/P$};
\node (GB) at (3,0) {$G/B$};
\draw[->] (bcY)--(cY) node[midway,right]{$\omega$};
\draw[->] (cY)--(GP) node[midway,right]{$\gamma$};
\draw[->] (cY)--(GB) node[midway,left]{$\pi$};
\draw[->] (bcY) to[bend right] node[midway,left]{$\bpi$} (GB);
\end{tikzpicture}
\end{equation}
\end{defi}

\noindent
Forgetting the last term gives a chain of morphisms
$$
\bcY_{(\alpha_1,\dots,\alpha_\ell)} \stackrel{\rho_\ell}{\longto}
\bcY_{(\alpha_1,\dots,\alpha_{\ell-1})} \stackrel{\rho_{\ell-1}}{\longto}
\cdots \longto
\bcY_{(\alpha_1)} \stackrel{\rho_{0}}{\longto}
G/B\,,
$$
which are known to be $\p^1$-bundles (see \cite[p.66]{BK:book}). 
In particular, $\bcYv$ is smooth, and more precisely, we have:

\begin{lemma}\cite[Theorem 3.4.3]{BK:book}
\label{lem:bott}
The morphism $\omega$ is a rational resolution. 
Namely, $\omega_* \cO_\bcY = \cO_\cY$ and
for $i>0$, $R^i \omega_* \cO_\bcY=0$.
\end{lemma}

\medskip

Since $\cXPv$, and thus $\cY_v$, are smooth in codimension $1$, their canonical classes are well-defined as Weil divisors, and since
they have rational singularities, we may compute them using
the following lemma:

\begin{lemma}\cite[Lemma 3.4.2]{BK:book}
\label{lem:pushkgeneral}
Let $\omega:\bcY \longto \cY$ be a rational resolution. Then $\cY$ is Cohen-Macaulay with dualizing sheaf
$\omega_* \cK_{\bcY}$.
\end{lemma}

Lemma~\ref{lem:bott} thus gives:
\begin{lemma}
 \label{lem:pushk}
 With the notation of Definition~\ref{def:bott}, we have
 $$\cK_\cYv = \omega_* {\cK}_\bcYv\mbox{ as sheaves and }
K_\cYv=\omega_*K_\bcYv \mbox{ in } \Cl(\cYv)
\,.$$
\end{lemma}

\medskip

Standard Weil divisors and line bundles are defined on those Bott-Samelson varieties as follows \cite[p.67]{BK:book}:

\begin{defi}
\label{def:divisors}
Let $k$ be such that $1 \leq k \leq \ell$ and let $\zeta$ be a character of $B$. Then:
\begin{itemize}
 \item $\bcD_k$ is the quotient by $B^{\ell+1}$ of
 $G \times P_{\alpha_1} \times \cdots \times B \times  \cdots \times P_{\alpha_\ell}$, where we replaced
 $P_{\alpha_k}$ by $B$. This is a divisor in $\bcY_\bv$, and there is an isomorphism
 $\sigma:\bcY_{(\alpha_1,\dots,\widehat{\alpha_k},\dots,\alpha_{\ell})} \longto \bcD_k \subset \bcY_\bv$.
 \item $\zeta$ and $k$ define a character $\zeta$ of $B^{\ell+1}$ by
 $\zeta(b_0,\dots,b_\ell)=\zeta(b_k)$. We define the line bundle $\bcL_k(\zeta)$ by its total space
 $\big ( G \times P_{\alpha_1} \times \cdots \times P_{\alpha_\ell} \times \C_{-\zeta} \big ) \Big / B^{\ell+1}$ and
 natural morphism to $\bcY_\bv$. 
\end{itemize}
\end{defi}

In the following lemma, we use heavier notation than the notation in the previous definition: we denote
by $\prescript{\bv}{}{\hspace{-0.8mm}\bcL_k(\zeta)}$ the bundle on $\bcY_\bv$ denoted by $\bcL_k(\zeta)$ above.

Given $\bv=(\alpha_1,\dots,\alpha_\ell)$ as above and $1\leq
k\leq\ell$, set $\bv(\leq k)=(\alpha_1,\dots,\alpha_k)$. 
The map $(g,p_1,\dots,p_\ell)\mapsto (g,p_1,\dots,p_k)$ induces a
proper surjection
$$
\rho_{\leq k}\,:\,\bcYv\longto\bcY_{\bv(\leq k)}.
$$
Note that $\rho_\ell=\rho_{\leq \ell-1}$.

\begin{lemma}
\label{lem:Lcompatible}
Let $\bv$ and $k$ be as above. Let $1\leq i\leq k$ and $\zeta\in X(T)$. 
 Then
$$\rho_{\leq k}^* (\prescript{\bv(\leq
  k)}{}{\hspace{-0.8mm}\bcL_k(\zeta)} )=
\prescript{\bv}{}{\hspace{-0.8mm}\bcL_k(\zeta)}
\qquad\mbox{and}\qquad 
\rho_{\leq k}^* (\cO(\bcD_i))=\cO(\bcD_i).$$
\end{lemma}

\noindent
This Lemma shows that there is no risk of confusion using the lighter notation $\bcL_k(\zeta)$.

\begin{lemma}
\label{lem:line-bundle}
Let $\bv$, $k$ and $\zeta$ be as above. Then, in $\Pic(\bcY_\bv)$, we have
$$
\bcL_k(\zeta) = \bcL_{k-1}(s_{\alpha_k}(\zeta)) +\scal{\zeta,\alpha_k^\vee} \cO(\bcD_k).
$$
\end{lemma}
\begin{proof}
Using Lemma~\ref{lem:Lcompatible}, it is sufficient to prove the case $k=\ell$.
This lemma is proved in \cite[Proposition 1]{demazure} with some
slightly different notation.
Therefore we give here a quick argument.
Let, as above, $\bv=(\alpha_1,\dots,\alpha_\ell)$ and $\bv'=(\alpha_1,\dots,\alpha_{\ell-1})$.
A fiber of $\rho_\ell$ is $P_{\alpha_\ell} / B \simeq \p^1$ and the bundle $\bcL_\ell(\zeta)$ restricts to
$\cO_{\p^1}(\scal{\zeta,\alpha_{\alpha_\ell}^\vee})$.
Thus
\begin{equation}
\label{eq:L}
\bcL_\ell(\zeta)=\rho_\ell^* L \otimes \cO(\scal{\zeta,\alpha_\ell^\vee}\bcD_\ell)
\end{equation}
for some line bundle $L$ on
$\bcY_\bv$ to be determined.

To determine $L$, we apply $\sigma^*$ to equation \eqref{eq:L}, where
$\sigma:\bcY_{\bv(\leq\ell-1)} \to \bcD_\ell \subset \bcY_\bv$ is as in Definition~\ref{def:divisors}. We have $\sigma^* \rho_\ell^* L=L$, and
$\sigma^* \cO(\bcD_\ell)=\cO_{\bcY_\bv}(\bcD_\ell)_{|\bcD_\ell}=\bcL_{\ell-1}(\alpha_\ell)$ since the normal bundle of $\bcD_\ell$
is $\bcL_{\ell-1}(\alpha_\ell)$. Thus we get
$$
L=\bcL_{\ell-1}(\zeta - \scal{\zeta,\alpha_\ell^\vee} \alpha_\ell) = \bcL_{\ell-1}(s_{\alpha_\ell}(\zeta))\,.
$$
\end{proof}

Since $\omega$ is proper, the pushforward $\omega_*:A^1(\bcYv)\longto
A^1(\cYv)$ is well defined. 

\begin{lemma}
\label{lem:omega}
Let $k$ be such that $1 \leq k \leq \ell$. Then, in $A^1(\cYv)$, we have
$$
\omega_* [\bcD_k] = \left \{
\begin{array}{ll}
[\cY_{s_{\alpha_1}\cdots \widehat{s_{\alpha_{k}}} \cdots s_{\alpha_\ell}}] &\mbox{ if }
s_{\alpha_1}\cdots \widehat{s_{\alpha_{k}}} \cdots s_{\alpha_\ell} \mbox{ is
                                                                     reduced ;
                                                                     
                                                                     }\\
0 &\mbox{ otherwise.}
\end{array}
\right .
$$
\end{lemma}
\begin{proof}
This follows from the fact that the restriction of $\omega$ to $\bcD_k$
is birational onto $\cY_{s_{\alpha_1}\cdots \widehat{s_{\alpha_{k}}} \cdots s_{\alpha_\ell}}$ if
the product $s_{\alpha_1}\cdots \widehat{s_{\alpha_k}} \cdots s_{\alpha_\ell}$ is reduced and contracts $\bcD_k$ otherwise.
\end{proof}

\subsection{Canonical divisor of a relative Schubert variety}
\label{sub:canonical}

We keep the notation of the previous subsection, in particular those of Definition~\ref{def:bott}.
We express the canonical classes of $\cY_v$ and $\bcY_\bv$:

\begin{prop}
\label{pro:kpi}
Assume that $\bv=(\alpha_1,\dots,\alpha_\ell)$ is a reduced expression
of $v\in W^P$ and let $\gamma_1=\alpha_\ell$, 
$\gamma_i=s_{\alpha_{\ell}}\cdots s_{\alpha_{\ell-i+2}}(\alpha_{\ell-i+1})$.
We have:
\begin{equation*}
-K_\bcYv = \bpi^* \cL_{G/B}\big (\rho+v(\rho)\big )+ \sum_{i=1}^\ell (\hc(\gamma_i)+1) [\bcD_i]
 \mbox{ in }\Cl(\bcYv)=\Pic(\bcYv).
\end{equation*}
Moreover, we have:
\begin{equation*}
-K_\cYv = \pi^* \cL_{G/B}\big (\rho+v(\rho)\big )+\sum_{v'} (\hc(\gamma)+1) [\cY_{v'}]
\mbox{ in }\Cl(\cYv)\,,
\end{equation*}
where the sum runs over the covering relations
\tikzdashedline{v}{v'}{\ggamma}
in the twisted Bruhat graph of $W^P$.
\end{prop}

\begin{proof}
An immediate induction using Lemma~\ref{lem:line-bundle} shows that
$$
\cL_k(\alpha_k)=\cL_\GB(s_{\alpha_1}\cdots s_{\alpha_k}(\alpha_k))+
\sum_{i=1}^k\scal{s_{\alpha_{i+1}}\cdots s_{\alpha_k}(\alpha_k),\alpha_i^\vee}\cO(\bcD_i).
$$
Applying
$s_{\alpha_1} \cdots s_{\alpha_i}$, we get the equality
$\scal{s_{\alpha_{i+1}}\cdots s_{\alpha_k}(\alpha_k),\alpha_i^\vee}=\scal{\beta_k,\beta_i^\vee}$,
since we have the relation $s_{\alpha_1}\cdots s_{\alpha_k}(\alpha_k)=-\beta_k$.
Finally
\begin{equation}
\label{eq:7}
\cL_k(\alpha_k)=\cL_\GB(-\beta_k)+\sum_{i=1}^k\scal{\beta_k,\beta_i^\vee}\cO(\bcD_i)\,.
\end{equation}

\noindent
The fibers of $\rho_\ell$ are isomorphic to $\p^1$ and thus the relative tangent bundle of $\rho_\ell$
is $\cL_\ell(\alpha_\ell)$. Therefore, the sequence $0 \to T \rho_\ell \to T\cY_\bv \to \rho_\ell^*(T\bcY_{\bv(\ell)}) \to 0$
reads:
\begin{equation}
  \label{eq:6}
  0\longto \cL_\ell(\alpha_\ell)\longto T \bcYv\longto
  \rho_\ell^*(T\bcY_{\bv(\ell)})\longto 0\,.
\end{equation}

Taking account that $\bcY_\emptyset=G/B$ and hence
$K_{\bcY_\emptyset}=\cL_\GB(-2\rho)$, an immediate induction gives
$$
K_{\bcYv}=\cL_\ell(-\alpha_\ell)+\rho_\ell^*(K_{\bcY_{\bv(\ell)}})=\cL_\GB(-2\rho)+\sum_{k=1}^\ell \cL_k(-\alpha_k).
$$
Injecting \eqref{eq:7}, we get
\begin{equation}
  \label{eq:8}
  \begin{array}{ll}
  K_{\bcYv}&=\cL_\GB(-2\rho+\sum_{k=1}^\ell
             \beta_k)+\sum_{k=1}^\ell\sum_{i=1}^k\scal{\beta_k,\beta_i^\vee}\cO(-\bcD_i)\\[1em]
&=\cL_\GB(-2\rho+\sum_{k=1}^\ell
             \beta_k)+\sum_{i=1}^\ell
\scal{\sum_{k=i}^\ell\beta_k,\beta_i^\vee}\;\cO(-\bcD_i)\,.
  \end{array}
\end{equation}

We already observed that $\{\beta_1,\dots,\beta_\ell\}=\Phi(v\inv)$. 
Then Lemma~\ref{lem:wrho} shows that $\sum_{k=1}^\ell
\beta_k=\rho-v(\rho)$. Moreover, Lemma~\ref{lem:sumgamma} shows that 
$\scal{\sum_{k=i}^\ell\beta_k,\beta_i^\vee}=\hc(\gamma_i)+1$,
and we get the given formula for $K_\bcYv$.

\bigskip
The formula for $K_\cYv$ then follows by applying
$\omega_*$, by Lemma \ref{lem:pushk}.  
Since $\cL_\GB(\zeta)$ is defined by a pullback,
$\omega_*(\cL_\GB(\zeta))=\pi^*(\cL_\GB(\zeta))$.
Lemma~\ref{lem:omega} allows to handle the terms $\cO(\bcD_i)$. 
We can then conclude using Lemma~\ref{lem:Bruhatorder} that realizes a
bijection between the covering relations $v'\bole v$ in the weak
Bruhat order and
$$
\{1\leq i\leq\ell \,:\omega \mbox{
  does not contract } \bcD_i\}.
$$
The given formula for $K_\cYv$ follows.
\end{proof}

\section{Proof of the main theorems}
\label{sec:branch}

\subsection{\texorpdfstring{Chow class of $Y(\bbv)$ for any $\bbv$}{subsec:Chow}}
\label{sec:Chow}

In this subsection, $\bbv=(v,\hv)$ denotes any element of $W^\bbP$ and
we consider $Y(\bbv)$ defined by \eqref{eq:defYinc}.
Observe that the projection $\cY_\hv\longto \hGP$ is a locally trivial
fibration with fiber $\overline{\hat P\hv\inv\hat B/\hat B}$.
Recall that $G/P$ is embedded in $\hGP$ and denote by $\cY_\hv^G$ the
preimage of $G/P$ in $\cY_\hv$. Then $\cY_\hv^G\simeq G\times^P \overline{\hat P\hv\inv\hat B/\hat B}$.

One can obtain $Y(\bbv)$ as a transverse  intersection:

\begin{lemma}
  \label{lem:Yinter}
Consider $\hGB \times \cYv$ and $\GB \times \hcYv$ as subvarieties
of $G/B \times \hGB \times \hGP$. Then $Y(\bbv)$ is the transverse intersection 
$$
Y(\bbv)\,=\,(\hGB \times \cYv)\, \cap \, (\GB \times \hcYv)\,.
$$
Similarly, $Y(\bbv)$ is also the transverse intersection
$(\hGB \times \cYv)\cap (\GB \times \cY_\hv^G)$ in $G/B \times \hGB \times \GP$. 
\end{lemma}

\begin{proof}
It is plain that $Y(\bbv)$ is equal to the two intersections of the lemma; we prove that these intersections are transverse.
Let $x \in \GB,\hx \in \hGB$ and $\hz \in \hGP$ such that $(x,\hx,\hz) \in (\hGB \times \cYv)\, \cap \, (\GB \times \hcYv)$.
Let $\Big ( T_{(x,\hx,\hz)}(\cY_v \times \hGB) \Big )^\bot$ denote the (conormal) space of linear forms on
$T_{(x,\hx,\hz)} (\GB \times \hGB \times \hGP)$ which vanish on the tangent space $T_{(x,\hx,\hz)}(\cY_v \times \hGB)$.
We have
$$
\begin{array}{rcl}
\Big ( T_{(x,\hx,\hz)}(\hGB \times \cY_v) \Big )^\bot & \subset & T_x^* \GB \times \{0\} \times T_\hz^* \hGP\mbox{ , and} \\
\Big ( T_{(x,\hx,\hz)}(\GB \times \cY_\hv) \Big )^\bot & \subset & \{0\} \times T_\hx^* \hGB \times T_\hz^* \hGP\ .
\end{array}
$$
It follows that the intersection
$\Big ( T_{(x,\hx,\hz)}(\cY_v \times \hGB) \Big )^\bot \cap \Big ( T_{(x,\hx,\hz)}(\cY_\hv \times \GB) \Big )^\bot$
of these conormal spaces is included in $\{0\} \times \{0\} \times T_\hz^* \hGP$, and thus reduced to $\{0\}$ since the projection
$\GB \times \cY_\hv \to \hGP$ is surjective. This means that the intersection is transverse.
\end{proof}

We now compute $[Y(\bbv)]$ in $A^*(\bbG/\bbB\times G/P)$ in the Schubert basis.

\begin{prop}
  \label{pro:classY}
With the above notation, in $A^*(\bbGB\times G/P)$, we have
\begin{equation}
[Y(\bbv)]=
  \sum
\sigma^{\bbv'\bbv\inv} \otimes
\delta^*(\tau_{\bbv'})\,,
\label{eq:tildeY}
\end{equation}
where the sum runs over the set of $\bbv'\in \bbW^\bbP$ such that $\bbv\wbole \bbv'$.
\end{prop}

\begin{proof}
By Lemma~\ref{lem:Yinter} and formula~\eqref{eq:Chowtransv}, we have

\begin{equation}
[Y(\bbv)]=[\hGB \times \cYv] \cdot [\GB \times \cY_\hv^G]\in
A^*(\bbG/\bbB \times G/P).\label{eq:10}
\end{equation}
Proposition~\ref{pro:incidence} gives a formula for $[\cY_v]\in
A^*(G/B\times G/P)$:
\begin{equation}
[\cY_{v}\times \hGB] = \sum
\sigma^{v'v\inv} \otimes 1 \otimes \tau_{v'},\label{eq:113}
\end{equation} 
where the sum runs over the set of $v'\in W^P$ such that $v\wbole v'$.

Consider now the class $[\cY_{\hv}^G]\in A^*(\hGB\times\GP)$ and 
the regular embedding  $i\,: \hGB\times \GP\longto \hGB\times \hGP$.
Let $i^*$ denote the associated Gysin map.
Since
$\cY_{\hv}^G$ is the transverse intersection of $\cY_{\hv}$ and
$\hGB\times \GP$ in $\hGB\times \hGP$, we have
\begin{equation}
  \label{eq:11}
  [\cY_{\hv}^G]=i^* [\cY_{\hv}] \in A^*(\hGB\times\GP).
\end{equation}
But $[\cY_{\hv}]$ has an expression given by
Proposition~\ref{pro:incidence} as a linear combinaison of terms
$\sigma^{\hv'\hv\inv}\otimes\tau_{\hv'}$. Using the relation 
$i^*(\sigma^{\hv'\hv\inv} \otimes \tau_{\hv'})=\sigma^{\hv'\hv\inv}\otimes\iota^*(\tau_{\hv'})$,
we get in $A^*(\hGB\times G/P)$
\begin{equation}
  \label{eq:12}
  [\cY_{\hv}^G]= \sum
\sigma^{\hv'\hv\inv}\otimes\iota^*(\tau_{\hv'})\,,
\end{equation}
and hence in $A^*(\bbG/\bbB \times G/P)$
\begin{equation}
  \label{eq:112}
  [\cY_{\hv}^G\times G/B]= \sum
1\otimes \sigma^{\hv'\hv\inv}\otimes\iota^*(\tau_{\hv'})\,,
\end{equation}
where the two last sums run over the set of $\hv'\in \hW^\hP$ such
that $\hv\wbole \hv'$.
Multiplying \eqref{eq:113} with~\eqref{eq:112} and using \eqref{equ:delta*}, one gets the formula of
the proposition.
\end{proof}

We now compute its pushforward by $\Pi$:

\begin{coro}\label{cor:PistarY}
Assume that $\ell(\bbv)=\dim(\hGP)-1$ and consider
the projection 
$\tilde\Pi\,:\, \bbG/\bbB\times G/P\longto \bbG/\bbB$.
Then, in $A^1(\bbG/\bbB)$, we have
$$
  \tilde\Pi_*([Y(\bbv)])=\sum
c(\bbv')\sigma^{s_{\bbalpha}},
$$
where the sum runs over the covering relations
\tikzlineleft{\bbv'}{\bbv}{\bbalpha} in the Bruhat graph of $\bbG/\bbP$.
\end{coro}
\begin{proof}
In general, for a class $\sum a_{\bbu,v} \sigma^\bbu  \otimes \tau_{v}$ in $A^*(\bbGB \times G/P)$, we have
\begin{equation}
\tilde \Pi_* 
\left (
\sum_{\bbu,v} a_{\bbu,v} \, \sigma^\bbu \otimes
\tau_{v}
\right)
=\sum_{\bbu} a_{\bbu,e} \, \sigma^\bbu \,.\label{eq:Pistar}
\end{equation}
But in \eqref{eq:tildeY}, the degree of $\tau_{\bbv'}$ equals
$\dim(G/P)$ if and only if $\ell(\bbv')=\dim(\hG/\hP)$ if and only if $\ell(\bbv'\bbv\inv)=1$.
Hence the terms in the sum~\eqref{eq:tildeY} that do not vanish after applying 
$\tilde \Pi_*$ are indexed by the $\bbv'=s_{\bbalpha}\bbv$ for some  simple root $\bbalpha\in\bbDelta$ 
such that $\ell(s_{\bbalpha}\bbv)=\ell(\bbv)$. 
The corresponding term is 
$\sigma^{s_\bbalpha} \otimes \delta^*(\tau_{s_{\bbalpha}\bbv})=c(s_{\bbalpha}\bbv)\sigma^{s_\bbalpha} \otimes [pt]$.
This implies the given formula.
\end{proof}

\subsection{\texorpdfstring{Canonical class of $Y(\bbv)$}{sub:canonicalY}}
\label{sub:canonicalY}

In the previous subsection, we computed the canonical class
of a relative Schubert variety, and we now express $Y(\bbv)$ as a fibered product of two relative Schubert varieties and explain
how to deduce the canonical class of $\Pi$.

Fix  reduced expressions $\bv$ and $\bhv$ of $v$ and $\hv$ respectively, and let us use
Definition~\ref{def:bott} for $v$ and $\hv$. 
We define
$\bcY_\hbv^G$ from $\bcY_\hbv$ in a way similar to what we did for $\cY_\hv^G$ from $\cY_\hv$.
We have morphisms:

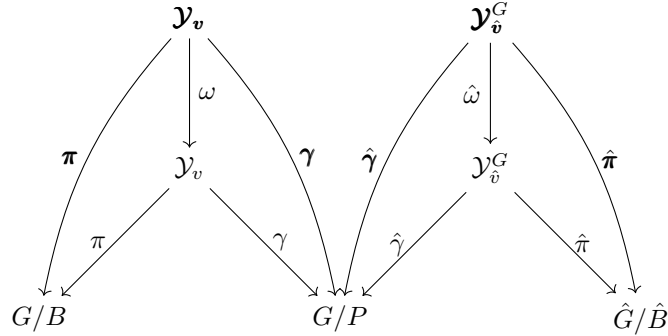
\begin{center}
\begin{tikzpicture}
\node (bcY) at (2,4) {$\bcY_\bv$};
\node (cY) at (2,2) {$\cY_v$};
\node (GP) at (4,0) {$G/P$};
\node (GB) at (0,0) {$G/B$};

\node (hbcY) at (6,4) {$\bcY_\hbv^G$};
\node (hcY) at (6,2) {$\cY_\hv^G$};
\node (hGB) at (8,0) {$\hGB$};

\draw[->] (bcY)--(cY) node[midway,right]{$\omega$};
\draw[->] (cY)--(GP) node[midway,right]{$\gamma$};
\draw[->] (cY)--(GB) node[midway,left]{$\pi$};
\draw[->] (bcY) to[bend right=15] node[midway,left]{$\bpi$} (GB);
\draw[->] (bcY) to[bend left=15] node[midway,right]{$\bgamma$} (GP);
\draw[->] (hbcY) to[bend right=15] node[midway,left]{$\hbgamma$} (GP);

\draw[->] (hbcY)--(hcY) node[midway,left]{$\homega$};
\draw[->] (hcY)--(GP) node[midway,left]{$\hgamma$};
\draw[->] (hcY)--(hGB) node[midway,right]{$\hpi$};
\draw[->] (hbcY) to[bend left=15] node[midway,right]{$\hbpi$} (hGB);
\end{tikzpicture}
\end{center}

Observe
that $Y(\bbv)$ as in \eqref{eq:defYinc} is
nothing else than $\cY_v \times_{G/P} \cY_\hv^G$ and
that $\Pi$ is the restriction of $\tilde \Pi = (\pi,\hpi)$ to $Y(\bbv)$.
We also let $\bbbv=(\bv,\hbv)$ and $\bYv = \bcY_\bv \times_{G/P} \bcY_\hbv^G$.
We have the following natural morphisms:

\begin{equation}
\label{equ:Y}
\begin{tikzpicture}[baseline=1.93cm]
\node (bcY) at (0,4) {$\bcY_\bv$};
\node (bY) at (2,4) {$\bYv$};
\node (hbcY) at (4,4) {$\bcY_\hbv^G$};
\node (cY) at (0,2) {$\cY_v$};
\node (Y) at (2,2) {$Y(\bbv)$};
\node (hcY) at (4,2) {$\cY_\hv^G$};
\node (GP) at (0,0) {$G/P$};
\node (GB) at (4,0) {$G/B \times \hGB$};
\node (GBP) at (-2.5,2) {$G/B \times G/P$};
\node (hGBP) at (6.5,2) {$\hGB \times G/P$};

\draw[->] (bY) -- (bcY) node[midway,above]{$\bepsilon$};
\draw[->] (bY) -- (hbcY) node[midway,above]{$\hbepsilon$};
\draw[->] (bcY) -- (cY) node[midway,left]{$\omega$};
\draw[->] (bY) -- (Y) node[midway,right]{$\Omega$};
\draw[->] (hbcY) -- (hcY) node[midway,right]{$\homega$};
\draw[->] (Y) -- (cY) node[midway,above]{$\epsilon$};
\draw[->] (Y) -- (hcY) node[midway,above]{$\hepsilon$};
\draw[->] (Y) -- (GP) node[midway,left]{$\Gamma$};
\draw[->] (Y) -- (GB) node[midway,right]{$\Pi$};

\draw[left hook->] (cY) -- (GBP);
\draw[right hook->] (hcY) -- (hGBP);
\end{tikzpicture}
\end{equation}

The first properties of the morphisms appearing in this diagram we
need are:

\begin{lemma}\label{lem:epsproper}
  With the above notation, 
  \begin{enumerate}
  \item $\epsilon$ is a locally trivial fibration with fiber
    $\overline{\hP\hv\inv \hB/\hB}=\cX_{w_{0,\hP}\hv}^\hB$.
\item $\hepsilon$ is a locally trivial fibration with fiber
    $\overline{P v\inv B/B}=\cX_{w_{0,P}v}^B$.
 \item $\bepsilon$ and $\hbepsilon$  are  locally trivial fibrations.
  \end{enumerate}
In particular, these four morphisms are flat and induce pullbacks
between Chow rings.
\end{lemma}

\begin{proof}
We already observed when we defined $\hcYvG$ at the beginning of Section \ref{sec:Chow}
that $\hcYvG$ is a locally trivial fibration over $\GP$ with fiber
$\overline{\hP\hv\inv \hB/\hB}$. Since $Y(\bbv)=\cYv\times_\GP\hcYvG$,
the first assertion follows. The second one works similarly. 

By the Bruhat decomposition, the morphism $G\longto G/P$ is locally trivial in Zariski
topology. Hence the maps $\hbcYvG\longto G/P$ and  $\bcYv\longto G/P$
are locally trivial fibrations. The last assertion follows.
\end{proof}

We can now express the canonical bundle of $Y(\bbv)$ in terms of the canonical bundles of $\cY_v$ and $\cY_\hv$:

\begin{prop}\label{prop:KY}
We have $\displaystyle K_{Y(\bbv)} = \epsilon^* K_{\cY_v} + \hepsilon^*
K_{\cY_\hv^G} - \Gamma^* K_{G/P}$ in $\Cl(Y(\bbv))$.
\end{prop}

\begin{proof}
Let us first prove this result at the level of Bott-Samelson varieties. We have by definition of $\bYv$ a fibered square of smooth varieties:
\begin{center}
\begin{tikzcd}
\bYv \arrow{r}{\bepsilon \times \hbepsilon} \arrow[labels=left]{d}{\Gamma \circ \Omega} \arrow[phantom]{rd}{\square} &
\bcYv \times \hbcYvG\arrow{d}{\bgamma \times \hbgamma} \\
G/P \arrow{r}{\Delta} & G/P \times \GP
\end{tikzcd}
\end{center}

\noindent
Let $\bbby=(\by,\hby) \in \bYv$ and let $u=\Gamma \circ \Omega(\bbby) \in G/P$. We thus have a fibered square of tangent spaces:
\begin{center}
\begin{tikzcd}
T_\bbby \bYv \arrow{r} \arrow{d} \arrow[phantom]{rd}{\square} &
T_\by \bcYv \oplus T_\hby \bcY_\hbv^G \arrow{d} \\
T_u G/P \arrow{r}{d\Delta} & T_uG/P \oplus T_u\GP
\end{tikzcd}
\end{center}

The bottom horizontal map being injective, so is the top horizontal
one. The right vertical map being surjective,
$(T_\by \bcY \oplus T_\hby \bcY_\hbv^G)/T_\bbby \bYv$ identifies with $( T_uG/P \oplus T_u\GP)/T_u
G/P$.
Moreover, on $\bYv$, the bundle $( TG/P \oplus T\GP)/TG/P$ is isomorphic to $\Omega^*\Gamma^*TG/P$.
We get the short exact sequence on $\bYv$
$$
0 \longto T_\bYv \longto T_{\bcY_\bv} \oplus T_{\bcY_\hbv^G} \longto \Omega^* \Gamma^* T_\GP \longto 0\,.
$$

\noindent
From this we deduce that
\begin{equation}
\label{equ:KbY}
\begin{array}{rcl}
K_\bYv & = & \big ( K_{\bcY_\bv} + K_{\bcY_\hbv^G} \big )_{|\bYv} - \Omega^* \Gamma^* K_\GP \\
& = & \bepsilon^* K_{\bcY_\bv} + \hbepsilon^*K_{\bcY_\hbv^G} - \Omega^* \Gamma^* K_\GP \ .
\end{array}
\end{equation}

We now apply $\Omega_*$. Lemma~\ref{lem:bott} yields $\Omega_* \cO_\bYv = \cO_{Y(\bbv)}$ and $\Omega_* K_\bYv = K_{Y(\bbv)}$.
From the first point and projection formula, we deduce in $\Cl(Y(\bbv))$
\begin{equation}
 \label{equ:KGP}
 \Omega_* \, \Omega^* \, \Gamma^* \, K_\GP \, = \, \Gamma^* \, K_\GP \,.
\end{equation}

Moreover, we have a fibered product
  \begin{tikzcd}
    \bYv \arrow{r}{(\bepsilon,\hbepsilon)} \arrow[labels=left]{d}{\Omega}
    \arrow[phantom]{rd}{\square} &
    \bcYv\times\hbcYvG \arrow{d}{\omega} \\
    Y(\bbv) \arrow[labels=below]{r}{(\epsilon,\hepsilon)} & \cYv\times\cY_\hv^G
  \end{tikzcd}
  . By the flat base change $(\epsilon,\hepsilon)$ and \cite[Proposition
  III.9.3]{hartshorne}, we deduce that
  \begin{equation}
    \label{equ:basechange}
    \Omega_* (\bepsilon^*,\hbepsilon^*) ({\mathcal K}_{\bcYv}\otimes{\mathcal K}_{\hbcYvG}) = (\epsilon^*,\hepsilon^*) \omega_* ({\mathcal K}_\bcYv\otimes{\mathcal K}_{\hbcYvG}) \,.
  \end{equation}

But, by \cite[Theorem~5.10]{KollarMori}, $\omega_* {\mathcal K}_\bcYv=
{\mathcal K}_\cYv$ and $\omega_* {\mathcal K}_{\hbcYvG}=
{\mathcal K}_{\cY_\hv^G}$. We deduce that 
 \begin{equation}
    \label{equ:Omegaepsstar}
    \Omega_* (\bepsilon^* K_\bcYv+ \hbepsilon^* K_{\hbcYvG})=  \epsilon^* K_\cYv+\hepsilon^* K_{\cY_\hv^G}.
  \end{equation}
\noindent
By the equalities \eqref{equ:KbY}, \eqref{equ:KGP}, \eqref{equ:basechange}, and \eqref{equ:Omegaepsstar}, we get the result of the proposition.
\end{proof}

For later use, we also compute $\Pi_*\circ\Gamma^*$ (recall from \eqref{equ:chi} that
$\euler{}\,$ is the coefficient of the point class):

\begin{lemma}
  \label{lemm:PiGamma}
The map $\Pi_*\circ\Gamma^*\,:A^*(G/P)\longto A^*(\bbG/\bbB)$ is
determined by
$$
\Pi_*\circ\Gamma^* (\xi) =
\sum
\euler{\GP}(\xi \cdot \delta^*(\tau_{\bbv'})) \sigma^{\bbv'\bbv\inv}
\quad \forall \xi \in A^*(G/P),
$$
where the sum runs over the set of $\bbv'\in \bbW^\bbP$ such
that $\bbv\wbole \bbv'$.
If $\xi \in A^d(\GP)$, then all the non-vanishing terms in this summand satisfy $\ell(\bbu)=d$.
\end{lemma}
\begin{proof}
Let $\xi \in A^*(\GP)$ and consider the following commutative diagram:

\begin{equation}
\label{eq:def_i}
\begin{tikzcd}
    Y(\bbv)\arrow{r}{i} \arrow{rd}{\Gamma} & \bbG/\bbB\times\GP \arrow{d}{p}\\
    &G/P.
\end{tikzcd}
\end{equation}

\noindent

By the projection formula, we have
$$
i_* \Gamma^* \xi = i_*i^*p^* \xi=p^*\xi \cup i_*[Y(\bbv)]\,.
$$
Proposition~\ref{pro:classY} gives $i_*[Y(\bbv)]$, and since
$p^* \xi=1\otimes \xi$, one gets
\begin{equation}
  \label{eq:318}
  i_* \Gamma^* \xi = \sum_{\bbu\bbv \stackrel{\bbP}{\to}_\bbu \bbv}
\sigma^{\bbv'\bbv\inv} \otimes
(\xi \cdot \delta^*(\tau_{\bbv'}))\,,
\end{equation}
where the sum runs over the set of $\bbv'\in \bbW^\bbP$ such
that $\bbv\wbole \bbv'$.
Since $\Pi_* \Gamma^* \xi = \tilde \Pi_* i_* \Gamma^* \xi$, we deduce the formula of the proposition.

Assume now that $\xi$ is homogeneous of degree $d$.
As in the proof of Corollary~\ref{cor:PistarY}, the terms to keep when
applying $\tilde \Pi_*$ to \eqref{eq:318} satisfy $\ell(\bbu)=d$.
\end{proof}

\subsection{Computation of the class of the branch divisor}
\label{sub:branch}

We compute $[B_\Pi] =\Pi_* K_\Pi$ using the already done computation of $K_{Y(\bbv)}$.

\begin{proof}[Proof of Theorem~\ref{theo:branch}:]
recall that $Y(\bbv) \subset \bbGB \times \GP$ and that $\Pi$ resp. $\Gamma$ is the projection of $Y(\bbv)$
on $\bbGB$ resp. $\GP$.
Consider also as before
$$
\tilde\Pi\,: \bbGB \times G/P \longto \bbGB\,,
$$
the restriction to $Y(\bbv)$ of which is $\Pi$.

By the definition of $K_\Pi$ in Section~\ref{sec:canonical} and Proposition~\ref{prop:KY}, we have 
\begin{equation}
  \label{eq:9}
  K_\Pi = \epsilon^* K_{\cY_v} + \hepsilon^*
K_{\cY_\hv^G} - \Gamma^* K_{G/P}-\Pi^*(K_\GB + K_\hGB).
\end{equation}

\medskip
Consider first $K_{\cY_\hv^G}$ in $A^1(\cY_\hv^G)$.
Let $j\,: \cY_\hv^G\longto \hcYv$ denote the inclusion.
As the pullback of the inclusion $\iota\,:G/P\longto\hGP$, it is a regular embedding and we can consider the associated Gysin map
$j^*$. 
At any point $y$ of  $\cY_\hv^G$ with projection $\bar y\in
\GP$ we have the exact sequence:
$$
0\longto T_y \cY_\hv^G\longto T_y \hcYv\longto T_{\bar y}\hGP/T_{\bar
  y}\GP\longto 0.
$$
We deduce that in $\Cl(\cY_\hv^G)$, we have
\begin{equation}
  \label{eq:314}
  K_{\hcYvG} = j^* K_{\hcYv} - \hgamma^*((K_\hGP)_{|\GP}) + \hgamma^*(K_\GP).
\end{equation}
Using the formulas $-K_\hGP=\cL_\hGP(2\rho^\hL)$, 
$-K_{\bbG/\bbB}=\cL_\bbGB(2\rho_\bbG)$, and Proposition~\ref{pro:kpi},
injecting \eqref{eq:314} in \eqref{eq:9}, one gets \vspace{-.5em}

\begin{equation}
\label{eq:315}
\begin{array}{c}
 K_\Pi = A + B - C - D \ \mbox{ with} \ \\[0.8em]
  \left \{
  \begin{array}{l}
   A = \Pi^* c_1(\cL_{\bbG/\bbB}(\rho_\bbG-\bbv(\rho_\bbG)))\\[0.2em]
   B = \Gamma^* c_1(\cL_\GP(2\rho^\hL_{|T})) \\[0.2em]
   C = \sum_\gamma (\hc(\gamma)+1 ) \epsilon^*([\cY_{v'}]) \\[0.2em]
   D = \sum_\hgamma (\hc(\hgamma)+1 ) \hepsilon^*\circ j^*([\cY_{\hv'}])\,,
  \end{array}
  \right .
\end{array}
\end{equation}

\noindent
where the sums $C$ and $D$ run respectively over the pictures
\tikzdashedline{v}{v'}{\ggamma} and \tikzdashedline{\hv}{\hv'}{\ghgamma}
in the twisted Bruhat graphs of $W^P$ and $W^\hP$ respectively.

\bigskip
We now pushforward by $\Pi$ each summand of \eqref{eq:315}.
First of all, Lemma~\ref{lemm:PiGamma} gives
$$
\Pi_* B=
\sum
\left (c_1(\cL_\GP(2\rho^\hL_{|T})) \cdot \delta^*(\tau_{s_{\bbalpha}\bbv}) \right )
\sigma^{s_{\bbalpha}}\,,
$$
where the sum runs over the covering relations
\tikzlineleft{s_\bbalpha\bbv}{\bbv}{\bbalpha} in the Bruhat graph of $\bbG/\bbP$.
This is the first term of the formula of Theorem \ref{theo:branch}. 

\medskip

For $\bbzeta\in X(\bbT)$, the projection formula gives 
$
  \Pi_*\circ\Pi^*c_1(\cL_{\bbG/\bbB}(\bbzeta))=\deg(\Pi)c_1(\cL_{\bbG/\bbB}(\bbzeta)),
$
and the Chevalley formula gives
$c_1(\cL_{\bbG/\bbB}(\bbzeta))=\sum_{\bbalpha\in\bbDelta}\scal{\bbzeta,\bbalpha^\vee}\sigma^{s_\bbalpha}$.
Thus
\begin{equation}
  \label{eq:316}
  \Pi_* A =
  c(\bbv) \sum_{\bbalpha\in\bbDelta} \scal{\rho_\bbG-\bbv(\rho_\bbG), \bbalpha^\vee} \sigma^{s_\bbalpha}\,.
\end{equation}

Consider first a term $\bbalpha\in\bbDelta$ such  that
$s_\bbalpha\bbv\in \bbv W_\bbP$. Then $s_\bbalpha\bbv =\bbv
s_\bbgamma$ with $\bbgamma=\bbv\inv(\bbalpha)\in\Phi(\bbL)$. 
Since  $\bbv\in W^\bbP$,
$\ell(\bbv)+1 \geq \ell(s_\bbalpha\bbv)=\ell(\bbv s_\bbgamma)=\ell(\bbv)+\ell(s_\bbgamma)$. Then
$\bbgamma$ is a simple root of $\bbL$. In particular,  
$\scal{\bbv(\rho_\bbG) , \bbalpha^\vee}=1$ and the
corresponding term in the sum~\eqref{eq:316} vanishes.

Assume now that $\ell(s_\bbalpha \bbv)=\ell(\bbv)+1$ and $s_\bbalpha \bbv \in W^\bbP$. 
Then $\scal{\bbv(\rho_\bbG)-\rho_\bbG, \bbalpha^\vee}=\hc(\bbgamma)-1$
for $\bbgamma = \bbv \inv (\bbalpha)$.
Assume finally that $\ell(s_\bbalpha \bbv)=\ell(\bbv)-1$ and $s_\bbalpha \bbv \in W^\bbP$.
Then, similarly, we have
$\scal{\rho_\bbG-\bbv(\rho_\bbG), \bbalpha^\vee}=\hc(\bbgamma)+1$,
for $\bbgamma=-\bbv \inv \bbalpha$. Summing up, we get
\begin{equation}
\label{eq:317}
\Pi_* \, A =
  \sum
-c(\bbv)(\hc(\bbgamma)-1) \sigma^{s_\bbalpha}
  + \sum
c(\bbv)(\hc(\bbgamma)+1) \sigma^{s_\bbalpha} \,,
\end{equation}
where the two sums run respectively over the pictures
\tikzdashedline{s_\bbalpha \bbv}{\bbv}{\gbbgamma} and \tikzdashedline {\bbv}{s_\bbalpha \bbv}{\gbbgamma}
in the twisted Bruhat graphs of $\bbW^\bbP$.

\bigskip
We now compute $\Pi_*\circ \epsilon^* ([\cY_{v'}])$, where $v'=
s_\alpha v\in W^P$, $\alpha\in\Delta$  and $\ell(v)=\ell(v')+1$.
Set $\bbv'=(v',\hv)$.
Set $Y'=Y(\bbv')=\cY_{v'}\times_\GP \cY_\hv^G$, so that we have
\begin{equation}
  \label{eq:4}
  \epsilon^* [\cY_{v'}]=[Y']\in\Cl(Y(\bbv)). 
\end{equation}
Consider the class $i_* [Y']$ of $Y'$ in $A^*(\bbGB\times G/P)$ (recall the diagram \eqref{eq:def_i}). Then
\begin{equation}
\label{eq:5}
\Pi_*(\epsilon^* [\cY_{v'}])=\tilde\Pi_*i_*[Y'].
\end{equation}
We apply Corollary~\ref{cor:PistarY} to $\bbv'$ to get
\begin{equation}
\tilde\Pi_* i_* [Y']
=\sum
c(s_{\bbalpha}\bbv')\sigma^{s_{\bbalpha}}\,,
\label{eq:118}
\end{equation}
where the sum runs over the covering relation
\tikzlineleft{s_\bbalpha\bbv'}{\bbv'}{\bbalpha}.
Then
\begin{equation}
\label{eq:335}
\Pi_* C = \sum (\hc(\gamma)+1) c(\bbv'') \sigma^{s_{\bbalpha}},
\end{equation}
where the sum runs over the subgraphs
\tikzv{\bbv}{\bbv''}{\bbv'}{\ggamma}{\bbalpha}
in the Bruhat graph
of $\bbG/\bbP$
with $\bbalpha\in\bbDelta$ and $\gamma\in \Phi^+(G)$.

\bigskip
We now compute a term
$
\Pi_*\circ \hepsilon^*\circ
j^*(\cY_{\hv'})
$  where $\hv'=
s_\halpha \hv\in \hW^\hP$, $\halpha\in\hat\Delta$  and $\ell(\hv)=\ell(\hv')+1$.
Since the intersection $\cY_{\hv'}\cap
\cY_\hv^G=\cY_{\hv'}^G$ is transverse, we have
$j^*([\cY_{\hv'}])=[\cY_{\hv'}^G]$.
Then 
$$
\Pi_*\circ \hepsilon^*\circ
j^*([\cY_{\hv'}])=\Pi_*([\cYv\times_\GP \cY_{\hv'}])
$$
is given by Corollary~\ref{cor:PistarY}. 
One gets that

\begin{equation}
  \label{eq:336}
  \Pi_* D = \sum (\hc(\hgamma)+1) c(\bbv'') \sigma^{s_{\bbalpha}},
\end{equation}
where the sum runs over the subgraphs
\tikzv{\bbv}{\bbv''}{\bbv'}{\ghgamma}{\bbalpha}
in the Bruhat graph of $\bbG/\bbP$
with $\bbalpha\in\bbDelta$ and $\hat\gamma\in \Phi^+(\hG)$.
Putting together \eqref{eq:335} and \eqref{eq:336}, we get

\begin{equation}
\label{eq:pi*CD}
\Pi_* (C+D) = \sum (\hc(\bbgamma)+1) c(\bbv'') \sigma^{s_{\bbalpha}},
\end{equation}
where the sum runs over the subgraphs
\tikzv{\bbv}{\bbv''}{\bbv'}{\gbbgamma}{\bbalpha}
in the Bruhat graph of $\bbG/\bbP$.

\bigskip
By Lemma~\ref{lemm:bruhat2}, the terms in \eqref{eq:pi*CD}
with $\bbv\neq \bbv''$ are exactly the terms in the second sum of
Theorem \ref{theo:branch}.
Let us now consider a term in \eqref{eq:pi*CD} with $\bbv = \bbv''$, meaning that
$s_\bbalpha \bbv = \bbv s_\bbgamma$. If $\ell(s_\bbalpha \bbv) < \ell(\bbv)$, then
the coefficient of $\sigma^{s_\bbalpha}$ in $\Pi_* A$ and in $\Pi_*(C+D)$ is $c(\bbv)(h(\bbgamma)+1)$, by
\eqref{eq:316} and \eqref{eq:pi*CD}.
So these terms cancel, and accordingly Theorem \ref{theo:branch} states that this coefficient vanishes.
If $\ell(s_\bbalpha \bbv) > \ell(\bbv)$, then
the coefficient of $\sigma^{s_\bbalpha}$ in $\Pi_* A$ is $-c(\bbv)(\hc(\bbgamma)-1)$
and in $\Pi_*(C+D)$ it is $0$. This is also compatible with Theorem \ref{theo:branch}.
\end{proof}

\subsection{Proof of our reduction formula}
\label{sub:proof-reduction}

In this section, we prove Theorem~\ref{main-theo}.

\subsubsection{The birational case}

We recall from \cite{ressayre-birational} the main construction in the proof of Theorem~\ref{theo1},
because this will allow introducing useful notation.
Recall that $\bbG=G\times\hG$, $\bbv=(v,\hv)$ and set
$\bbL=L\times\hL$ and $\bbzeta=(\zeta,\hzeta)$.
Consider the variety $X=\bbG/\bbB$ endowed with the diagonal
$G$-action.
Let $\cL \to X$ be the line bundle
defined by $\bbzeta$.
By Borel-Weil theorem, we have

\begin{equation}
m_{G \subset \hG}(\zeta,\hzeta)=\dim(\Ho(X,\cL)^G).\label{eq:mGgeom}
\end{equation}

Let $C \subset X$ be the irreducible component of the fixed point set
$X^\tau$ of $\tau$ in $X$ containing the point $x_0:=\bbv\inv\bbB/\bbB$.
Let $B_\bbL=\bbB \cap \bbL$.
Observe that the stabilizer of $x_0$ in $\bbL$ is $B_\bbL$, so that 
the variety $C$ is isomorphic to
$\bbL/B_\bbL$. Moreover, $\cL_{|C}$ is isomorphic to the line bundle defined by
$\bbv\inv\bbzeta$. In particular

\begin{equation}
m_{L \subset
  \hL}(v\inv\zeta,\hv\inv\hzeta)=\dim(\Ho(C,\cL)^L).
\label{eq:mLgeom}
\end{equation}

\noindent
The Bia\l{}ynicki-Birula cell $C^+$ is
$$
C^+ = \{ x \in X \, | \, \lim_{t \to 0} \tau(t)x \in C \}=\bbP \bbv\inv \bbB/\bbB\,.
$$
Consider 
\begin{equation}
\label{equa-Y}
\begin{tikzcd}
G \times^P \overline{C^+} \arrow{d} \arrow{r}{\Pi} & X \\
G/\Pl
\end{tikzcd}
\end{equation}
Observe that $Y(\bbv)=G \times^P \overline{C^+}$ endowed with  $\Pi\,:\, Y(\bbv)\longto X$ is nothing but the incidence
variety defined in \eqref{eq:defYinc}. 

The proof in \cite{ressayre-birational} goes on proving the following, which does not need the hypothesis $c(\bbv)=1$:
\begin{theo}
\label{theo-geometrie}
Assume that $v\inv(\zeta)_{|S} + \hv\inv(\hzeta)_{|S}$ is trivial.
With the notation of \eqref{equa-Y}, there is a natural isomorphism
$$\Ho(Y(\bbv),\Pi^* \cL)^G \simeq \Ho(C,\cL_{|C})^L\,.$$
\end{theo}

In \cite{ressayre-birational}, to prove Theorem~\ref{theo-geometrie},
we first prove that $\Ho(C,\cL_{|C})^L=\Ho(C^+,\cL_{|C^+})^P$ using
the action of $\tau$. Since $\Ho(C^+,\cL_{|C^+})^P=\Ho(G \times^P
C^+,\cL)^G$, it remains to prove that each rational $G$-invariant on
$Y(\bbv)$ that is regular on the open set $G
\times^P C^+$ extends to $Y(\bbv)$. By normality of the Schubert
varieties, it is sufficient to prove that $\sigma$ has no pole along
the boundary. This last step is made by quite explicit local
computation.\\

The fiber of $\Pi$ over any general point in $x$ is isomorphic to the
intersection of two subvarieties of $G/P$ of Chow class
$\tau_v$ and $\iota^* \tau_\hv$: therefore the assumption $c(\bbv)=1$ implies that $\Pi$ is birational and gives further
$\Ho(Y(\bbv),\Pi^* \cL) \simeq \Ho(X,\cL)$, from which Theorem~\ref{theo1} follows.

\subsubsection{The degree 2 case}

We  now deal with the case $c(\bbv)=2$. Then $\Pi$ is no longer birational, but it is generically finite
of degree $2$, and we apply the projection formula as in Section
\ref{sec:genfinite}.
We denote by $B_\Pi \subset X$ the branch divisor.
Then, in the Schubert basis of $X=\bbG/\bbB$, 
$[B_\Pi]$ expands as

\begin{equation}
\label{eq:1defni}
[B_\Pi] = \sum_\alpha n_\alpha \ \sigma^{s_\alpha} \otimes 1 + \sum_\halpha n_\halpha \ 1 \otimes \sigma^{s_\halpha} \in
A^*(\bbG/\bbB)\,,
\end{equation}
for some well defined integers $n_\alpha$ and $n_\halpha$.
These coefficients are described by Theorem~\ref{theo:branch}.

We proved in Lemma~\ref{degree2} that $\frac 1 2[B_\Pi]\in \Pic(X)$.  
Recall that  $\theta$ and $\hat\theta$ in $X(T)$ and  $X(\hT)$ are
defined by 
\begin{equation}
  \label{eq:2bis}
  \theta=\sum_\alpha\frac{n_\alpha} 2\varpi_\alpha\qquad\hat\theta =\sum_\halpha\frac{n_\halpha}2{\varpi}_\halpha\,,
\end{equation}
in such a way that $\frac 1 2[B_\Pi]$ is the line bundle on
$X=G/B\times \hG/\hB$ associated to the character $(\theta,\hat\theta)$
of $B\times\hB$.

\begin{proof}[Proof of Theorem~\ref{main-theo}]
Observe that $C$ can be seen as a subvariety of
$Y(\bbv)=G\times^P\overline{C^+}$ via the closed immersion $x \mapsto [e,x]$,
and that $\Pi^*(\cL)_{|C}$ is equal to $\cL_{|C}$. 
In particular equality~\eqref{eq:mLgeom} can be rewriten as
$$
m_{L \subset \hL}(v\inv\zeta,\hv\inv\hzeta)=\dim
\Ho(C,\Pi^*(\cL)_{|C})^L.
$$
Now, by Theorem~\ref{theo-geometrie}, there is an isomorphism
$$
\Ho(C,\Pi^*(\cL)_{|C})^L \simeq \Ho(Y(\bbv),\Pi^*\cL)^G.
$$
By the projection formula and Proposition~\ref{prop:mor2},
$$\Ho(Y(\bbv),\Pi^*\cL)=\Ho(X,\cL \otimes \Pi_* \cO_{Y(\bbv)} )=\Ho (X,\cL) \oplus
\Ho \Big (X,\cL \big (-\frac12[B_\Pi] \big ) \Big ).$$
But by the Borel-Weil Theorem
$$
\Ho \Big (X,\cL \big (-\frac12[B_\Pi] \big ) \Big )=V_{\zeta-\theta} \otimes V_{\hzeta-\hat\theta}.
$$
Considering the dimensions of $G$-invariant subspaces, we get the theorem.
\end{proof}

\section{Properties of the branch divisor class}

\subsection{The Belkale-Kumar product}
\label{sec:BKprod} 

Fix a one-parameter subgroup $\tau$ of $T$ and set $P=P(\tau)$.
For $w\in W^P$, recall that $\tau^w=\tau_{w^\vee}$ is a class of degree $\ell(w)$, see Section \ref{sec:notations}.
For $w\in W^P$, define the BK degree of $\tau^w\in A^*(G/P)$ to be
$$
\BKdeg(\tau^w):=\scal{w\inv(\rho)-\rho,\tau}.
$$ 
Let $w_1,w_2$ and $w_3$ in $W^P$. By \cite[Proposition 17]{bk}, if
$c(w_1^\vee,w_2^\vee,w_3)\neq 0$ that is if $\tau^{w_3}$ appears in
the product $\tau^{w_1} \cdot \tau^{w_2}$ then 
\begin{equation}
  \label{eq:13}
  \BKdeg(\tau^{w_3})\leq \BKdeg(\tau^{w_1})+\BKdeg(\tau^{w_2}).
\end{equation}
In other words the BK degree filters the Chow ring. 
Let $\odot$ denotes the associated graded product on $A^*(G/P)$.
For later use, note that $\odot$ coincides with the usual  product
when $G/P$ is cominuscule (See \cite{bk}).

\medskip

Assume now that $G$ is embedded in $\hG$ and set $\hP=\hP(\tau)$. 
Recall that $\iota\,:\,G/P\longto\hG/\hP$ denotes the inclusion. 
In \cite{RR} a graded ring morphism is defined
$$
\iota^\odot\,:\,(A^*(\hG/\hP),\odot)\longto (A^*(G/P),\odot).
$$
As in \eqref{eq:def_c}, define the integers $c^\odot(v,\hv)$ by
\begin{equation}
  \label{eq:14}
  \forall \hv\in W^\hP\,,\ \iota^\odot(\tau_\hv)=\sum_{v\in W^P}c^\odot(v,\hv)\tau^{v}.
\end{equation}

Fix $\bbv=(v,\hv)\in W^\bbP$ such that $\ell(\bbv)=\dim(\hG/\hP)$.
Consider now the incidence variety $\Pi\,:\,Y(\bbv)=G\times^P\overline{C^+}\longto
\bbG/\bbB$ as above. 
Note that $c(\bbv)$ is not zero if and only if $\Pi$ is dominant. By
equivariance this means that there exists a point $x\in C^+$ such the
tangent map of $\Pi$ at $[e:x]$ is invertible. 

By \cite[Proposition 2.3]{RR} and the definition of $\iota^\odot$, 
the coefficient $c^\odot(\bbv)$ is not zero if and only if there exists
$y_0\in C$ such that  the tangent map of $\Pi$ at $[e:y_0]$ is invertible.
Moreover, if $c^\odot(\bbv) \neq 0$, then $\tau(\C^*)$ acts trivially on the line bundle $(K_\Pi)_{|C}$, so that
for $y\in C^+$ with $y_0=\lim_{t\to 0}\tau(t)y$, the
tangent map of $\Pi$ at $[e:y]$ is invertible if and only if the
tangent map of $\Pi$ at $[e:y_0]$ is.
Writing $y_0=(l,\hl).\bbv^\inv \bbB/\bbB$,
the tangent map of $\Pi$ at $[e:y_0]$ is invertible if and only if the
intersection
$$
\iota(l v^\inv B v P/P)\cap (\hat l \hv^\inv \hB \hv \hP/\hP)
$$
is transverse at $P/P$ in $\hG/\hP$.
In \cite{bk}, this condition is called Levi-movability. To sum up, we obtain:

\begin{lemma}
\label{lem:levi-movability}
In the setting of Theorem~\ref{theo:branch},
assume that 
$$
\tau_v\odot \iota^\odot(\tau_{\hv})=k[pt]\,,
$$
with $k>0$.
Then, for general $(l,\hl)$ in $\bbL$, the intersection
$\iota(lv\inv BvP/P) \cap  ( \hl \hv\inv \hB \hv \hP/\hP)$
is transverse at $P/P$ and the tangent map of $\Pi$ at $[e:y_0]$ is invertible, with
$y_0=(l,\hl).\bbv^\inv \bbB/\bbB$.
\end{lemma}

\subsection{The branch divisor class in the Levi-movable case}
\label{sec:eigenconeBPi}

Consider the support 
$$
\Gamma(G,\hG)=\{(\zeta,\hzeta)\in X(T)^+\times X(\hT)^+\;:\; m_{G \subset
  \hG}(\zeta,\hzeta)\neq 0\}
$$
of the branching multiplicities.
It is known (see \cite{elash}) to be a finitely generated semigroup and
hence it generates a convex polyhedral cone $\Gamma_\QQ(G,\hG) \subset
(X(T) \oplus X(\hT)) \otimes \QQ$, that we call Horn cone since the introduction. 
Let $S$ be the center of the Levi subgroup $L$ of $G$.
Given $\bbv=(v,\hv)\in W^\bbP$ such that $c(\bbv)\neq 0$, the points
$(\zeta,\hzeta)$ in $\Gamma_\QQ(G,\hG)$ satisfy
$$
\scal{v\inv(\zeta),\tau} + \scal{\hv\inv(\hzeta),\tau}\leq 0
$$ 
for any dominant (for $G$) one-parameter subgroup $\tau$ of $S$.
In particular, the set of pairs $(\zeta,\hzeta)\in (X(T)\times X(\hT))\cap\Gamma_\QQ(G,\hG)$ 
such that $v\inv(\zeta)_{|S} + \hv\inv(\hzeta)_{|S}$ is trivial is a face
$\cF(\bbv)$ of the Horn cone. 

\medskip

Under the Levi-movability assumption, the branch divisor class $\frac 1 2[B_\Pi]$ is
the minimal element of $\cF(\bbv)$ that does not satisfy the
conclusion of Theorem~\ref{theo1}:

\begin{prop}
  \label{pro:posC}
In the setting of Theorem~\ref{theo:branch},
assume that $c^\odot(\bbv)>0$. Then $\frac 12 [B_\Pi] \in \cF(\bbv)$, and we have,
for any $(\zeta,\hzeta)\in \cF(\bbv)$, 
$$
m_{G \subset \hG}(\zeta,\hzeta)\leq m_{L \subset \hL}(v\inv(\zeta),\hv\inv(\hzeta))\,.
$$

\noindent
Assume more specifically that $c^\odot(\bbv)=2$.
Then the set of pairs $(\zeta,\hzeta)\in \cF(\bbv)$ where this
inequality is strict is exactly the set
$$\frac 12[B_\Pi]+(\cF(\bbv)\cap\Gamma(G,\hG))\,.$$
\end{prop}
\begin{proof}
First, let us prove that $[B_\Pi]$ belongs to the face $\cF(\bbv)$.
If $c^\odot(\bbv)=1$, then $\frac 12 [B_\Pi]=0$ and there is nothing to prove. We assume $c^\odot(\bbv)>1$.
Recall that $\cO(B_\Pi)=\cL(\theta,\hat\theta) \in \Pic(\bbG/\bbB)$. We have $(\theta,\htheta) \in \cF(\bbv)$ if and only if
$v\inv(\theta)_{|S} + \hv\inv(\htheta)_{|S}$ is trivial
that is if $S$ acts trivialy on $\cO(B_\Pi)_{|C}$. 

We pick $(l,\hl) \in \bbL$ as provided by Lemma~\ref{lem:levi-movability} and we set
$y_0=(l,\hl)\bbv\inv \bbB/\bbB$.
In particular $[e:y_0]\in Y(\bbv)=G\times^P \overline{C^+}$ is an isolated
point in the fiber $\Pi\inv(y_0)$.
By transversality and degree assumption, the fiber
$\Pi\inv(y_0)$ cannot be reduced to this point. 
In particular, this fiber is not connected.
By Proposition~\ref{pro:supportBPi}, $y_0$ does not belong to the support
$\Supp(B_\Pi)$ of $B_\Pi$.

Since $\cO(B_\Pi)$ is the line bundle on the smooth variety $X$
associated to some divisor with $\Supp(B_\Pi)$ as support, $\cO(B_\Pi)$
has a canonical section $\sigma$ that does not vanish at $y_0$.
By unicity up to  scalar multiplication of this section and $G$-invariance of $\Supp(B_\Pi)$, it has to be an eigenvector for the action
of $G$. Since $G$ is semi-simple this section is a $G$-invariant element of $\oH^0(X,\cO(B_\Pi))$, meaning
that $\sigma$ is a $G$-equivariant section of $\cO(B_\Pi)$.

But $S$ fixes $y_0$ and hence acts on the fiber $\cO(B_\Pi)_{y_0}$. This
action has to fix $\sigma(y_0)$ and hence is trivial.

\medskip

Let $f \in \oH^0(\bbG/\bbB,\cL)^G$ that restricts to $0$ on $C$. Then $f$ vanishes on $C^+$ and thus $\Pi^* f$ vanishes on $Y(\bbv)$.
Thus, $\oH^0(\bbG/\bbB,\cL)^G$ injects into $\oH^0(C,\cL_{|C})^L$, and the inequality of the proposition follows.

Assuming that $c^\odot(\bbv)=2$, Theorem~\ref{main-theo} shows that the inequality will be strict if and
only if $m_{G \subset \hG}(\zeta-\theta,\hzeta-\hat\theta) \neq 0$, which means
that $(\zeta,\hzeta)$ belongs to $(\theta,\htheta)+(\cF(\bbv)\cap\Gamma(G,\hG))$.
\end{proof}

\begin{coro}
\label{cor:alt}
In the setting of Theorem~\ref{main-theo}, assume that $c^\odot(v,\hv)=2$. Then
$$
m_{G \subset \hG}(\zeta,\hzeta) =\sum_{k\geq 0}(-1)^k
m_{L \subset \hL}(v\inv(\zeta-k\theta),\hv\inv(\hzeta-k\hat\theta))\,.
$$
\end{coro}

\begin{proof}
Direct induction from Theorem~\ref{main-theo}. Indeed
Proposition~\ref{pro:posC} allows to apply the theorem to each weight
$(\zeta-k\theta,\hzeta-k\hat\theta)$.
\end{proof}

\section{Examples}
\label{sec:examples}

\subsection{Type A}
\label{sec:notation-example}

In this section, we assume that $G$ has type $A$.
In type $A_{n-1}$, it is more convenient to work with representations
of $\GL_n(\CC)$ instead of $\SL_n(\CC)$.
For $\GL_n(\CC)$, we have
$$
\Lambda^+_{\GL_n}:=X(T)^+=\{\sum_{i=1}^n\lambda_i\varepsilon_i\,:\,\lambda_1\geq\cdots\geq\lambda_n\in\ZZ\},
$$
with notation as in \cite{bourb}. 
Given $\lambda\in X(T)^+$, set
$\bar\lambda=\sum_{i=1}^{n-1}(\lambda_i-\lambda_n)\varepsilon_i$. 
Then $\bar\lambda$ is a partition and the $\GL_n(\CC)$-representation $V_\lambda$
of highest weight $\lambda$ is isomorphic to the $\SL_n(\CC)$-representation
of highest weight $\bar \lambda$
 as an $\SL_n(\CC)$-representation.
For $\lambda,\mu,\nu\in \Lambda^+_{\GL_n}$, set
$m_{\GL_n}(\lambda,\mu,\nu)=\dim(V_\lambda\otimes V_\mu\otimes V_\nu)^{\GL_n(\CC)}$.
Because of the action of the center of $\GL_n(\CC)$, we have
$$
m_{\GL_n}(\lambda,\mu,\nu)\neq 0\Longrightarrow |\lambda|+|\mu|+|\nu|=0,
$$
where $|\lambda|=\sum_i\lambda_i$ and $|\mu|$ and $|\nu|$ are defined similarly.
Set $\nu^*=\sum_i-\nu_{n+1-i}\varepsilon_i$ so that $V_{\nu^*}$
is the dual representation of $V_\nu$. If $\lambda$, $\mu$
and $\nu^*$ are partitions (that is $\lambda_n,\mu_n,\nu^*_n\geq 0$)
then

\begin{equation}
m_{\GL_n}(\lambda,\mu,\nu)=c_{\lambda,\mu}^{\nu^*}
\label{eq:mlr}
\end{equation}
is the LR coefficient.
The Horn semigroup of $\GL_n(\CC)$ is
$$
\Gamma(\GL_n(\CC))=\{(\lambda,\mu,\nu)\in (\Lambda^+_{\GL_n})^3\,:\, m_{\GL_n}(\lambda,\mu,\nu)\neq 0\}.
$$
Let $1\leq r\leq n-1$ and consider $\Gr(r,n)=G/P$. Given a subset $I\subset \{1,\ldots,n\}$ 
with $r$ elements, define $v_I\in W^P$ by $v_I(\{1,\dots,r\})=I$ and set
$I^\vee = \{n+1-i\, | \, i \in I\}$. Recall from \eqref{eq:15} in the introduction that $c_{I,J}^K$
denotes the structure constants of the Chow ring of the Grassmannian $\Gr(r,n)$ equipped
with its Schubert basis. 
For any subsets $I,J,K\subset \{1,\ldots,n\}$ 
with $r$ elements, we have

\begin{equation}
c(v_I,v_J,v_K)=c_{I^\vee,J^\vee}^K.
\label{eq:cvlr}
\end{equation}

Recall also that given a partition $\lambda$ and a subset $I$, $\lambda_I$ is the
partition whose parts are $\lambda_i$ with $i \in I$.
Assume that $c(v_I,v_J,v_K)\neq 0$.
The associated Horn inequality~\eqref{eq:166} expressed in terms of
$\Gamma(\GL_n(\CC))$ becomes

\begin{equation}
(\lambda,\mu,\nu)\in \Gamma(\GL_n(\CC)) \Longrightarrow
|\lambda_I|+|\mu_J|+|\nu_K|\leq 0.
\label{eq:HornGL}
\end{equation}

\noindent{\bf Two interpretations.} Fix $0<r<n$.
Let $I=\{i_1<\cdots<i_r\}$ be a subset of $\{1,\dots,n\}$ with $r$
elements. Let 
$$
\lambda(I)=i_r-r\geq\cdots\geq i_1-1
$$
be the associated partition. 

An amazing property of the LR coefficients is that they are also the
multiplicities of the tensor product decomposition for the linear
groups. Namely, 
$c_{I,J}^K$ is also the multiplicity
$c_{\lambda(I),\lambda(J)}^{\lambda(K)}$ of $V_{\SL_r}(\lambda(K))$
in $V_{\SL_r}(\lambda(I))\otimes V_{\SL_r}(\lambda(J))$.
Hence

\begin{equation}
  \label{eq:Lessieur}
  c_{I,J}^K=c_{\lambda(I),\lambda(J)}^{\lambda(K)}.
\end{equation}

In this section, we only consider the case $G \subset \hG=G \times G$,
for a given classical group $G$ of type $A$ or $B$. In particular, we use Notation~\ref{nota:tensor_product_case}.
In each case, we have a vector space $V$ endowed with a basis
$(e_1,\dots,e_n)$ and eventually a symmetric bilinear form.
In type $B$, we choose the bilinear form
$\beta$ such that $\beta(e_i,e_j) =1$ if and only if $j=n+1-i$, and equals $0$ otherwise. With this convention,
the set of diagonal resp. triangular matrices in $G$
is a maximal torus $T$ resp. Borel subgroup $B$.

The classical Grassmannians $\Gr=G/P$ occur.
The Schubert varieties in $\Gr$ are indexed by subsets $I$ of $\{1,\dots,n\}$ in the following way:
$$
X_I=\overline{B \cdot V_I} \mbox{ with } V_I:=\Span(e_i\,:\,i\in I)\;, 
\quad{\rm and}\quad
\tau_I=[X_I]\in A^{\codim(X_I)}(\Gr).
$$
This also gives a bijection between a set of subsets of
$\{1,\dots,n\}$ and $W^P$.
Moreover, set $\overline \imath =n+1-i$.
and $I^\vee = \{\overline \imath \, | \, i \in I\}$.
This operation is the Poincaré duality: $\tau^{I}=\tau_{I^\vee}$.

In type A, we also consider the two steps flag manifolds
$\Fl(p,q;n)=\{F_1\subset F_2\subset V\,:\,\dim(F_1)=p,\;\dim(F_2)=q\}$.
The Schubert varieties are then indexed by ``flags'' $(I_p \subset I_q)$ of subsets of
$\{1,\dots,n\}$: $X_{(I_p \subset I_q)} = \overline {B \cdot (V_{I_p},V_{I_q})}$.

\subsection{\texorpdfstring{A detailed example in $\Gr(3,6)$}{sub:g36}}
\label{sub:g36}

In this subsection, $G=\SL_6$, $G/P=\Gr(3,6)$, and $\hG=G^2$.
We let $I=J=K=\{2,4,6\}\in W^P$.
Set $\bbv=(v_I,v_J,v_K)$.
The inversion set $\Phi(v_I)$ of $v_I$ is depicted with
black nodes:

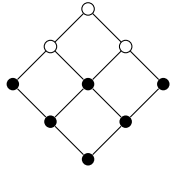
\begin{center}
\begin{tikzpicture}[scale=0.5, transform shape]
\tikzstyle{vide}=[circle,draw]
\tikzstyle{plein}=[circle,fill]
\tikzstyle{ligne}=[] 

\node (1) at (2,4) [vide] {};
\node (2) at (1,3) [vide] {};
\node (3) at (3,3) [vide] {};
\node (4) at (0,2) [plein] {};
\node (5) at (2,2) [plein] {};
\node (6) at (4,2) [plein] {};
\node (7) at (1,1) [plein] {};
\node (8) at (3,1) [plein] {};
\node (9) at (2,0) [plein] {};

\draw [ligne] (1) -- (2);
\draw [ligne] (1) -- (3);
\draw [ligne] (2) -- (4);
\draw [ligne] (2) -- (5);
\draw [ligne] (3) -- (5);
\draw [ligne] (3) -- (6);
\draw [ligne] (4) -- (7);
\draw [ligne] (5) -- (7);
\draw [ligne] (5) -- (8);
\draw [ligne] (6) -- (8);
\draw [ligne] (7) -- (9);
\draw [ligne] (8) -- (9);
\end{tikzpicture}
\end{center}

Since $\lambda(I^\vee)=\lambda(J^\vee)=21$ and $\lambda(K)=321$, we
have $c(\bbv)=c_{21,\ 21}^{321}=2$, by \eqref{eq:cvlr} and \eqref{eq:Lessieur}.
We describe geometrically the divisors $R_\Pi$ and $B_\Pi$ in this example, but we start with
 Lemma~\ref{lem:bitriangle} about configurations of triangles in the plane.

\begin{defi}
\label{def:bitriangle}
Given a $3$-dimensional vector space $E$, a sextuple $(A_1,A_2,A_3,B_1,B_2,B_3)$ of points in the projective plane $\p E$
is called a bitriangle if the three sets of points
$\{A_i,A_j,B_k\}$ for $i,j,k$ distinct in $\{1,2,3\}$ are colinear. In the generic situation,
$A_1,A_2,A_3$ define a triangle and the three points $B_1,B_2,B_3$ are on the sides of this triangle.

The variety of all bitriangles in $\p E$ will be denoted by $\cT(E)$. Given two vector spaces $E$ and $E'$ of dimension $3$, and
two bitriangles
\begin{equation}
T=(A_1,A_2,A_3,B_1,B_2,B_3) \in \cT(E) \ ,\
T'=(A'_1,A'_2,A'_3,B'_1,B'_2,B'_3) \in \cT(E')\,,\label{eq:defTriangles}
\end{equation}
 a morphism
$u:T \longto T'$ is a linear map
$u:E \longto E'$ such that for all $i$, $u(A_i) \subset A'_i$ and $u(B_i) \subset B'_i$ (here $A_i,B_i$ resp.
$A'_i,B'_i$ are considered as
$1$-dimensional subspaces of $E$ resp. $E'$). Note that the set of
morphisms of bitriangles is a subspace of the vector space of linear maps $E \to E'$.
\end{defi}

Observe that $\cT(E)$ is irreducible of dimension $9$ and that the modality of the action of $\PSL(E)$ is one.

\begin{lemma}
\label{lem:bitriangle}
There is exactly one divisor $\cD$ in $\cT(E) \times \cT(E')$ such that for all pairs $(T,T')$ in $\cD$, the space of morphisms
from $T$ to $T'$ is not reduced to $\{0\}$. It may be described as the divisor of pairs of isomorphic bitriangles.
\end{lemma}
\begin{proof}
First of all, the variety $\cD$ of pairs of isomorphic bitriangles is a divisor,
since the modality of the action of $\PSL(E)$ on the space of bitriangles is $1$. Moreover,
tautologically, if $T$ and $T'$ are isomorphic, then the space of morphisms $T \longto T'$ is not reduced to $\{0\}$.

On the other
hand, let $u:T \longto T'$ be a morphism of bitriangles, with $T$
and $T'$ as in \eqref{eq:defTriangles}.
If $u$ has rank $3$, then $T$ and $T'$ are isomorphic and $(T,T') \in \cD$.

Assume first that $T$ is generic.
If $u$ has rank $2$, then at least $5$ of the six points $A_1,A_2,A_3,B_1,B_2,B_3$ have a well defined image in $\p E'$ and all these
images belong to the line in $\p E'$ defined by the image of $u$. It follows that $5$ of the $6$ points
$A'_1,A'_2,A'_3,B'_1,B'_2,B'_3$ are on a line, and the set of such bitriangles has codimension $2$ in $\cT(E')$.
If $u$ has rank $1$, similarly, we don't find a divisor in $\cT(E')$ since $3$ of the six points $A'_1,A'_2,A'_3,B'_1,B'_2,B'_3$
must be equal in this case.
The case where $T$ is degenerate is similar.
\end{proof}

\begin{prop}
\label{pro:rx}
The hypersurface $B_\Pi$ is irreducible, it is equal to the variety $\Delta$ defined below in (\ref{equa:Delta}), and 
$$
\cO_{(G/B)^3}(B_\Pi) = \cL_{(G/B)^3}(2\varpi_2+2\varpi_4,2\varpi_2+2\varpi_4,2\varpi_2+2\varpi_4)\,.
$$
The hypersurface $R_\Pi$ is also irreducible, it is the preimage of $B_\Pi$ by $\Pi$,
and it also has the description given below in (\ref{equ:rpi-bitriangle}).
\end{prop}
\begin{proof}
Let $U \subset (G/B)^3$ be the open subset of triples $(\drapeau{1}{\bullet},\drapeau{2}{\bullet},\drapeau{3}{\bullet})$
such that $\C^6=\drapeau{1}{2} \oplus \drapeau{2}{2} \oplus \drapeau{3}{2}$,
$\drapeau{1}{4} \cap \drapeau{2}{4} \cap \drapeau{3}{4} = \{0\}$, and
$\drapeau{i}{4} \cap \drapeau{j}{2} = \{0\}$ for $i \neq j$.
We first investigate the intersection $B_\Pi \cap U$. In fact, let
$(\drapeau{1}{\bullet},\drapeau{2}{\bullet},\drapeau{3}{\bullet}) \in U$: the fiber
$\Pi^{-1}(\drapeau{1}{\bullet},\drapeau{2}{\bullet},\drapeau{3}{\bullet})$ can be explicitly described as follows.

\medskip

Let $V_3 \in \Gr(3,6)$. 
Then $(V_3;\drapeau{1}{\bullet},\drapeau{2}{\bullet},\drapeau{3}{\bullet})$ is an element of $Y(\bbv)$ if and only if
\begin{equation}
\label{equ:g36}
\forall i \in \{1,2,3\}\ , \ \dim V_3 \cap \drapeau{i}{2} \geq 1 \mbox{ and } \dim V_3 \cap \drapeau{i}{4} \geq 2 \,.
\end{equation}
From $\dim V_3 \cap \drapeau{i}{2} \geq 1$ and the genercity assumption, we deduce that $V_3$ can be written as
\begin{equation}
\label{equa:v3}
V_3=L_1 \oplus L_2 \oplus L_3\,,
\end{equation}
with $L_i \subset \drapeau{i}{2}$ a subspace of dimension $1$.
Moreover, under our genericity assumption, we have an isomorphism
$\drapeau{i}{2} \simeq \C^6/\drapeau{k}{4}$ if $i \neq k$. If $i,j,k$ are distinct, we get isomorphisms
$\drapeau{i}{2} \simeq \C^6 / \drapeau{k}{4}$ and $\C^6 / \drapeau{k}{4} \simeq \drapeau{j}{2}$, and we denote
$\varphi_{i,j}:\drapeau{i}{2} \longto \drapeau{j}{2}$ the isomorphism obtained by composition.

Observe that $L_1 \oplus L_2 \oplus L_3$ will meet $\drapeau{k}{4}$ in dimension $2$ if and only if
the subspace of $\C^6$ generated by $L_i,L_j$ and $\drapeau{k}{4}$ has dimension $5$, which is equivalent to the equality
$L_j=\varphi_{i,j}(L_i)$.

Thus, $V_3=L_1 \oplus L_2 \oplus L_3$ satisfies (\ref{equ:g36}) if and only if for $i \neq j$, $L_j=\varphi_{i,j}(L_i)$.
In other words, $L_2=\varphi_{1,2}(L_1), L_3=\varphi_{2,3}(L_2)$, and
$L_1 \subset \drapeau{1}{2}$ is an eigenline of $\varphi_{3,1} \circ \varphi_{2,3} \circ \varphi_{1,2}$.

\medskip

By Proposition \ref{pro:supportBPi}, $B_\Pi \cap U$ is the set of triples
$(\drapeau{1}{\bullet},\drapeau{2}{\bullet},\drapeau{3}{\bullet})$
such that $\varphi_{3,1} \circ \varphi_{2,3} \circ \varphi_{1,2}$ has
only one eigenvalue or is the identity.

We denote by $\Delta$ this divisor:

\begin{equation}
\label{equa:Delta}
\Delta = \{ (\drapeau{1}{\bullet},\drapeau{2}{\bullet},\drapeau{3}{\bullet}) \, | \,
\varphi_{3,1} \circ \varphi_{2,3} \circ \varphi_{1,2} \mbox{ has only one eigenvalue} \}
\end{equation}

\medskip

We now show that there is only one divisor in the ramification of $\Pi$. To show this, we observe that for a general element
$(V_3,\drapeau{1}{\bullet},\drapeau{2}{\bullet},\drapeau{3}{\bullet})$ of a divisor of $Y(\bbv)$,
$\dim(V_3 \cap \drapeau{i}{2})=1$ and $\dim(V_3 \cap \drapeau{i}{4})=2$ for $i \in \{1,2,3\}$.
Therefore, this defines a bitriangle
$T$ in $\p V_3$, with vertices $A_i=V_3 \cap \drapeau{i}{2}$ and edges $E_i=V_3 \cap \drapeau{i}{4}$, and this defines
also a bitriangle $T'$ in $\C^6/V_3$, with vertices
$A'_i=p(V_3\drapeau{i}{2})$ and edges $E'_i=p(\drapeau{i}{4})$, where $p:\C^6 \longto \C^6/V_3$ denotes the projection.

In this way, we factorize the morphism $Y(\bbv) \to \Gr(3,6)$ as a composition of rational maps
$Y(\bbv) \stackrel{f}{\dasharrow} \cT(E) \times_{\Gr(3,6)} \cT(Q) \longto \Gr(3,6)$, where $E$ and $Q$ denote the tautological bundles on
$\Gr(3,6)$, and $\cT(E)$ and $\cT(Q)$ denote the relative varieties of bitriangles therein. The rational map $f$ is defined
in codimension $1$ and equidimensional. It follows that $f(R_\Pi)$ is defined and at least
a divisor in $\cT(E) \times \cT(Q)$.

\smallskip

On the other hand, an element $u \in T_{V_3} \Gr(3,6) \simeq \Hom(V_3,\C^6/V_3)$ belongs to the tangent space to the Schubert
variety defined by $\drapeau{i}{\bullet}$ if and only if $u(A_i) \subset A'_i$ and $u(E_i) \subset E'_i$. Since this must hold
for all $i$ in $\{1,2,3\}$, $u$ defines a morphism of bitriangles (recall Definition \ref{def:bitriangle}). It thus follows from
Lemma \ref{lem:bitriangle} that the image of $R_\Pi$ by $f$ is the divisor of isomorphic bitriangles in $\cT(E) \times \cT(Q)$:
\begin{equation}
\label{equ:rpi-bitriangle}
R_\Pi = \{ y \in Y : f(y) \in \cT(E) \times \cT(Q) \mbox{ is a pair of isomorphic bitriangles.} \}
\end{equation}

It follows that $R_\Pi$ is irreducible, equal to $\Pi^{-1}(B_\Pi)$, and that $B_\Pi=\Delta$.
Moreover, the given formula for $\cO_{(G/B)^3}(B_\Pi)$ can be proved using Corollary \ref{coro:branch} and can also
be seen more geometrically as a consequence of \eqref{equa:Delta}.
\end{proof}

We now check directly Theorem \ref{main-theo} in this case. First we give a full
description of the face in the Horn cone defined by $\bbv$.
Recall that $I=J=K=\{2,4,6\}$.

\begin{lemma}
\label{lem:face}
The span of the face $\cF_{IJK}$ of $\Gamma(\GL_6)$ defined by the equation
$|\lambda_I|+|\mu_J|+|\nu_K|=0$ has  equations
\begin{equation}
\label{equ:face}
\left \{
\begin{array}{l}
|\lambda|+|\mu|+|\nu|=0\,;\\
\lambda_1=\lambda_2\,, \lambda_3=\lambda_4 \,, \mbox{and } \lambda_5=\lambda_6\,; \\
\mu_1=\mu_2\,, \mu_3=\mu_4 \,, \mbox{and } \mu_5=\mu_6 \,;\\
\nu_1=\nu_2\,, \nu_3=\nu_4 \,, \mbox{and } \nu_5=\nu_6 \,.
\end{array}
\right .
\end{equation}
\end{lemma}
\begin{proof}
The equalities (\ref{equ:face}) clearly imply
$|\lambda_I|+|\mu_J|+|\nu_K|=0$.
 On the other hand, if $|\lambda_I|+|\mu_J|+|\nu_K|=0$,
we get
$$
0=2|\lambda_I|+2|\mu_J|+2|\nu_K|\leq |\lambda|+|\mu|+|\nu|=0.
$$

This proves that the middle inequality is an equality. So (\ref{equ:face}) hold.
\end{proof}

\begin{rema}
\label{rem:facevsfacereg}
  The face $\cF_{IJK}$ in Lemma~\ref{lem:face} is contained in no face corresponding to some
LR coefficient equal to one.
\end{rema}

\begin{proof}
  Let $(I',J',K')$ of the same cardinality $r$ such that
  $c(I',J',K')\neq 0$ and the elements of $\cF_{IJK}$ satisfy
  $|\lambda_{I'}|+|\mu_{J'}|+|\nu_{K'}|=0$.
The description of the span of $\cF_{IJK}$ in Lemma~\ref{lem:face}
implies that the linear form $(\lambda,\mu,\nu)\mapsto
|\lambda_{I'}|+|\mu_{J'}|+|\nu_{K'}|$ belongs to the span of the 10
linear forms vanished in equations~\eqref{equ:face}: there are rational numbers 
$n,d,a_1,a_3,a_5,b_1,b_3,b_5,c_1,c_3,c_5$ such that
$$
\begin{array}[]{l@{\,}l@{\,}l}
  n(|\lambda_{I'}|+|\mu_{J'}|+|\nu_{K'}|)=&&d(|\lambda|+|\mu|+|\nu|)\\
  &+&a_1(\lambda_1-\lambda_2)+a_3(\lambda_3-\lambda_4)+a_5(\lambda_5-\lambda_6)\\
  &+&b_1(\mu_1-\mu_2)+b_3(\mu_3-\mu_4)+b_5(\mu_5-\mu_6)\\
  &+&c_1(\nu_1-\nu_2)+c_3(\nu_3-\nu_4)+c_5(\nu_5-\nu_6)
\end{array}
$$
for any $\lambda,\mu,\nu$. We may assume $n=6$.
Choosing $\lambda=\mu=\nu=1^6$, one gets $d=r$.
Since $0<r<6$, the linear form
$6(|\lambda_{I'}|+|\mu_{J'}|+|\nu_{K'}|)-r(|\lambda|+|\mu|+|\nu|)$ has positive coefficient on $\lambda_i$
if and only if $i \in I'$, and negative coefficient otherwise.
It follows that $I'$, $J'$,
and $K'$ all intersect $\{1,2\}$, $\{3,4\}$ and $\{5,6\}$ in exactly one element.
In particular $r=3$.  
Now, the degree condition implied by the non-vanishing of the Schubert
coefficient $c(I',J',K')$ implies that $I'=J'=K'=\{2,4,6\}$. In particular, $c_{I'J'}^{K'}=2$.
\end{proof}

By the remark, Theorem~\ref{introtheo:birational} cannot be applied to
the points in this face. 
However, applying our reduction formula,
we compute explicitly those coefficients:

\begin{lemma}
Let $\lambda,\mu,\nu\in \Lambda^+_{\GL_3}$ such that  $|\lambda|+|\mu|+|\nu|=0$.
Then $m_{\GL_3}(\lambda,\mu,\nu)$ is equal to the number of integers in the interval
$$
[\max(\mu_1-\lambda_2,\mu_2,-\nu_3-\lambda_1,\mu_1+\nu_1,-\nu_2-\lambda_2,\mu_1+\mu_2+\nu_2),
\min(\mu_1,-\nu_3-\lambda_2,\mu_1+\mu_2+\nu_1)]\,.
$$

Moreover, set $\lambda^2=(\lambda_1,\lambda_1,\lambda_2,\lambda_2,\lambda_3,\lambda_3)$,
and define similarly $\mu^2$ and $\nu^2$. Then
$$
m_{\GL_6}(\lambda^2,\mu^2,\nu^2) = \frac{m_{\GL_3}(\lambda,\mu,\nu)(m_{\GL_3} (\lambda,\mu,\nu)+1)}{2}\,.
$$
\end{lemma}
\begin{proof}
The first equality follows from \cite{kt:hive}. In fact, it is explained in \cite[Proposition 9]{PR:unexpected}
that $c_{\lambda,\mu}^{\nu}$ is the number of integers in the interval
$$
[\max(\mu_1-\lambda_2,\mu_2,\nu_1-\lambda_1,\mu_1-\nu_3,\nu_2-\lambda_2,\mu_1+\mu_2-\nu_2),
\min(\mu_1,\nu_1-\lambda_2,\mu_1+\mu_2-\nu_3)]\,.
$$
Using the fact that $V_{(\nu_1,\nu_2,\nu_3)}^* \simeq V_{(-\nu_3,-\nu_2,-\nu_1)}$, we get our formula.
We want to compute $m_{\GL_6}(\lambda^2,\mu^2,\nu^2)$ using Corollary \ref{cor:alt} with $\bbv=(v_I,v_J,v_K)$ and
$I=J=K=\{2,4,6\}$. By Proposition \ref{pro:rx}, one half of the branch divisor corresponds to the triple of partitions
$((2,2,1,1,0,0),(2,2,1,1,0,0),(-1,-1,-2,-2,-3,-3)).$
Thus, a term
$m_{L}(v\inv(\gamma-k\theta),\hv\inv(\hgamma-k\hat\theta))$ in Corollary \ref{cor:alt}
is equal in our context to
$$\Big (m_{\GL_3}(\lambda-k(210),\mu-k(210),\nu+k(123)) \Big )^2.$$

Observe that by the first point, $m_{\GL_3} (\lambda-k(210),\mu-k(210),\nu+k(123))$ is $m_{\GL_3}(\lambda,\mu,\nu)-k$ or $0$.
Indeed, when $\lambda,\mu,\nu$ are replaced respectively by $\lambda-(210),\mu-(210),\nu+(123)$, all the integers that appear
in the $\max$ decrease by $1$ and all the integers that appear in the $\min$ decrease by $2$.
Corollary~\ref{cor:alt} therefore gives
$$
m_{\GL_6} (\lambda^2,\mu^2,\nu^2) = \sum_{k=0}^{m_{\GL_3}(\lambda,\mu,\nu)} (-1)^{m_{\GL_3} (\lambda,\mu,\nu)+k}\, k^2\,.
$$
This is equal to $\frac{m_{\GL_3} (\lambda,\mu,\nu)(m_{\GL_3} (\lambda,\mu,\nu)+1)}{2}$, as stated.
\end{proof}

\begin{rema}
\label{rem:polynomial}
As we can observe here, the \lr coefficients are close to being polynomials in the coefficients of the weights. This is a general phenomenon.
Kostant \cite[Theorem 6.2]{kostant} proves that the multiplicity of a weight in an irreducible representation is a partition
function. Steinberg \cite{steinberg} deduces that the multiplicity of an irreducible submodule in the tensor
product of two representations is again a partition function. Rassart \cite[Theorem 4.1]{rassart} deduces that the
Horn cone can be subdivided into subcones where the \lr coefficients are polynomial in the three weights (this holds in type $A$;
in general the coefficients are only quasi-polynomial).
\end{rema}

\subsection{A bigger example}
\label{sub:big}
In the following example, the face of $\Gamma(G)$ associated to the
LR-coefficient equal to two is contained in a face associated to some
LR-coefficient equal to one. Hence, one can apply both Theorems~\ref{introtheo:birational}
and \ref{main-theo-A}. 
The conclusion is that, even if
Theorem~\ref{introtheo:birational} apply, our Theorem~ \ref{main-theo-A} is useful.
\\

Here $n=10$, $r=6$ and $G/P=\Gr(6,10)$. Set $I=\{2,4,5,6,9,10\}$, $J=\{3,4,6,7,9,10\}$
and $K=\{2, 4, 5, 7, 8, 10\}$. 
Set $\bbv=(v_I,v_J,v_K)$ in such a way that $c(\bbv)=c_{\lambda(I^\vee), \lambda(J^\vee)}^{\lambda(K)}= c_{3222,2211}^{433221}=2$.
Then, applying Theorem~\ref{theo:branch}, one gets:

\begin{equation}
\label{eq:expleBPi2}
\frac 12 [B_\Pi] = \cL_{(G/B)^3}(\varpi_2+2\varpi_6, \varpi_4+\varpi_7, \varpi_2+\varpi_5+\varpi_8).
\end{equation}

First, we prove a statement analogous to Lemma \ref{lem:face}:
\begin{lemma}
\label{lem:face2}
The span of the face $\cF_{IJK}$  defined by the equation $|\lambda_I|+|\mu_J|+|\nu_{K}|=0$ in
$\Gamma(\GL_{10})$ has  equations
\begin{equation}
\label{equ:face2}
\left \{
\begin{array}{l}
|\lambda|+|\mu|+|\nu|=0\,;\\
(\lambda_5+\lambda_6) + (\mu_7+\mu_{10}) + (\nu_5+\nu_8)=0\,;\\
\lambda_1=\lambda_2\,, \lambda_3=\lambda_4 \,, \mbox{and } \lambda_7=\lambda_{8}=\lambda_9=\lambda_{10}\,;\\
\mu_1=\mu_2=\mu_3=\mu_4 \,, \mu_5=\mu_6 \,, \mbox{and } \mu_8=\mu_9 \,;\\
\nu_1=\nu_2\,, \nu_3=\nu_4 \,, \nu_6=\nu_7 \,, \mbox{and } \nu_9=\nu_{10} \,.
\end{array}
\right .
\end{equation}
\end{lemma}
\begin{proof}
Let $F$ be the linear space defined by the  equations~\eqref{equ:face2}. We observe that the following triples $(\lambda,\mu,\nu)$ span $F$ and satisfy
$m_{\GL_{10}}(\lambda,\mu,\nu)=1$:
$$
\begin{array}{cccc}
&(1^{10},0,(-1)^{10}) & (1^{10},0,(-1)^{10}) \\
(1^6,(-1)^{10},1^4)&(1^5,(-1)^{10},1^5)&(1^2,(-1)^{10},1^8)&(1^6,(-1)^4(-2)^6,1^{10})\\
(1^6,(-1)^9-2,1^5)&(1^6, (-1)^7(-2)^3,1^7)&(1^6, (-1)^6(-2)^4,1^8)&(1^5, (-1)^7(-2)^3,1^8)\\
(1^6,(-1)^7(-2)^3,2^21^3)&(1^6, (-1)^4(-2)^6,2^21^6)&(1^6, (-1)^4(-2)^3(-3)^3,2^51^3)&(2^41^2, (-2)^7(-3)^3,2^51^3)
\end{array}
$$

Conversely, let $(\lambda,\mu,\nu)$ satisfying
$|\lambda_I|+|\mu_J|+|\nu_{K}|=0$ and $m_{\GL_{10}}(\lambda,\mu,\nu)\neq 0$.
To have nicer formulas, we write $a$ instead of $10$.
We have
\begin{equation}
\label{equ:ineq}
\begin{array}{rclcl}
0 & = & 2|\lambda_I|+2|\mu_J|+2|\nu_{K}| &\leq&
                                           |\lambda_{\{12\}}|+|\lambda_{\{34\}}|+2\lambda_5+2\lambda_6+|\lambda_{\{789a\}}|\\
&&&+&|\mu_{\{1234\}}|+|\mu_{\{56\}}|+2\mu_7+|\mu_{\{89\}}|+2\mu_{a}\\
&&&+&|\nu_{\{12\}}|+|\nu_{\{34\}}|+2\nu_5+|\nu_{\{67\}}|+2\nu_{8}+|\nu_{\{9a\}}|\\
&&&=&|\lambda_{\{56\}}|+|\mu_{\{7a\}}|+|\nu_{\{58\}}|\\
&&&\leq&0
\end{array}
\end{equation}
The last inequality is the Horn inequality  associated to the LR coefficient
$c(v_{\{56\}}, v_{\{7a\}}, v_{\{58\}})=c_{44,2}^{64}=1$.
It follows that all inequalities in (\ref{equ:ineq}) must be
equalities, which implies that $(\lambda,\mu,\nu)$ belongs to $F$.
\end{proof}

\begin{rema}
We do not know about a general procedure to compute the face given by a triple $I,J,K$ with $c(v_I,v_J,v_K)>1$. When $c(v_I,v_J,v_K)=1$, it is known
that the corresponding face has codimension $1$, but when $c(v_I,v_J,v_K)>1$ it seems an interesting problem to determine the linear span of the
corresponding face, or at least its dimension.
\end{rema}

\begin{rema}
In this case, Lemma~\ref{lem:face2} shows that the face $\cF_{IJK}$  is contained
in some regular face associated to some LR coefficient equal to
one (namely $c(v_{\{56\}}, v_{\{7a\}}, v_{\{58\}}) =1$). Then both Theorems~\ref{introtheo:birational} and
\ref{introtheo:LRalt} can be applied to the points of the face. 
These two statements are not concurrent but complementary. 
Indeed, by applying these results consecutively, one gets an expression
of each LR coefficient on the face of Lemma~\ref{lem:face2} as an alternating
sum of products of three LR coefficients. 
The use of a two step flag variety allows to recover this result from
Corollary~\ref{cor:alt} more conceptually. 
Indeed, set $I_1=\{5,6\}$, $J_1=\{7,10\}$ and $K_1=\{5,8\}$.
Set also $I_2=I$, $J_2=J$, and $K_2=K$.
The pairs $(I_1 \subset I_2)$, $(J_1 \subset J_2)$ and $(K_1 \subset K_2)$ define three Schubert classes in $\Fl(2,6;10)$.

Consider the fibration $\Fl(2,6;10)\longto \Gr(6,10)$ with fiber
$\Gr(2,6)$.
The three given Schubert varieties in  $\Fl(2,6;10)$ map onto
$X_{I_2}$, $X_{J_2}$ and $X_{K_2}$ respectively with fibers isomorphic
to
$X_{\{34\}}$, $X_{\{46\}}$ and $X_{\{35\}}$ in $\Gr(2,6)$.
The associated Schubert coefficient in $\Gr(6,10)$ is $2$ and in
$\Gr(2,6)$ it is $c_{22,1}^{32}=1$. Then (see e.g. \cite{richmond}),
the associated Schubert coefficient for $\Fl(2,6;10)$ is $2$.

Similarly, consider  the fibration $\Fl(2,6;10)\longto \Gr(2,10)$ with fiber
$\Gr(4,8)$.
Here we find a Schubert coefficient 1 in the base and 2 in the
fiber. 

Recalling the BK coefficients $c^\odot$ from \eqref{eq:14}, these two assertions imply that 
\begin{equation}
\label{eq:bkFl}
c^\odot \big ( (I_1 \subset I_2),(J_1 \subset J_2),(K_1 \subset K_2) \big )=2
\end{equation}
and we are in position to apply Corollary~\ref{cor:alt}.
Since the incidence variety $\tilde Y$ in $\Fl(2,6;10)\times (G/B)^3$
maps birationally on that corresponding to $(I_2,J_2,K_2)$ in
$\Gr(6,10)$ these two varieties have the same branch divisor class given
by \eqref{eq:expleBPi2}:
$$
\frac 1 2[B_\Pi]=\cL(3^22^4,2^41^3,3^22^31^3).
$$
Let $(\lambda,\mu,\nu)\in(\Lambda^+_{\GL_{10}})^3$ on the span of the
considered face (see Lemma~\ref{lem:face2}):

\begin{equation}
\begin{array}{l}
  \lambda=\lambda_1^2\lambda_3^2\lambda_5\lambda_6\lambda_7^4\\
\mu=\mu_1^4\mu_5^2\mu_7\mu_8^2\mu_{10}\\
\nu=\nu_1^2\nu_3^2\nu_5\nu_6^2\nu_8\nu_9^2
\end{array}\label{eq:lmnonface}
\end{equation}
The restriction of $\cL(\lambda)\otimes\cL(\mu)\otimes\cL(\nu)$ to $C$
is 
$$
\begin{array}{lll}
 \lambda_{\{56\}}=\lambda_5\lambda_6&\lambda_{\{249a\}}=\lambda_1\lambda_3\lambda_7^2&\lambda_{\{1378\}}=\lambda_1\lambda_3\lambda_7^2\\
\mu_{\{7a\}}=\mu_7\mu_{10}&\mu_{\{3469\}}=\mu_1^2\mu_5\mu_8&\mu_{\{1258\}}=\mu_1^2\mu_5\mu_8\\
\nu_{\{58\}}=\nu_5\nu_8&\nu_{\{247a\}}=\nu_1\nu_3\nu_6\nu_9&\nu_{\{1369\}}=\nu_1\nu_3\nu_6\nu_9.
\end{array}
$$
Similarly, the restriction of $\frac 1 2[B_\Pi]$ to $C$ is as an
$(\SL_2\times\SL_4\times\SL_4)$-linearized line bundle
$$
\begin{array}{lll}
2\ 2&3\ 2\ 0\ 0\ &3\ 2\ 0\ 0\\
1\ 0&2\ 2\ 1\ 0&2\ 2\ 1\ 0\\
2\ 1&3\ 2\ 1\ 0&3\ 2\ 1\ 0.
\end{array}
$$
It has the same invariant sections as the $(\GL_2\times\GL_4\times\GL_4)$-linearized line bundle
$$
\begin{array}{lll}
  0\ 0&\ 3\ 2\ 0\ 0&\ 3\ 2\ 0\ 0\\
1\ 0&\ 2\ 2\ 1\ 0&\ 2\ 2\ 1\ 0\\
0\mbox{-}1&\mbox{-}1\mbox{-}2\mbox{-}3\mbox{-}4&\mbox{-}1\mbox{-}2\mbox{-}3\mbox{-}4.
\end{array}
$$
Hence Corollary~\ref{cor:alt} gives:
$$
  \begin{array}{r@{}c@{}l}
  m_{\GL_{10}}(\lambda,\mu,\nu)=\sum_k &(-1)^k\,&
m_{\GL_2}(\lambda_5\lambda_6,\mu_7\moins k\,\mu_{10},\nu_5\,\nu_8\plus
                                                  k)\\
&\times&
m_{\GL_4}(\lambda_1\moins 3k\,\lambda_3\moins 2k\,\lambda_7^2,
(\mu_1\moins 2k)^2\,\mu_5\moins k\,\mu_8,
\nu_1\plus k\,\nu_3\plus 2k\,\nu_6\plus 3k\,\nu_9\plus 4k)^2\,.
  \end{array}
$$

For $\lambda=[41, 41, 36, 36, 35, 24, 0, 0, 0, 0]$, $\mu=[-41, -41,
-41, -41, -48, -48, -49, -65, -65, -72]$ and $\nu=[49, 49, 42, 42, 40,
25, 25, 22, 2, 2]$, this formula gives 
$$
6= 1\times 3^2-1\times 2^2+1\times 1^2.
$$

Given a point as in \eqref{eq:lmnonface}, one can also  apply Theorem~\ref{introtheo:birational} to the triple
  $(I_1,J_1,K_1)$ to get 
$$
m_{\GL_{10}}(\lambda,\mu,\nu)=m_{\GL_2}(\lambda_5\lambda_6,\mu_7\mu_{10},\nu_5\nu_8)
m_{\GL_8}(\lambda_1^2\lambda_3^2\lambda_7^4,\mu_1^4\mu_5^2\mu_8^2,\nu_1^2\nu_3^2\nu_6^2\nu_9^2).
$$
But, for $\tilde I=\{2,4,7,8\}$, $\tilde J=\{3,4,6,8\}$ and $\tilde
K=\{2,4,6,8\}$ we have $c(\tilde I,\tilde J,\tilde K)=c_{32,221}^{4321}=2$ and one can
apply Theorem~\ref{introtheo:LRalt} to the second factor. 
Since 
  $\frac 12[B_\Pi]=(\varpi_2+\varpi_4,
  \varpi_4+\varpi_6,\varpi_2+\varpi_6)=(2^21^2,2^41^2,2^21^4)$, one
  gets that $m_{\GL_{10}}(\lambda,\mu,\nu)$ is equal to
$$
m_{\GL_2}(\lambda_5\lambda_6,\mu_7\mu_{10},\nu_5\nu_8)
\sum_{k\geq 0}(-1)^k
m_{\GL_4}(\lambda_1\moins 2k\,\lambda_3\moins k\,\lambda_7^2,
\mu_1^2\plus k\,\mu_5\plus 2k\,\mu_8\plus 3k,
\nu_1\moins2k\,\nu_3\moins k\,\nu_6\moins k\,\nu_9)^2.
$$
Applied to the explicit example above, one gets
$$
6=1\times(3^2-2^2+1^2).
$$
\end{rema}

\subsection{\texorpdfstring{Examples in type $B$}{expleB}}
\label{expleB}

In this section, we give some examples with $G=\Spin(2n+1)$. New
phenomena occur compared to the type $A$. First, the Grassmannians are not always
cominuscule and these examples show that the Levi-movability assumption
is needed in Proposition~\ref{pro:posC} and Corollary~\ref{cor:alt}.
Second, the Horn semigroup is not saturated.

Many results in this section were obtained with the help of a computer, see \cite{programmation}.
Given a weight
$\zeta$, we associate a partition $\lambda$ to it with parts $\lambda_i$ such that $\zeta = \sum \lambda_i \epsilon_i$,
with the notation of \cite{bourb}, (note that, given $\zeta$
all the coefficients $\lambda_i$ are in $\Z$ or they are all in $\frac12 + \Z$).

Number the simple roots as follows:
\pgfkeys{/Dynkin diagram,
edge length=8mm,
root radius=1mm,
mark=o}
\dynkin[labels={1,2,n-2,n-1,n}] B{}.
Note that the Horn inequality associated to a triple $(I,J,K)$ of subsets of
$\{1,\ldots,n\} \cup \{n+2,\ldots,2n+1\}$ such that $c(I,J,K) \neq 0$ is
\begin{equation}
\label{eq:HornBn}
|\lambda_{I\cap [1;n]}| - |\lambda_{\overline{I\cap [n+2;2n+1]}}|
+|\mu_{J\cap [1;n]}| - |\mu_{\overline{J\cap [n+2;2n+1]}}|
+|\nu_{K\cap [1;n]}| - |\nu_{\overline{K\cap [n+2;2n+1]}}|\leq 0.
\end{equation}

Let $n=4$. In the Chow ring of the quadric $Q^7$, for $I=J=\{7\}$ and $K=\{6\}$, we have
$\sigma_{\{7\}}\odot\sigma_{\{7\}}\odot\sigma_{\{6\}}=2[pt]$ and \eqref{eq:HornBn} means
\begin{equation}
\label{eq:HornEx1}
\lambda_3+\mu_3+\nu_4 \geq 0.
\end{equation}

Using \cite{programmation}, we get
$\theta=(\varpi_2,\varpi_2, \varpi_3)=(1^2,1^2,1^3)$ and \eqref{eq:HornEx1} is an equality, as predicted by
Proposition \ref{pro:posC}.
Here $L^{\mathrm{ss}}=\Spin(7)$. We put superscript $\varpi^G$ and $\varpi^L$ to make the difference between weights of $G$ and $L$.
We have $v_I\invv(\varpi_2^G)=v_I\invv(\epsilon_1^G+\epsilon_2^G)=\epsilon_2^G+\epsilon_3^G$,
and its restriction to $L$ is $\epsilon_1^L+\epsilon_2^L$, namely $\varpi_2^L$. Similarly,
$(v_K\invv(\varpi_3^G))^L=2\varpi_3^L$. Using the equality $m_{L}(\varpi_2, \varpi_2, 2\varpi_3)=1$
and the following Remark~\ref{rema:dense-orbit}, we get:
$$m_{L}(kv_I\inv \theta^1,kv_J\invv\theta^2,kv_K\invv\theta^3)
=m_{L}(k\varpi_2, k\varpi_2, 2k\varpi_3)=1.$$
Corollary~\ref{cor:alt} gives $m_{G}(k\varpi_2,k\varpi_2, k\varpi_3)=1$ if $k$ is even and $0$ if $k$ is odd,
consistently with the non-saturation of the tensor semigroup
$\Gamma(G):=\{(\zeta_1,\zeta_2,\zeta_3) \in (\Lambda_G^+)^3 \,:\,
m_{G}(\zeta_1,\zeta_2,\zeta_3) \neq 0\}$.

\begin{rema}
\label{rema:dense-orbit}
The group $\Spin_7$ has a dense orbit in $\OG(2,7) \times \OG(2,7) \times \OG(3,7)$.
\end{rema}
\begin{proof}
We consider the isotropic space $\C^2_a$ generated by $e_1$ and $e_2$, the isotropic space $\C^2_b$ generated by
$e_6$ and $e_7$, and the non-isotropic space
$\C^3$ generated by $e_3,e_4,e_5$.
Let $V_3$ be the isotropic space generated by $v_1=e_1+e_3,v_2=e_5-e_7$ and $v_3=e_2+e_4-e_6$.
We claim that the isotropy of the triple $(\C^2_a,\C^2_b,V_3)$
in $\Spin_7$ is $1$-dimensional, which proves the statement since
$\dim(\OG(2,7) \times \OG(2,7) \times \OG(3,7))=20=\dim(\Spin_7)-1$.

Let indeed $g$ be in the isotropy. Since $g$ stabilizes $\C^2_a$ and $\C^2_b$,
it stabilizes $\C^3$ as it is the subspace orthogonal to $\C^2_a \oplus \C^2_b$. The action of $g$ on
$\C^2_b \simeq (\C^2_a)^*$ is the dual of that of $g$ on $\C^2_a$.

Observe that $g$ must stabilize the line generated by $v_1$, since it is the intersection of $V_3$ and $\C^2_a \oplus \C^3$.
Similarly it stabilizes the line generated by $v_2$. It follows that there exists a scalar $\lambda$ such that
$g \cdot e_1 = \lambda e_1, g \cdot e_3=\lambda e_3, g \cdot e_7 = \lambda\inv e_7$ and $g \cdot e_5 = \lambda \inv e_5$.
Then $g$ preserves $\C^2_a \cap e_7^\bot$, which is the line generated by $e_2$. It follows that there exists a scalar $\mu$
such that $g \cdot e_2 = \mu e_2, g \cdot e_6 = \mu \inv e_6$ and $g \cdot e_4=e_4$. Since $g$ preserves the line generated by $v_3$,
we have $\mu=1$. Thus $g$ is uniquely defined by $\lambda$ and the isotropy has dimension $1$.
\end{proof}

\medskip
In the Chow ring of $\OG(4,9)$, for $I=J=\{3, 6, 8, 9\}$ and
$K=
\{1,3,6,8\}$, we have
$\sigma_I\odot\sigma_J\odot\sigma_{K}=2[pt]$. Inequality \eqref{eq:HornBn} means
\begin{equation}
\label{eq:HornEx2}
\lambda_1+\lambda_2+\lambda_4+\mu_1+\mu_2+\mu_4+\nu_2+\nu_4 \geq \lambda_3+\mu_3+\nu_1+\nu_3.
\end{equation}
Using \cite{programmation}, we get $\theta=(\varpi_3,\varpi_3,\varpi_1+\varpi_3)=(1^3\,0,1^3\,0,2\,1^2\,0)$
and once again, \eqref{eq:HornEx2} is an equality.

Here $L^{\mathrm{ss}}=\SL(4)$. We have
$v_I\invv(\theta^1)=v_J\invv(\theta^2)=v_I\invv(\epsilon_1^G+\epsilon_2^G+\epsilon_3^G)
=-\epsilon_4^G-\epsilon_3^G+\epsilon_1^G$ and this restricts to
$2\epsilon_1^L+\epsilon_2^L=\varpi_1^L+\varpi_2^L$. We also have $v_K\invv(\theta^3)^L=\varpi_1^L+\varpi_2^L+\varpi_3^L$. We get
$m_{L}(kv_I\invv\theta^1,kv_J\invv\theta^2,kv_K\invv\theta^3)=c_{2k\, k,\,2k\,k}^{3k\,2k\,k}=k+1$. 
Corollary~\ref{cor:alt} gives $m_{G}(k\varpi_3,
k\varpi_3, k\varpi_1+\varpi_3)=\lceil \frac k 2 \rceil$.

\medskip
In the Chow ring of $\OG(4,9)$, for $I=\{3, 6, 8, 9\}$ and $J=K=\{2, 4, 7, 9\}$, we have
$\sigma_I\sigma_J\sigma_{K}=2[pt]$ and
$\sigma_I\odot\sigma_J\odot\sigma_{K}=0$. Inequality \eqref{eq:HornBn} means
\begin{equation}
\label{eq:HornEx3}
\lambda_1+\lambda_2+\lambda_4+\mu_1+\mu_3+\nu_1+\nu_3 \geq \lambda_3+\mu_2+\mu_4+\nu_2+\nu_4.
\end{equation}Using \cite{programmation}, we get
$\theta=(\varpi_3,\varpi_2+\varpi_4, \varpi_2+\varpi_4)=(1^3\,0,\frac 1 2(3^2\,1^2), \frac 1 2(3^2\,1^2))$,
and this does not
belong to the face. This shows that the Levi-movability assumption in Proposition \ref{pro:posC} is necessary.
The restriction of $\frac 1 2 [B_\Pi]$ to $C$ is  $(21,321,321)$. Corollary~\ref{cor:alt} does not apply.

\appendix

\section{Appendix: puzzles}

\subsection{Some computations}

In this appendix, we give an example for recent results \cite{kp,ressayre-quiver,richmond}
about the Schubert coefficients when $G/P$ is a partial flag variety.
It is a counterexample to \cite[Theorem 3]{kp}. 
Since  \cite[Theorem 3]{kp} is ``proved'' by constructing a bijection
between hives, we picture the concerned hives.

Even more precisely, we consider the flag variety $\Fl(3,6;9)$ parametrizing flags $(V_3 \subset V_6)$ in $\C^9$ with
$\dim(V_3)=3$ and $\dim(V_6)=6$.
We have $\Fl(3,6;9)=G/P$ and $\Gr(6,9)=\SL_9/Q$ with $P \subset Q \subset G=\SL_9$ the corresponding parabolic subgroups.

Consider the Schubert class $\omega=(\{3,6,9\} \subset
\{2,3,5,6,8,9\})$ for $\Fl(3,6;9)$. The projection of the Schubert
variety $X_\omega$ on $\Gr(6,9)$ is $X_I$ with
$I=\{2,3,5,6,8,9\}$. Its intersection with the fiber $\Gr(3,6)$ is
$X_J$ with $J=\{2,4,6\}$.

\begin{exple}
\label{example-calculs}
We have the equalities
$$c_{\Gr(3,6)}(J,J,J)=2\, , \, c_{\Gr(6,9)}(I,I,I)=3 \, , \, \mbox{and } c_{\Fl(3,6;9)}(\omega,\omega,\omega)=6\, .$$
\end{exple}
\begin{rema}
\label{rema-erreur}
Since $6=3*2$, this example is compatible with \cite[Theorem 7.14]{dw}, \cite[Theorem 1.1]{richmond},
and \cite[Theorem A]{ressayre-quiver}.
However, it shows that \cite[Theorem 3]{kp} cannot be correct, since
the right hand side of \cite[Theorem 3]{kp} is equal to $2*2*2=8$.

Actually, the map in  \cite[Theorem~2]{kp} is not surjective.  The
"inverse map" they construct should  certainly be a left-inverse, but
it is never  verified that it is a right-inverse.

In fact, let us denote by $R$ the puzzle displayed
on the left in Appendix \ref{sec:2puz} (with a Rhombus over the edge drawn with a thick line) and by $T$ the puzzle on the right (with a
Triangle). Performing the map $(D_{23},D_{13},D_{12})$ of \cite[Theorem 2]{kp}
on the six puzzles of Appendix \ref{sec:6puz}, we get the six triples
$(R,R,R)$, $(T,R,R)$, $(T,T,R)$, $(T,T,T)$, $(R,R,T)$, $(R,T,T)$.
So the map in \cite[Theorem 2]{kp} is injective, but its image misses the two
triples $(R,T,R)$ and $(T,R,T)$.

As a conclusion, Theorems 2 and 3 in \cite{kp} are not correct, whereas we have no objection for the
other statements in that paper.
\end{rema}

\newpage 

\subsection{\texorpdfstring{The puzzles for $c_{\Gr(3,6)}(J,J,J)=2$ with $J=\{2,4,6\}$}{subsec:2puz}}
\label{sec:2puz}

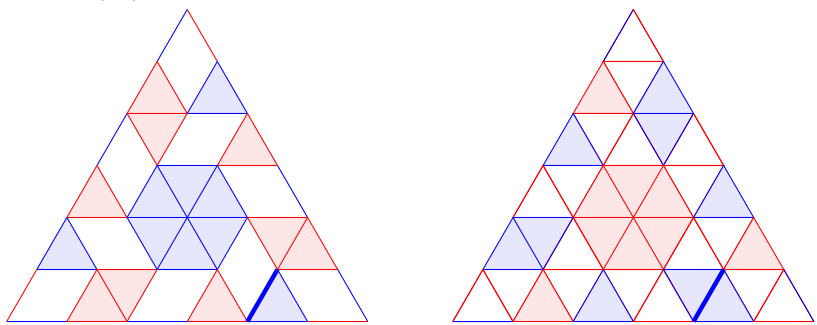
\begin{center}
\begin{tabular}{c@{\hspace{1cm}}c}
   \makebox{
\begin{tikzpicture}[yscale=1.73,scale=0.4]
\path[color=white, fill=white!10](0, 0) -- (-1, -1)-- (1, -1)-- (0, 0);
\path[color=white, fill=white!10](0, -2) -- (-1, -1)-- (1, -1)-- (0, -2);
\draw[color=blue, fill=none] (0, 0) -- (-1, -1);
\draw[color=red, fill=none] (1, -1) -- (0, 0);
\path[] (0, 0)-- (-1, -1)  -- (1, -1) -- (0, 0) ;
\path[] (0, -2)-- (-1, -1) -- (1, -1) -- (0, -2) ;
\path[color=red, fill=red!10](-1, -1) -- (-2, -2)-- (0, -2)-- (-1, -1);
\path[color=red, fill=red!10](-1, -3) -- (-2, -2)-- (0, -2)-- (-1, -3);
\draw[color=red, fill=none] (-2, -2) -- (0, -2);
\draw[color=red, fill=none] (-1, -1) -- (-2, -2);
\draw[color=red, fill=none] (0, -2) -- (-1, -1);
\path[] (-1, -1)-- (-2, -2) -- (0, -2) -- (-1, -1) ;
\path[] (-1, -3)-- (-2, -2) -- (0, -2) -- (-1, -3) ;
\path[color=blue, fill=blue!10](1, -1) -- (0, -2)-- (2, -2)-- (1, -1);
\path[color=white, fill=white!10](1, -3) -- (0, -2)-- (2, -2)-- (1, -3);
\draw[color=blue, fill=none] (0, -2) -- (2, -2);
\draw[color=blue, fill=none] (1, -1) -- (0, -2);
\draw[color=blue, fill=none] (2, -2) -- (1, -1);
\path[] (1, -1)-- (0, -2)  -- (2, -2)  -- (1, -1)  ;
\path[] (1, -3)-- (0, -2) -- (2, -2) -- (1, -3) ;
\path[color=white, fill=white!10](-2, -2) -- (-3, -3)-- (-1, -3)-- (-2, -2);
\path[color=white, fill=white!10](-2, -4) -- (-3, -3)-- (-1, -3)-- (-2, -4);
\draw[color=blue, fill=none] (-2, -2) -- (-3, -3);
\draw[color=red, fill=none] (-1, -3) -- (-2, -2);
\path[] (-2, -2)-- (-3, -3)  -- (-1, -3) -- (-2, -2) ;
\path[] (-2, -4)-- (-3, -3) -- (-1, -3) -- (-2, -4) ;
\path[color=white, fill=white!10](0, -2) -- (-1, -3)-- (1, -3)-- (0, -2);
\path[color=blue, fill=blue!10](0, -4) -- (-1, -3)-- (1, -3)-- (0, -4);
\draw[color=blue, fill=none] (-1, -3) -- (1, -3);
\draw[color=red, fill=none] (0, -2) -- (-1, -3);
\path[] (0, -2)-- (-1, -3) -- (1, -3)  -- (0, -2) ;
\path[] (0, -4)-- (-1, -3) -- (1, -3) -- (0, -4) ;
\path[color=red, fill=red!10](2, -2) -- (1, -3)-- (3, -3)-- (2, -2);
\path[color=white, fill=white!10](2, -4) -- (1, -3)-- (3, -3)-- (2, -4);
\draw[color=red, fill=none] (1, -3) -- (3, -3);
\draw[color=red, fill=none] (2, -2) -- (1, -3);
\draw[color=red, fill=none] (3, -3) -- (2, -2);
\path[] (2, -2)-- (1, -3) -- (3, -3) -- (2, -2) ;
\path[] (2, -4)-- (1, -3) -- (3, -3) -- (2, -4) ;
\path[color=red, fill=red!10](-3, -3) -- (-4, -4)-- (-2, -4)-- (-3, -3);
\path[color=white, fill=white!10](-3, -5) -- (-4, -4)-- (-2, -4)-- (-3, -5);
\draw[color=red, fill=none] (-4, -4) -- (-2, -4);
\draw[color=red, fill=none] (-3, -3) -- (-4, -4);
\draw[color=red, fill=none] (-2, -4) -- (-3, -3);
\path[] (-3, -3)-- (-4, -4) -- (-2, -4) -- (-3, -3) ;
\path[] (-3, -5)-- (-4, -4) -- (-2, -4) -- (-3, -5) ;
\path[color=blue, fill=blue!10](-1, -3) -- (-2, -4)-- (0, -4)-- (-1, -3);
\path[color=blue, fill=blue!10](-1, -5) -- (-2, -4)-- (0, -4)-- (-1, -5);
\draw[color=blue, fill=none] (-2, -4) -- (0, -4);
\draw[color=blue, fill=none] (-1, -3) -- (-2, -4);
\draw[color=blue, fill=none] (0, -4) -- (-1, -3);
\path[] (-1, -3)-- (-2, -4)  -- (0, -4)  -- (-1, -3)  ;
\path[] (-1, -5)-- (-2, -4) -- (0, -4) -- (-1, -5) ;
\path[color=blue, fill=blue!10](1, -3) -- (0, -4)-- (2, -4)-- (1, -3);
\path[color=blue, fill=blue!10](1, -5) -- (0, -4)-- (2, -4)-- (1, -5);
\draw[color=blue, fill=none] (0, -4) -- (2, -4);
\draw[color=blue, fill=none] (1, -3) -- (0, -4);
\draw[color=blue, fill=none] (2, -4) -- (1, -3);
\path[] (1, -3)-- (0, -4)  -- (2, -4)  -- (1, -3)  ;
\path[] (1, -5)-- (0, -4) -- (2, -4) -- (1, -5) ;
\path[color=white, fill=white!10](3, -3) -- (2, -4)-- (4, -4)-- (3, -3);
\path[color=red, fill=red!10](3, -5) -- (2, -4)-- (4, -4)-- (3, -5);
\draw[color=red, fill=none] (2, -4) -- (4, -4);
\draw[color=blue, fill=none] (4, -4) -- (3, -3);
\path[] (3, -3)-- (2, -4) -- (4, -4) -- (3, -3)  ;
\path[] (3, -5)-- (2, -4) -- (4, -4) -- (3, -5) ;
\path[color=blue, fill=blue!10](-4, -4) -- (-5, -5)-- (-3, -5)-- (-4, -4);
\path[color=white, fill=white!10](-4, -6) -- (-5, -5)-- (-3, -5)-- (-4, -6);
\draw[color=blue, fill=none] (-5, -5) -- (-3, -5);
\draw[color=blue, fill=none] (-4, -4) -- (-5, -5);
\draw[color=blue, fill=none] (-3, -5) -- (-4, -4);
\path[] (-4, -4)-- (-5, -5)  -- (-3, -5)  -- (-4, -4)  ;
\path[] (-4, -6)-- (-5, -5) -- (-3, -5) -- (-4, -6) ;
\path[color=white, fill=white!10](-2, -4) -- (-3, -5)-- (-1, -5)-- (-2, -4);
\path[color=red, fill=red!10](-2, -6) -- (-3, -5)-- (-1, -5)-- (-2, -6);
\draw[color=red, fill=none] (-3, -5) -- (-1, -5);
\draw[color=blue, fill=none] (-1, -5) -- (-2, -4);
\path[] (-2, -4)-- (-3, -5) -- (-1, -5) -- (-2, -4)  ;
\path[] (-2, -6)-- (-3, -5) -- (-1, -5) -- (-2, -6) ;
\path[color=blue, fill=blue!10](0, -4) -- (-1, -5)-- (1, -5)-- (0, -4);
\path[color=white, fill=white!10](0, -6) -- (-1, -5)-- (1, -5)-- (0, -6);
\draw[color=blue, fill=none] (-1, -5) -- (1, -5);
\draw[color=blue, fill=none] (0, -4) -- (-1, -5);
\draw[color=blue, fill=none] (1, -5) -- (0, -4);
\path[] (0, -4)-- (-1, -5)  -- (1, -5)  -- (0, -4)  ;
\path[] (0, -6)-- (-1, -5) -- (1, -5) -- (0, -6) ;
\path[color=white, fill=white!10](2, -4) -- (1, -5)-- (3, -5)-- (2, -4);
\path[color=white, fill=white!10](2, -6) -- (1, -5)-- (3, -5)-- (2, -6);
\draw[color=blue, fill=none] (2, -4) -- (1, -5);
\draw[color=red, fill=none] (3, -5) -- (2, -4);
\path[] (2, -4)-- (1, -5)  -- (3, -5) -- (2, -4) ;
\path[] (2, -6)-- (1, -5) -- (3, -5) -- (2, -6) ;
\path[color=red, fill=red!10](4, -4) -- (3, -5)-- (5, -5)-- (4, -4);
\path[color=white, fill=white!10](4, -6) -- (3, -5)-- (5, -5)-- (4, -6);
\draw[color=red, fill=none] (3, -5) -- (5, -5);
\draw[color=red, fill=none] (4, -4) -- (3, -5);
\draw[color=red, fill=none] (5, -5) -- (4, -4);
\path[] (4, -4)-- (3, -5) -- (5, -5) -- (4, -4) ;
\path[] (4, -6)-- (3, -5) -- (5, -5) -- (4, -6) ;
\path[color=white, fill=white!10](-5, -5) -- (-6, -6)-- (-4, -6)-- (-5, -5);
\draw[color=blue, fill=none] (-6, -6) -- (-4, -6);
\draw[color=red, fill=none] (-5, -5) -- (-6, -6);
\path[] (-5, -5)-- (-6, -6) -- (-4, -6)  -- (-5, -5) ;
\path[color=red, fill=red!10](-3, -5) -- (-4, -6)-- (-2, -6)-- (-3, -5);
\draw[color=red, fill=none] (-4, -6) -- (-2, -6);
\draw[color=red, fill=none] (-3, -5) -- (-4, -6);
\draw[color=red, fill=none] (-2, -6) -- (-3, -5);
\path[] (-3, -5)-- (-4, -6) -- (-2, -6) -- (-3, -5) ;
\path[color=white, fill=white!10](-1, -5) -- (-2, -6)-- (0, -6)-- (-1, -5);
\draw[color=blue, fill=none] (-2, -6) -- (0, -6);
\draw[color=red, fill=none] (-1, -5) -- (-2, -6);
\path[] (-1, -5)-- (-2, -6) -- (0, -6)  -- (-1, -5) ;
\path[color=red, fill=red!10](1, -5) -- (0, -6)-- (2, -6)-- (1, -5);
\draw[color=red, fill=none] (0, -6) -- (2, -6);
\draw[color=red, fill=none] (1, -5) -- (0, -6);
\draw[color=red, fill=none] (2, -6) -- (1, -5);
\path[] (1, -5)-- (0, -6) -- (2, -6) -- (1, -5) ;
\path[color=blue, fill=blue!10](3, -5) -- (2, -6)-- (4, -6)-- (3, -5);
\draw[color=blue, fill=none] (2, -6) -- (4, -6);
\draw[color=blue, fill=none] (3, -5) -- (2, -6);
\draw[color=blue, fill=none] (4, -6) -- (3, -5);
\path[] (3, -5)-- (2, -6)  -- (4, -6)
 -- (3, -5) ;

\path[color=white, fill=white!10](5, -5) -- (4, -6)-- (6, -6)-- (5, -5);
\draw[color=red, fill=none] (4, -6) -- (6, -6);
\draw[color=blue, fill=none] (6, -6) -- (5, -5);
\path[] (5, -5)-- (4, -6) -- (6, -6) -- (5, -5)  ;

\draw[color=blue,line width=0.6mm](3,-5) -- (2,-6);
\end{tikzpicture}}
&\makebox{
\begin{tikzpicture}[yscale=1.73,scale=0.4]
\draw[color=red](0, 0) -- (-1, -1)-- (1, -1)-- (0, 0);
\draw[color=red](0, -2) -- (-1, -1)-- (1, -1)-- (0, -2);
\draw[color=blue, fill=none] (0, 0) -- (-1, -1);
\draw[color=red, fill=none] (1, -1) -- (0, 0);
\path[] (0, 0)-- (-1, -1)  -- (1, -1) -- (0, 0) ;
\path[] (0, -2)-- (-1, -1) -- (1, -1) -- (0, -2) ;
\path[color=red, fill=red!10](-1, -1) -- (-2, -2)-- (0, -2)-- (-1, -1);
\draw[color=red](-1, -3) -- (-2, -2)-- (0, -2)-- (-1, -3);
\draw[color=red, fill=none] (-2, -2) -- (0, -2);
\draw[color=red, fill=none] (-1, -1) -- (-2, -2);
\draw[color=red, fill=none] (0, -2) -- (-1, -1);
\path[] (-1, -1)-- (-2, -2) -- (0, -2) -- (-1, -1) ;
\path[] (-1, -3)-- (-2, -2) -- (0, -2) -- (-1, -3) ;
\path[color=blue, fill=blue!10](1, -1) -- (0, -2)-- (2, -2)-- (1, -1);
\path[color=blue, fill=blue!10](1, -3) -- (0, -2)-- (2, -2)-- (1, -3);
\draw[color=blue, fill=none] (0, -2) -- (2, -2);
\draw[color=blue, fill=none] (1, -1) -- (0, -2);
\draw[color=blue, fill=none] (2, -2) -- (1, -1);
\path[] (1, -1)-- (0, -2)  -- (2, -2)  -- (1, -1)  ;
\path[] (1, -3)-- (0, -2) -- (2, -2) -- (1, -3) ;
\path[color=blue, fill=blue!10](-2, -2) -- (-3, -3)-- (-1, -3)-- (-2, -2);
\draw[color=red](-2, -4) -- (-3, -3)-- (-1, -3)-- (-2, -4);
\draw[color=blue, fill=none] (-3, -3) -- (-1, -3);
\draw[color=blue, fill=none] (-2, -2) -- (-3, -3);
\draw[color=blue, fill=none] (-1, -3) -- (-2, -2);
\path[] (-2, -2)-- (-3, -3)  -- (-1, -3)  -- (-2, -2)  ;
\path[] (-2, -4)-- (-3, -3) -- (-1, -3) -- (-2, -4) ;
\draw[color=red](0, -2) -- (-1, -3)-- (1, -3)-- (0, -2);
\path[color=red, fill=red!10](0, -4) -- (-1, -3)-- (1, -3)-- (0, -4);
\draw[color=red, fill=none] (-1, -3) -- (1, -3);
\draw[color=blue, fill=none] (1, -3) -- (0, -2);
\path[] (0, -2)-- (-1, -3) -- (1, -3) -- (0, -2)  ;
\path[] (0, -4)-- (-1, -3) -- (1, -3) -- (0, -4) ;
\draw[color=red](2, -2) -- (1, -3)-- (3, -3)-- (2, -2);
\draw[color=red](2, -4) -- (1, -3)-- (3, -3)-- (2, -4);
\draw[color=blue, fill=none] (2, -2) -- (1, -3);
\draw[color=red, fill=none] (3, -3) -- (2, -2);
\path[] (2, -2)-- (1, -3)  -- (3, -3) -- (2, -2) ;
\path[] (2, -4)-- (1, -3) -- (3, -3) -- (2, -4) ;
\draw[color=red](-3, -3) -- (-4, -4)-- (-2, -4)-- (-3, -3);
\path[color=blue, fill=blue!10](-3, -5) -- (-4, -4)-- (-2, -4)-- (-3, -5);
\draw[color=blue, fill=none] (-4, -4) -- (-2, -4);
\draw[color=red, fill=none] (-3, -3) -- (-4, -4);
\path[] (-3, -3)-- (-4, -4) -- (-2, -4)  -- (-3, -3) ;
\path[] (-3, -5)-- (-4, -4) -- (-2, -4) -- (-3, -5) ;
\path[color=red, fill=red!10](-1, -3) -- (-2, -4)-- (0, -4)-- (-1, -3);
\path[color=red, fill=red!10](-1, -5) -- (-2, -4)-- (0, -4)-- (-1, -5);
\draw[color=red, fill=none] (-2, -4) -- (0, -4);
\draw[color=red, fill=none] (-1, -3) -- (-2, -4);
\draw[color=red, fill=none] (0, -4) -- (-1, -3);
\path[] (-1, -3)-- (-2, -4) -- (0, -4) -- (-1, -3) ;
\path[] (-1, -5)-- (-2, -4) -- (0, -4) -- (-1, -5) ;
\path[color=red, fill=red!10](1, -3) -- (0, -4)-- (2, -4)-- (1, -3);
\path[color=red, fill=red!10](1, -5) -- (0, -4)-- (2, -4)-- (1, -5);
\draw[color=red, fill=none] (0, -4) -- (2, -4);
\draw[color=red, fill=none] (1, -3) -- (0, -4);
\draw[color=red, fill=none] (2, -4) -- (1, -3);
\path[] (1, -3)-- (0, -4) -- (2, -4) -- (1, -3) ;
\path[] (1, -5)-- (0, -4) -- (2, -4) -- (1, -5) ;
\path[color=blue, fill=blue!10](3, -3) -- (2, -4)-- (4, -4)-- (3, -3);
\draw[color=red](3, -5) -- (2, -4)-- (4, -4)-- (3, -5);
\draw[color=blue, fill=none] (2, -4) -- (4, -4);
\draw[color=blue, fill=none] (3, -3) -- (2, -4);
\draw[color=blue, fill=none] (4, -4) -- (3, -3);
\path[] (3, -3)-- (2, -4)  -- (4, -4)  -- (3, -3)  ;
\path[] (3, -5)-- (2, -4) -- (4, -4) -- (3, -5) ;
\path[color=blue, fill=blue!10](-4, -4) -- (-5, -5)-- (-3, -5)-- (-4, -4);
\draw[color=red](-4, -6) -- (-5, -5)-- (-3, -5)-- (-4, -6);
\draw[color=blue, fill=none] (-5, -5) -- (-3, -5);
\draw[color=blue, fill=none] (-4, -4) -- (-5, -5);
\draw[color=blue, fill=none] (-3, -5) -- (-4, -4);
\path[] (-4, -4)-- (-5, -5)  -- (-3, -5)  -- (-4, -4)  ;
\path[] (-4, -6)-- (-5, -5) -- (-3, -5) -- (-4, -6) ;
\draw[color=red](-2, -4) -- (-3, -5)-- (-1, -5)-- (-2, -4);
\draw[color=red](-2, -6) -- (-3, -5)-- (-1, -5)-- (-2, -6);
\draw[color=blue, fill=none] (-2, -4) -- (-3, -5);
\draw[color=red, fill=none] (-1, -5) -- (-2, -4);
\path[] (-2, -4)-- (-3, -5)  -- (-1, -5) -- (-2, -4) ;
\path[] (-2, -6)-- (-3, -5) -- (-1, -5) -- (-2, -6) ;
\path[color=red, fill=red!10](0, -4) -- (-1, -5)-- (1, -5)-- (0, -4);
\draw[color=red](0, -6) -- (-1, -5)-- (1, -5)-- (0, -6);
\draw[color=red, fill=none] (-1, -5) -- (1, -5);
\draw[color=red, fill=none] (0, -4) -- (-1, -5);
\draw[color=red, fill=none] (1, -5) -- (0, -4);
\path[] (0, -4)-- (-1, -5) -- (1, -5) -- (0, -4) ;
\path[] (0, -6)-- (-1, -5) -- (1, -5) -- (0, -6) ;
\draw[color=red](2, -4) -- (1, -5)-- (3, -5)-- (2, -4);
\path[color=blue, fill=blue!10](2, -6) -- (1, -5)-- (3, -5)-- (2, -6);
\draw[color=blue, fill=none] (1, -5) -- (3, -5);
\draw[color=red, fill=none] (2, -4) -- (1, -5);
\path[] (2, -4)-- (1, -5) -- (3, -5)  -- (2, -4) ;
\path[] (2, -6)-- (1, -5) -- (3, -5) -- (2, -6) ;
\path[color=red, fill=red!10](4, -4) -- (3, -5)-- (5, -5)-- (4, -4);
\draw[color=red](4, -6) -- (3, -5)-- (5, -5)-- (4, -6);
\draw[color=red, fill=none] (3, -5) -- (5, -5);
\draw[color=red, fill=none] (4, -4) -- (3, -5);
\draw[color=red, fill=none] (5, -5) -- (4, -4);
\path[] (4, -4)-- (3, -5) -- (5, -5) -- (4, -4) ;
\path[] (4, -6)-- (3, -5) -- (5, -5) -- (4, -6) ;
\draw[color=red](-5, -5) -- (-6, -6)-- (-4, -6)-- (-5, -5);
\draw[color=blue, fill=none] (-6, -6) -- (-4, -6);
\draw[color=red, fill=none] (-5, -5) -- (-6, -6);
\path[] (-5, -5)-- (-6, -6) -- (-4, -6)  -- (-5, -5) ;
\path[color=red, fill=red!10](-3, -5) -- (-4, -6)-- (-2, -6)-- (-3, -5);
\draw[color=red, fill=none] (-4, -6) -- (-2, -6);
\draw[color=red, fill=none] (-3, -5) -- (-4, -6);
\draw[color=red, fill=none] (-2, -6) -- (-3, -5);
\path[] (-3, -5)-- (-4, -6) -- (-2, -6) -- (-3, -5) ;
\path[color=blue, fill=blue!10](-1, -5) -- (-2, -6)-- (0, -6)-- (-1, -5);
\draw[color=blue, fill=none] (-2, -6) -- (0, -6);
\draw[color=blue, fill=none] (-1, -5) -- (-2, -6);
\draw[color=blue, fill=none] (0, -6) -- (-1, -5);
\path[] (-1, -5)-- (-2, -6)  -- (0, -6)  -- (-1, -5)  ;
\draw[color=red](1, -5) -- (0, -6)-- (2, -6)-- (1, -5);
\draw[color=red, fill=none] (0, -6) -- (2, -6);
\draw[color=blue, fill=none] (2, -6) -- (1, -5);
\path[] (1, -5)-- (0, -6) -- (2, -6) -- (1, -5)  ;
\path[color=blue, fill=blue!10](3, -5) -- (2, -6)-- (4, -6)-- (3, -5);
\draw[color=blue, fill=none] (2, -6) -- (4, -6);
\draw[color=blue, fill=none] (3, -5) -- (2, -6);
\draw[color=blue, fill=none] (4, -6) -- (3, -5);
\path[] (3, -5)-- (2, -6)  -- (4, -6)  -- (3, -5)  ;
\draw[color=red](5, -5) -- (4, -6)-- (6, -6)-- (5, -5);
\draw[color=red, fill=none] (4, -6) -- (6, -6);
\draw[color=blue, fill=none] (6, -6) -- (5, -5);
\path[] (5, -5)-- (4, -6) -- (6, -6) --
(5, -5)  ;

\draw[color=blue,line width=0.6mm](3,-5) -- (2,-6);
\end{tikzpicture}}
\end{tabular}
\end{center}

\subsection{\texorpdfstring{The puzzles for $c_{\Gr(6,9)}(I,I,I)=3$ with $I=\{2,3,5,6,8,9\}$}{subsec:3puz}}
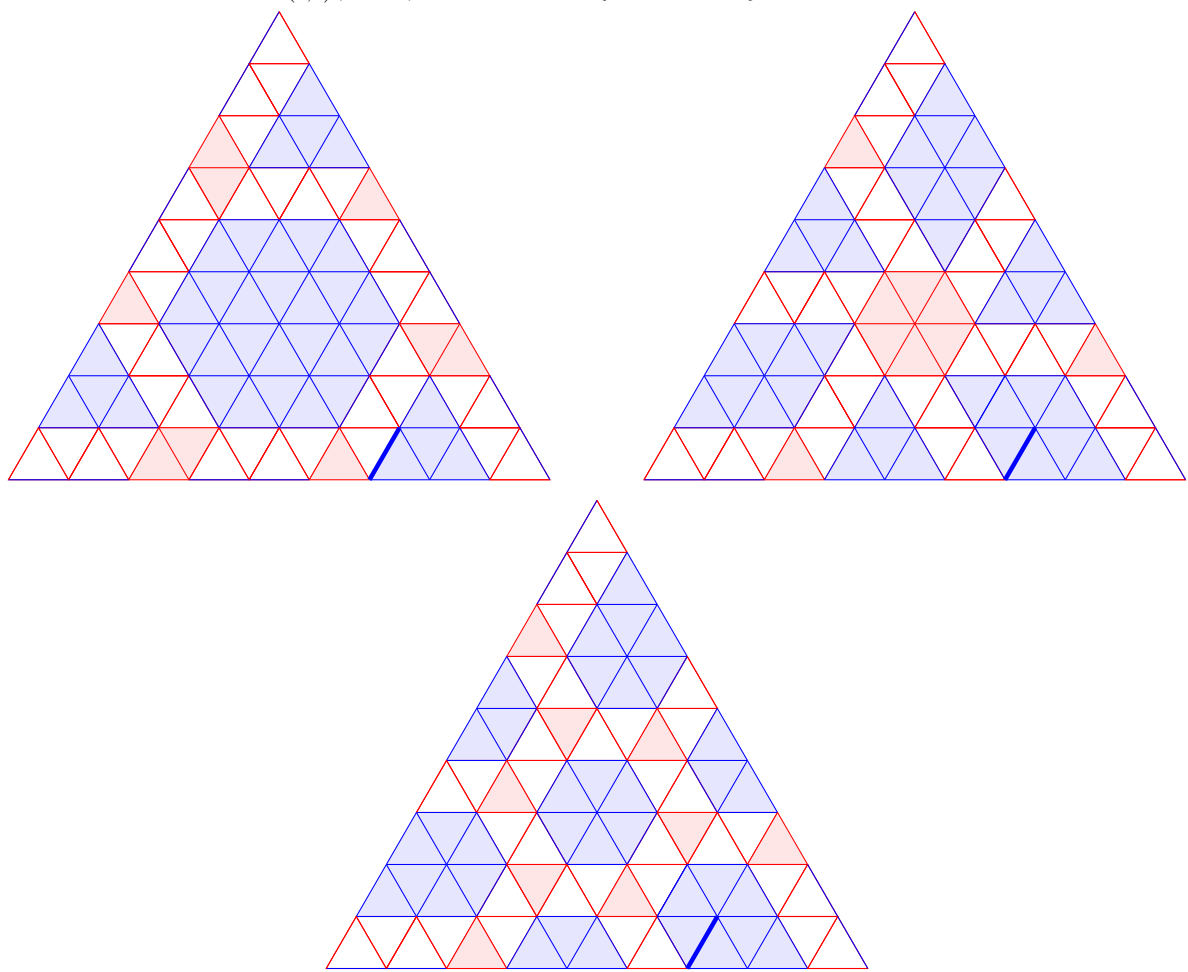
\label{sec:3puz}

\input{g69}

\subsection{\texorpdfstring{The puzzles for $c_{\cF(3,6;9)}(\omega,\omega,\omega)=6$ with
$\omega=(\{3,6,9\} \subset \{2,3,5,6,8,9\})$}{subsec:6puz}}
\label{sec:6puz}

\newpage

\input{f369}

\vspace{10mm}

\end{document}